\documentclass{elsarticle}
  
\usepackage[margin=1in]{geometry}

\usepackage{bbm}
\usepackage{graphicx}
\usepackage{pstool}
\usepackage{wrapfig}
\usepackage{caption}
\usepackage{subcaption}
\usepackage[section]{placeins} 

\usepackage{fancyhdr} 
\usepackage{lastpage} 
\usepackage{extramarks} 

\usepackage{xcolor} 
\usepackage{enumerate} 
\usepackage{paralist} 
\usepackage{amsmath, amsthm, amssymb, mathtools}
\usepackage{mathabx, pifont, stmaryrd} 
\usepackage[explicit]{titlesec} 
\usepackage{etoolbox} 
\usepackage{bibentry} 
\makeatletter\let\saved@bibitem\@bibitem\makeatother 
\usepackage[colorlinks, bookmarksopen, bookmarksnumbered,
            citecolor=red,urlcolor=red]{hyperref} 
\makeatletter\let\@bibitem\saved@bibitem\makeatother 
\usepackage{rotating}

\usepackage{algorithm}
\usepackage{algorithmic}

\newtheorem{remark}{Remark}

\theoremstyle{definition}

\usepackage{bm} 
\usepackage{amsfonts} 

\DeclareMathOperator*{\argmin}{arg\,min}

\newcommand{\func}[3]{\ensuremath{#1 : #2 \rightarrow #3}}

\newcommand{\norm}[1]{\ensuremath{\left\| #1 \right\|}}

\newcommand{\suchthat}{\mathrel{}\middle|\mathrel{}}

\newcommand{\pder}[2]{\ensuremath{\frac{\partial #1}{\partial #2}}}

\newcommand{\Dcal}{\ensuremath{\mathcal{D}}}
\newcommand{\Ecal}{\ensuremath{\mathcal{E}}}

\newcommand{\Gcal}{\ensuremath{\mathcal{G}}}
\newcommand{\Hcal}{\ensuremath{\mathcal{H}}}

\newcommand{\Ocal}{\ensuremath{\mathcal{O}}}
\newcommand{\Pcal}{\ensuremath{\mathcal{P}}}

\newcommand{\Tcal}{\ensuremath{\mathcal{T}}}

\newcommand{\Vcal}{\ensuremath{\mathcal{V}}}
\newcommand{\Wcal}{\ensuremath{\mathcal{W}}}

\newcommand{\Gbb}{\ensuremath{\mathbb{G}}}

\newcommand{\Rbb}{\ensuremath{\mathbb{R} }}

\newcommand\Rbm{{\ensuremath{\bm{R}}}}

\newcommand\rbm{{\ensuremath{\bm{r}}}}

\newcommand\ubm{{\ensuremath{\bm{u}}}}

\newcommand\xbm{{\ensuremath{\bm{x}}}}
\newcommand\ybm{{\ensuremath{\bm{y}}}}
\newcommand\zbm{{\ensuremath{\bm{z}}}}

\newcommand\mubold{{\ensuremath{\boldsymbol{\mu}}}}

\newcommand\phibold{{\ensuremath{\boldsymbol{\phi}}}}

\newcommand\zerobold{\ensuremath{\mathbf{0}}}

\usepackage{tikz}
\usepackage{pgfplots}
\usepackage{pgfplotstable, filecontents, booktabs}
\pgfplotsset{compat=1.9}

\usetikzlibrary{pgfplots.groupplots}
\usepgfplotslibrary{fillbetween}
\usetikzlibrary{calc,fit,matrix,arrows,automata,positioning,shapes}
\usetikzlibrary{arrows.meta}

\pgfplotsset{select coords between index/.style 2 args={
    x filter/.code={
        \ifnum\coordindex<#1\fi
        \ifnum\coordindex>#2\fi
    }
}}

\tikzset{
 invisible/.style={opacity=0},
 visible on/.style={alt={#1{}{invisible}}},
 alt/.code args={<#1>#2#3}{%
   \alt<#1>{\pgfkeysalso{#2}}{\pgfkeysalso{#3}}
 },
}

\newcommand{\colorbarMatlabParula}[5]{
\begin{tikzpicture}
\begin{axis}[
   hide axis, scale only axis,
   height=0pt, width=0pt,
   colormap={parula}{rgb255=(62,38,168) rgb255=(62,39,172) rgb255=(63,40,175) rgb255=(63,41,178) rgb255=(64,42,180) rgb255=(64,43,183) rgb255=(65,44,186) rgb255=(65,45,189) rgb255=(66,46,191) rgb255=(66,47,194) rgb255=(67,48,197) rgb255=(67,49,200) rgb255=(67,50,202) rgb255=(68,51,205) rgb255=(68,52,208) rgb255=(69,53,210) rgb255=(69,55,213) rgb255=(69,56,215) rgb255=(70,57,217) rgb255=(70,58,220) rgb255=(70,59,222) rgb255=(70,61,224) rgb255=(71,62,225) rgb255=(71,63,227) rgb255=(71,65,229) rgb255=(71,66,230) rgb255=(71,68,232) rgb255=(71,69,233) rgb255=(71,70,235) rgb255=(72,72,236) rgb255=(72,73,237) rgb255=(72,75,238) rgb255=(72,76,240) rgb255=(72,78,241) rgb255=(72,79,242) rgb255=(72,80,243) rgb255=(72,82,244) rgb255=(72,83,245) rgb255=(72,84,246) rgb255=(71,86,247) rgb255=(71,87,247) rgb255=(71,89,248) rgb255=(71,90,249) rgb255=(71,91,250) rgb255=(71,93,250) rgb255=(70,94,251) rgb255=(70,96,251) rgb255=(70,97,252) rgb255=(69,98,252) rgb255=(69,100,253) rgb255=(68,101,253) rgb255=(67,103,253) rgb255=(67,104,254) rgb255=(66,106,254) rgb255=(65,107,254) rgb255=(64,109,254) rgb255=(63,110,255) rgb255=(62,112,255) rgb255=(60,113,255) rgb255=(59,115,255) rgb255=(57,116,255) rgb255=(56,118,254) rgb255=(54,119,254) rgb255=(53,121,253) rgb255=(51,122,253) rgb255=(50,124,252) rgb255=(49,125,252) rgb255=(48,127,251) rgb255=(47,128,250) rgb255=(47,130,250) rgb255=(46,131,249) rgb255=(46,132,248) rgb255=(46,134,248) rgb255=(46,135,247) rgb255=(45,136,246) rgb255=(45,138,245) rgb255=(45,139,244) rgb255=(45,140,243) rgb255=(45,142,242) rgb255=(44,143,241) rgb255=(44,144,240) rgb255=(43,145,239) rgb255=(42,147,238) rgb255=(41,148,237) rgb255=(40,149,236) rgb255=(39,151,235) rgb255=(39,152,234) rgb255=(38,153,233) rgb255=(38,154,232) rgb255=(37,155,232) rgb255=(37,156,231) rgb255=(36,158,230) rgb255=(36,159,229) rgb255=(35,160,229) rgb255=(35,161,228) rgb255=(34,162,228) rgb255=(33,163,227) rgb255=(32,165,227) rgb255=(31,166,226) rgb255=(30,167,225) rgb255=(29,168,225) rgb255=(29,169,224) rgb255=(28,170,223) rgb255=(27,171,222) rgb255=(26,172,221) rgb255=(25,173,220) rgb255=(23,174,218) rgb255=(22,175,217) rgb255=(20,176,216) rgb255=(18,177,214) rgb255=(16,178,213) rgb255=(14,179,212) rgb255=(11,179,210) rgb255=(8,180,209) rgb255=(6,181,207) rgb255=(4,182,206) rgb255=(2,183,204) rgb255=(1,183,202) rgb255=(0,184,201) rgb255=(0,185,199) rgb255=(0,186,198) rgb255=(1,186,196) rgb255=(2,187,194) rgb255=(4,187,193) rgb255=(6,188,191) rgb255=(9,189,189) rgb255=(13,189,188) rgb255=(16,190,186) rgb255=(20,190,184) rgb255=(23,191,182) rgb255=(26,192,181) rgb255=(29,192,179) rgb255=(32,193,177) rgb255=(35,193,175) rgb255=(37,194,174) rgb255=(39,194,172) rgb255=(41,195,170) rgb255=(43,195,168) rgb255=(44,196,166) rgb255=(46,196,165) rgb255=(47,197,163) rgb255=(49,197,161) rgb255=(50,198,159) rgb255=(51,199,157) rgb255=(53,199,155) rgb255=(54,200,153) rgb255=(56,200,150) rgb255=(57,201,148) rgb255=(59,201,146) rgb255=(61,202,144) rgb255=(64,202,141) rgb255=(66,202,139) rgb255=(69,203,137) rgb255=(72,203,134) rgb255=(75,203,132) rgb255=(78,204,129) rgb255=(81,204,127) rgb255=(84,204,124) rgb255=(87,204,122) rgb255=(90,204,119) rgb255=(94,205,116) rgb255=(97,205,114) rgb255=(100,205,111) rgb255=(103,205,108) rgb255=(107,205,105) rgb255=(110,205,102) rgb255=(114,205,100) rgb255=(118,204,97) rgb255=(121,204,94) rgb255=(125,204,91) rgb255=(129,204,89) rgb255=(132,204,86) rgb255=(136,203,83) rgb255=(139,203,81) rgb255=(143,203,78) rgb255=(147,202,75) rgb255=(150,202,72) rgb255=(154,201,70) rgb255=(157,201,67) rgb255=(161,200,64) rgb255=(164,200,62) rgb255=(167,199,59) rgb255=(171,199,57) rgb255=(174,198,55) rgb255=(178,198,53) rgb255=(181,197,51) rgb255=(184,196,49) rgb255=(187,196,47) rgb255=(190,195,45) rgb255=(194,195,44) rgb255=(197,194,42) rgb255=(200,193,41) rgb255=(203,193,40) rgb255=(206,192,39) rgb255=(208,191,39) rgb255=(211,191,39) rgb255=(214,190,39) rgb255=(217,190,40) rgb255=(219,189,40) rgb255=(222,188,41) rgb255=(225,188,42) rgb255=(227,188,43) rgb255=(230,187,45) rgb255=(232,187,46) rgb255=(234,186,48) rgb255=(236,186,50) rgb255=(239,186,53) rgb255=(241,186,55) rgb255=(243,186,57) rgb255=(245,186,59) rgb255=(247,186,61) rgb255=(249,186,62) rgb255=(251,187,62) rgb255=(252,188,62) rgb255=(254,189,61) rgb255=(254,190,60) rgb255=(254,192,59) rgb255=(254,193,58) rgb255=(254,194,57) rgb255=(254,196,56) rgb255=(254,197,55) rgb255=(254,199,53) rgb255=(254,200,52) rgb255=(254,202,51) rgb255=(253,203,50) rgb255=(253,205,49) rgb255=(253,206,49) rgb255=(252,208,48) rgb255=(251,210,47) rgb255=(251,211,46) rgb255=(250,213,46) rgb255=(249,214,45) rgb255=(249,216,44) rgb255=(248,217,43) rgb255=(247,219,42) rgb255=(247,221,42) rgb255=(246,222,41) rgb255=(246,224,40) rgb255=(245,225,40) rgb255=(245,227,39) rgb255=(245,229,38) rgb255=(245,230,38) rgb255=(245,232,37) rgb255=(245,233,36) rgb255=(245,235,35) rgb255=(245,236,34) rgb255=(245,238,33) rgb255=(246,239,32) rgb255=(246,241,31) rgb255=(246,242,30) rgb255=(247,244,28) rgb255=(247,245,27) rgb255=(248,247,26) rgb255=(248,248,24) rgb255=(249,249,22) rgb255=(249,251,21) },
   colorbar horizontal,
   point meta min=#1, point meta max=#5,
   colorbar style={width=10cm, xtick={#1,#2,#3,#4,#5}}
]
\addplot [draw=none] coordinates {(0,0)};
\end{axis}
\end{tikzpicture}
}

\usepackage{subcaption}
\usepackage{pdflscape}

\usepackage{setspace}

\newbool{fastcompile}
\setbool{fastcompile}{false}

\begin{document}
\title{High-order implicit shock tracking boundary conditions for flows with parametrized shocks}

\author[rvt1]{Tianci Huang\fnref{fn1}}
\ead{thuang5@nd.edu}

\author[rvt1]{Charles Naudet\fnref{fn2}}
\ead{cnaudet@nd.edu}

\author[rvt1]{Matthew J. Zahr\fnref{fn3}\corref{cor1}}
\ead{mzahr@nd.edu}

\address[rvt1]{Department of Aerospace and Mechanical Engineering, University
               of Notre Dame, Notre Dame, IN 46556, United States}
\cortext[cor1]{Corresponding author}

\fntext[fn1]{Graduate Student, Department of Aerospace and Mechanical
             Engineering, University of Notre Dame}
\fntext[fn2]{Graduate Student, Department of Aerospace and Mechanical
             Engineering, University of Notre Dame}
\fntext[fn3]{Assistant Professor, Department of Aerospace and Mechanical
             Engineering, University of Notre Dame}


\begin{abstract}
High-order implicit shock tracking (fitting) is a class of high-order,
optimization-based numerical methods to approximate solutions of conservation
laws with non-smooth features by aligning elements of the computational mesh
with non-smooth features. This ensures the non-smooth features are perfectly
represented by inter-element jumps and high-order basis functions approximate
smooth regions of the solution without nonlinear stabilization, which leads
to accurate approximations on traditionally coarse meshes. In this work, we
introduce a robust implicit shock tracking framework specialized
for problems with parameter-dependent lead shocks (i.e., shocks separating a
farfield condition from the downstream flow), which commonly arise in high-speed
aerodynamics and astrophysics applications. After a shock-aligned mesh is
produced at one parameter configuration, all elements upstream of the
lead shock are removed and the nodes on the lead shock are positioned
for new parameter configurations using the implicit shock tracking solver.
The proposed framework can be used for most many-query applications involving
parametrized lead shocks such as optimization, uncertainty quantification,
parameter sweeps, ``what-if'' scenarios, or parameter-based continuation.
We demonstrate the robustness and flexibility of the framework
using a one-dimensional space-time Riemann problem, and two- and 
three-dimensional supersonic and hypersonic benchmark problems.
\end{abstract}


\maketitle

\section{Introduction}
\label{sec:intro}
Bow shocks are strong, detached, curved shocks that frequently arise in
aerospace and astrophysics applications. The strength, geometry, and
stand-off distance of bow shocks are highly dependent on the problem
configuration including farfield conditions, the geometry of the
body, and fluid properties. As the flow speed increases, accurate
resolution of the lead shock in a computational setting is paramount
to predict quantities of interest integrated over the body,
particularly aerodynamic heating
\cite{gnoffo_computational_2003,candler_unstructured_2007}.
This has forced computational fluid dynamics researchers and
practitioners to invest substantial effort and resources to
generate hexahedral meshes with very tight grid spacing
near shocks and elements aligned to the curvature of the
lead shock
\cite{lee_spurious_1999,gnoffo_computational_2003,candler_current_2009,gnoffo_multi-dimensional_2009,candler_advances_2015}.
Given the substantial amount of user-intensive effort
\cite{gnoffo_computational_2003,candler_current_2009}
required to mesh and simulate a single configuration, \textit{many-query}
analyses such as optimization, uncertainty quantification, parameter
sweeps, ``what-if'' scenarios, or even parameter-based continuation,
remain a significant challenge because they cause the bow shock to move,
which requires modifications to the mesh.


Shock capturing is a popular and effective class of approaches to stabilize
higher-than-first-order methods in the vicinity of shocks on a fixed
computational grid. Limiters, which are used to limit the solution gradient
near shocks, are commonly used with second-order finite volume methods
\cite{van1979towards} and high-order discontinuous Galerkin (DG) methods
\cite{cockburn_rungekutta_2001}. These methods are commonly 
used in real-world flow simulations; however, they place stringent demands
on the computational mesh for high-speed flows (hexahedral elements,
tight grid spacing near shocks, alignment of elements with lead shock)
\cite{candler_advances_2015} and are not well-suited for many-query
studies involving parametrized shocks. Weighted essentially non-oscillatory
(WENO) methods
\cite{harten_uniformly_1987, liu_weighted_1994, jiang_efficient_1996}
use stencil-based high-order reconstruction near shocks to mitigate spurious 
oscillations and can lead to crisp shocks on structured meshes. They have
shown good agreement with experiments for high-speed flows
\cite{shen_simulation_2010, lee_hypersonic_2010}, although
they are not well-suited for complex domains that require unstructured meshes.
For high-order methods, artificial viscosity approaches can smoothly resolve 
steep gradients with sub-cell accuracy and is the preferred shock capturing
approach for finite-element-based methods 
\cite{persson_sub-cell_2006,barter_shock_2010,
      fernandez_physics-based_2018,ching_shock_2019}.
The combination of high-order DG methods with artificial viscosity has
even been shown to reduce the sensitivity of hypersonic flow simulations
to the choice of numerical flux and grid alignment \cite{ching_shock_2019}.
However, artificial viscosity models usually suffer from a relatively
strong dependence on a large amount of empirical parameters that must be tuned 
\cite{yu_artificial_2020} and require substantial grid refinement near shocks
where accuracy has been dropped to first order. As such, these methods are
not ideally suited for many-query analyses involving parametrized shocks.

An alternative approach is \emph{shock tracking} or \emph{shock fitting}, 
where the computational mesh is moved to align faces of mesh elements with
solution discontinuities to represent them geometrically with the inter-element
jump in the solution basis without requiring additional stabilization.
This leads to accurate solutions on coarse meshes when using high-order 
methods and avoids instabilities often associated with shock capturing
methods such as carbuncles. The traditional approach to shock tracking
\cite{moretti2002thirty, salas2009shock} largely consists of explicitly
identifying shock locations and using the Rankine-Hugoniot conditions to
determine the shock motion and states across the shock. Although research
in this area has experienced a resurgence in recent years
\cite{ciallella_extrapolated_2020, bonfiglioli_unsteady_2016, 
geisenhofer_extended_2020, daquila_novel_2021, italian2021CnF1, 
italian2021CnF2, italian2022cmame, italian2022CPC}, these methods
are most suited for flows with relatively simple shock geometries
because they require explicit meshing of the shocks and specialized
strategies to track the shocks separately from the remainder of the flow.

A different approach to shock tracking known as
\textit{implicit shock tracking}, which includes the
High-Order Implicit Shock Tracking (HOIST) method
\cite{zahr_optimization-based_2018,zahr_implicit_2020, zahr_high-order_2020}
and the Moving Discontinuous Galerkin Method with Interface Condition
Enforcement (MDG-ICE)
\cite{corrigan_moving_2019,kercher_least-squares_2020,kercher_moving_2021},
has overcome some of the challenges of traditional approaches
to shock tracking.
These methods discretize the conservation law on a
shock-agnostic mesh and pose an optimization problem whose solution is the
discontinuity-aligned mesh and corresponding discretized flow solution.
That is, shock tracking is achieved implicitly through the solution of an
optimization problem. The meshing challenge from traditional shock
tracking approaches is largely circumvented, which leads to a general
approach that is independent of the underlying conservation law.
These methods have been shown to reliably and effectively
solve inviscid and viscous, steady and unsteady, inert and reacting
flows of varying degrees of shock complexity
\cite{corrigan_moving_2019, zahr_high-order_2020, shi2022implicit,
huang2022robust, corrigan2019scitech, corrigan2019avi, corrigan2021avi}.
Because these methods inherently and automatically align the computational
grid with shocks in the domain, they are well-suited to many-query studies
involving parametrized shocks.

We propose a novel framework based on implicit shock tracking,
specifically the HOIST method, to efficiently and robustly conduct
many-query analyses of problems involving parameter-dependent
\textit{lead shocks}. We define a lead shock as a shock (or,
more generally, any non-smooth feature) separating the boundary
state from the downstream flow. Bow shocks and primary blast
waves are examples of lead shocks; however, in the current
setting, even the head of a rarefaction wave in a shock tube qualifies.
The approach is initialized by applying the HOIST method on a
shock-agnostic mesh at one parameter configuration of interest
to generate a shock-aligned mesh and the corresponding flow field.
Because the grid is aligned with the lead shock, the
solution in all elements upstream of the shock will be a constant
equal to the boundary state. As such, all elements upstream of the shock
are removed to produce a reduced mesh that will be used for all subsequent
parameter configurations, with the farfield boundary condition directly
applied on the new shock boundary. In our numerical experiments, this
can reduce the size of the computational mesh by three-fold, which is
additional computational savings on top of the coarse meshes required
by implicit shock tracking approaches \cite{huang2022robust}.

With only the portion of the domain downstream of the lead shock modeled,
there is no reason to directly optimize for all nodal coordinates in the mesh
if the lead shock is the only non-smooth feature in the domain. Instead,
we only optimize the positions of the nodes on the shock boundary
with all other nodal coordinates determined by boundary constraints
and partial differential equation (PDE) based smoothing using the
approach in \cite{huang2021mshreduct}. For problems with secondary
non-smooth features in addition to the lead shock, all downstream
nodal positions are optimized to eliminate the need for nonlinear
stabilization and provide highly accurate solutions. In addition to
substantially reducing the overall degrees of freedom (DoFs) of the
implicit shock tracking discretization, it also improves robustness
and accelerates convergence because the overall tracking problem is
easier and a high-quality initial guess is provided from
the solution at previous parameter configurations.
The proposed framework can be used for most \textit{many-query}
applications involving parametrized lead shocks such as optimization,
uncertainty quantification, parameter sweeps, ``what-if'' scenarios,
or parameter-based continuation. In the continuation setting, we
outline and demonstrate a procedure to leverage partially converged
solves at intermediate stages to improve the efficiency of the approach.

The remainder of the paper is organized as follows.
Section~\ref{sec:govern} introduces the governing system of inviscid
conservation laws, its reformulation on a fixed reference domain,
and its discretization using a discontinuous Galerkin method.
Section~\ref{sec:opt-form} provides a brief summary of the HOIST 
formulation with targeted mesh optimization proposed in 
\cite{zahr_implicit_2020, huang2022robust, huang2021mshreduct}. 
Section~\ref{sec:bowshk} introduces the specialized HOIST solver
for abstract many-query problems involving parametrized shocks
and its specialization to parameter-based continuation with partially
converged intermediate stages. Finally, Section~\ref{sec:numexp}
demonstrates the robustness and flexibility of the proposed approach
for Mach continuation of two- and three-dimensional supersonic and
hypersonic problems, and for a parameter sweep of a Riemann problem
(Euler equations) parametrized by its initial condition.

%

\section{Governing equations and high-order discretization}
\label{sec:govern}
In this section, we introduce a system of steady conservation laws
whose solution will be assumed to contain a lead shock (Section~\ref{sec:govern:claw}).
Next, we transform the system of conservation laws to a reference domain such that
domain deformations appear explicitly in the governing equations
(Section~\ref{sec:govern:transf}), and discretize the transformed
equations using a high-order DG method (Section~\ref{sec:govern:dg}).

\subsection{System of conservation laws}
\label{sec:govern:claw}
Consider a general system of $m$ inviscid conservation laws over
$\Omega \subset \Rbb^d$
\begin{equation} \label{eqn:claw-phys}
 \nabla\cdot F(U) = S(U) \quad \text{in}~~\Omega,
\end{equation}
where $\func{U}{\Omega}{\Rbb^m}$ is implicitly defined as the solution
of (\ref{eqn:claw-phys}), $\func{F}{\Rbb^m}{\Rbb^{m\times d}}$ is the flux
function, and $\func{S}{\Rbb^m}{\Rbb^m}$ is the source term. In general,
the solution $U(x)$ may contain discontinuities, in which case the
conservation law (\ref{eqn:claw-phys}) holds away from the discontinuities
and the Rankine-Hugoniot conditions hold at discontinuities.
In this work, we focus on problems containing a \textit{lead shock} and, for
concreteness, focus on the compressible Euler equations (Section~\ref{sec:numexp}).
However, the method developed in this work applies to any conservation law of the
form (\ref{eqn:claw-phys}) whose solution contains a lead shock.

\subsection{Transformed system of conservation laws on a fixed reference domain}
\label{sec:govern:transf}
Let $\Gbb$ be the collection of diffeomorphisms from a reference domain
$\Omega_0$ to the physical domain $\Omega$, i.e., for $\Gcal\in\Rbb$, we have
\begin{equation} \label{eqn:dom-map}
 \func{\Gcal}{\Omega_0}{\Omega}, \quad
 \Gcal : X \mapsto \Gcal(X).
\end{equation}
Following the approach in \cite{zahr_optimization-based_2018},
for any $\Gcal\in\Gbb$, the conservation law on the physical
domain $\Omega$ is transformed to a conservation law on
the reference domain $\Omega_0$ as
\begin{equation} \label{eqn:claw-ref}
 \bar\nabla\cdot\bar{F}(\bar{U};G) =
 \bar{S}(\bar{U};g) \quad \text{in}~~\Omega_0,
\end{equation}
where $\bar\nabla$ is the gradient operator on the reference domain, 
$\func{\bar{U}}{\bar\Omega}{\Rbb^m}$ is the transformed solution,
$\func{\bar{F}}{\Rbb^m\times\Rbb^{d\times d}}{\Rbb^{m\times d}}$ is
the transformed flux function, $\func{\bar{S}}{\Rbb^m\times\Rbb}{\Rbb^m}$
is the transformed source term, and
\begin{equation}
 G = \bar\nabla \Gcal, \qquad g = \det G
\end{equation}
are the deformation gradient and Jacobian, respectively, of the mapping
$\Gcal\in\Gbb$. For any $X\in\Omega_0$, the transformed and physical
solution are related as
\begin{equation}
 \bar{U}(X) = U(\Gcal(X)),
\end{equation}
and the transformed flux and source term are defined as
\begin{equation}
 \bar{F} : (\bar{W}; \Theta) \mapsto (\det\Theta) F(\bar{W}) \Theta^{-T},
 \qquad
 \bar{S} : (\bar{W}; q) \mapsto q S(\bar{W}).
\end{equation}

\subsection{Discontinuous Galerkin discretization of
            the transformed conservation law}
\label{sec:govern:dg}
We use a standard nodal discontinuous Galerkin method 
to discretize the transformed conservation law (\ref{eqn:claw-ref}).
Let $\Ecal_h$ represent a discretization of the reference domain
$\Omega_0$ into non-overlapping, potentially curved, computational elements.
The DG trial space of discontinuous piecewise polynomials associated with 
the mesh $\Ecal_h$ is defined as 
\begin{equation}
 \Vcal_h^p = \left\{v \in [L^2(\Omega_0)]^m \suchthat
         \left.v\right|_K \in [\Pcal_p(K)]^m,
         ~\forall K \in \Ecal_h\right\},
\end{equation}
where $\Pcal_p(K)$ is the space of polynomial functions of degree at most
$p \geq 1$ on the element $K$. The space of globally
continuous piecewise polynomials of degree $q$ associated with the mesh
$\Ecal_h$ is defined as
\begin{equation} \label{eqn:gfcnsp}
  \Wcal_h = \left\{v \in C^0(\Omega_0) \suchthat
         \left.v\right|_K \in \Pcal_q(K),~\forall K \in \Ecal_h\right\};
\end{equation}
we discretize the domain mapping ($\Gcal\in\Gbb$) with the corresponding
vector-valued space $[\Wcal_h]^d$.
Taking the DG test space to be $\Vcal_h^{p'}$, where $p'\geq p$, the
DG formulation is: given $\Gcal_h\in[\Wcal_h]^d$, find
$\bar{U}_h\in\Vcal_h^p$ such that for all $\bar\psi_h\in\Vcal_h^{p'}$, we have
\begin{equation} \label{eqn:claw-weak-elem}
 \int_{\partial K} \bar\psi_h^+ \cdot
   \bar\Hcal(\bar{U}_h^+,\bar{U}_h^-,N_h;\bar\nabla\Gcal_h) \, dS -
 \int_K \bar{F}(\bar{U}_h; \bar\nabla\Gcal_h):\bar\nabla \bar\psi_h \, dV =
 \int_K \bar\psi_h \cdot \bar{S}(\bar{U}_h; \det(\bar\nabla\Gcal_h)) \, dV,
\end{equation}
where $N_h$ is the unit outward normal to
element $K\in\Ecal_h$, $\bar{W}_h^+$ ($\bar{W}_h^-$) denotes the interior
(exterior) trace of $\bar{W}_h$ to the element $K$ for $\bar{W}_h\in\Vcal_h^s$
(any $s\in\{0,1,\dots\}$), and $\bar\Hcal$ is the numerical flux function
associated with the reference inviscid flux $\bar{F}$; see
\cite{huang2022robust} for additional details.

After introducing a basis for the test, trial, and domain mapping spaces,
the governing DG equations (\ref{eqn:claw-weak-elem}) with $p'=p$ reduces
to the algebraic residual form
\begin{equation} \label{eqn:disc_dg}
  \func{\rbm}{\Rbb^{N_{\ubm}} \times \Rbb^{N_{\xbm}}}{\Rbb^{N_{\ubm}}}, \qquad 
  \rbm : (\ubm,\xbm) \mapsto \rbm(\ubm,\xbm), 
\end{equation} 
where $N_{\ubm}=\mathrm{dim}\Vcal^p_h$ and
$N_{\xbm}=\mathrm{dim}([\Wcal_h]^d)$; $\ubm\in\Rbb^{N_{\ubm}}$
is the vector of flow field coefficients and $\xbm\in\Rbb^{N_{\xbm}}$
is the vector of nodal coordinates of the mesh (also called mesh DoFs).
Notice that for a fixed mesh $\xbm$, (\ref{eqn:disc_dg}) is a
standard DG discretization residual. In addition, we define the
algebraic enriched residual associated with a test space of
degree $p' = p+1$ as 
\begin{equation} \label{eqn:disc_dg1}
 \func{\Rbm}{\Rbb^{N_{\ubm}} \times \Rbb^{N_{\xbm}}}{\Rbb^{N_{\ubm}'}}, \qquad 
 \Rbm : (\ubm,\xbm) \mapsto \Rbm(\ubm,\xbm),
\end{equation}
where $N_{\ubm}'=\mathrm{dim}\Vcal_h^{p'}$, which will later be used to 
construct the HOIST objective function. 

\section{The High-Order Implicit Shock Tracking (HOIST) method}
\label{sec:opt-form}
In this section, we provide a brief summary of the HOIST method with 
targeted mesh optimization \cite{zahr_implicit_2020, huang2022robust}, 
which allows only a selected portion of the mesh DoFs to be optimized
while aligning mesh faces with non-smooth solution features and
preserving boundaries. This is achieved by partitioning the mesh DoFs
$\xbm\in\Rbb^{N_\xbm}$ as
\begin{equation} \label{eqn:mshpart1}
 \xbm = (\xbm_\mathrm{c}, \xbm_\mathrm{u}), \qquad
 \xbm_\mathrm{u} = (\ybm, \xbm_\mathrm{s})
\end{equation}
where $\xbm_\mathrm{c}\in\Rbb^{N_\xbm^\mathrm{c}}$ are the constrained
DoFs and $\xbm_\mathrm{u}\in\Rbb^{N_\xbm^\mathrm{u}}$ are the unconstrained
DoFs. Following \cite{huang2022robust}, $\xbm_\mathrm{u}$ can be
freely chosen (e.g., to align element faces with non-smooth features),
whereas $\xbm_\mathrm{c}$ is uniquely determined from $\xbm_\mathrm{u}$.
Following \cite{huang2021mshreduct}, unconstrained DoFs are further
partitioned into optimized DoFs $\ybm\in\Rbb^{N_\ybm}$
and smoothed DoFs $\xbm_\mathrm{s}\in\Rbb^{N_\xbm^\mathrm{s}}$, where
$\ybm$ will be optimized for shock alignment and
$\xbm_\mathrm{s}$ will be determined through PDE-based smoothing.
For abstraction, we let $\phibold$ be a parametrization of the mesh DoFs
\begin{equation} \label{eqn:mapparam}
 \func{\phibold}{\Rbb^{N_\ybm}}{\Rbb^{N_\xbm}}, \qquad
 \phibold: \ybm \mapsto \phibold(\ybm),
\end{equation} 
that maps the optimized mesh DoFs to all mesh DoFs, i.e.,
$\xbm = \phibold(\ybm)$. The parametrization must be constructed
to ensure
\begin{inparaenum}[1)]
 \item optimized mesh DoFs ($\ybm$) can move freely,
 \item nodes on fixed domain boundaries can only slide along
  those boundaries by computing the constrained mesh DoFs ($\xbm_\mathrm{c}$)
  from $\ybm$, and
 \item smoothed mesh DoFs ($\xbm_\mathrm{s}$) are determined through
  PDE-based smoothing (e.g., linear elasticity equations with pure Dirichlet
  boundary conditions).
\end{inparaenum}
A complete description of the construction of $\phibold$ can be found in
\cite{zahr_implicit_2020, huang2022robust} (boundary preservation only)
and \cite{huang2021mshreduct} (boundary preservation and PDE-based smoothing).

The HOIST method is formulated as an optimization problem 
over the DG solution coefficients and the optimized mesh DoFs as
\begin{equation} \label{eqn:pde-opt}
 (\ubm^\star,\ybm^\star) \coloneqq
 \argmin_{\ubm\in\Rbb^{N_\ubm},\ybm\in\Rbb^{N_\ybm}} f(\ubm,\phibold(\ybm)) \quad \text{subject to:} \quad \rbm(\ubm,\phibold(\ybm)) = \zerobold,
\end{equation}
where $\func{f}{\Rbb^{N_\ubm}\times\Rbb^{N_\xbm}}{\Rbb}$ is the
objective function defined in \cite{zahr_implicit_2020,huang2022robust}
and $\xbm^\star = \phibold(\ybm^\star)$ are the nodal coordinates 
of the discontinuity-aligned mesh.
The objective function consists of two terms as
\begin{equation}\label{eqn:obj0}
 f : (\ubm,\xbm) \mapsto f_\text{err}(\ubm,\xbm) + \kappa^2 f_\text{msh}(\xbm),
\end{equation}
where $\func{f_\mathrm{err}}{\Rbb^{N_\ubm}\times\Rbb^{N_\xbm}}{\Rbb}$
is the alignment term and $\func{f_\mathrm{msh}}{\Rbb^{N_\xbm}}{\Rbb}$
is a mesh quality term, defined as
\begin{equation}
 f_\mathrm{err} : (\ubm,\xbm) \mapsto \frac{1}{2}\norm{\Rbm(\ubm,\xbm)}_2^2,
 \qquad
 f_\mathrm{msh} : \xbm \mapsto \frac{1}{2}\norm{\Rbm_\mathrm{msh}(\xbm)}_2^2,
\end{equation}
$\func{\Rbm_\mathrm{msh}}{\Rbb^{N_\xbm}}{\Rbb^{|\Ecal_h|}}$ is 
an element-wise mesh distortion residual defined in \cite{huang2022robust},
and $\kappa$ is a penalty parameter that balances the alignment and mesh
quality objectives; see \cite{huang2022robust} for a complete definition
of $\Rbm_\mathrm{msh}$ and an adaptive algorithm to set $\kappa$.
The norm of the enriched DG residual has proven to
be an effective alignment indicator \cite{zahr_implicit_2020,huang2022robust}
because it penalizes non-physical oscillations that will arise on
meshes that do not fit solution discontinuities.

A sequential quadratic programming (SQP) method with a modified
Levenberg-Marquardt Hessian approximation introduced in
\cite{corrigan_moving_2019, zahr_implicit_2020} is used to
solve the optimization problem (\ref{eqn:pde-opt}). The SQP
solver simultaneously converges the optimized mesh DoFs ($\ybm$)
and the DG solution coefficients ($\ubm$) to their optimal
values, i.e., a high-order DG solution ($\ubm^\star$) on a
discontinuity-aligned mesh ($\phibold(\ybm^\star)$). This is
accomplished by combining the solution and optimized mesh DoFs 
into a single vector of optimization variables
$\zbm = (\ubm,\ybm)\in\Rbb^{N_\zbm}$ ($N_\zbm=N_\ubm+N_\ybm$)
and generating a sequence of iterates as
\begin{equation} \label{eqn:sqp_step0}
 \zbm_{i+1} = \zbm_i + \alpha_i\Delta\zbm_i
\end{equation}
for $i=0,1,\dots$, where $\zbm_i\in\Rbb^{N_\zbm}$ is the vector of
optimization variables at iteration $i$, $\Delta\zbm_i\in\Rbb^{N_\zbm}$
is the search direction at the $i$th iteration, and $\alpha_i\in\Rbb_{>0}$
is the step length at the $i$th iteration. The search direction is determined
at each iteration by solving a quadratic approximation to (\ref{eqn:pde-opt})
at $\zbm_i$ and the step length is determined via a line search of an
$\ell_1$ merit function. A complete description of the HOIST method and
solver can be found in \cite{zahr_implicit_2020,huang2022robust}.

\section{The HOIST method for flows with parametrized shocks}
\label{sec:bowshk}
In this section, we introduce a framework for solving flows with
parametrized lead shocks using the HOIST method. This approach can be
leveraged for \textit{many-query} analyses such as optimization,
uncertainty quantification, parameter sweeps, or ``what-if'' scenarios,
or to drive a continuation strategy for more complex flow regimes.

\subsection{Parametrized setting}
\label{sec:bowshk:param}
While not introduced as such for brevity, all terms in the conservation
law (\ref{eqn:claw-phys}) as well as its transformed (\ref{eqn:claw-ref})
and discrete variants (\ref{eqn:disc_dg})-(\ref{eqn:disc_dg1})
depend on \textit{problem data} such as boundary condition,
domain geometry, material properties, etc. Many-query analyses and
continuation-based solvers inherently vary one or more of these
parameters for some purpose, e.g., find an optimal solution (optimization),
compute a probability distribution (uncertainty quantification), solve
a complex problem using a sequence of easier ones (continuation).
To handle these settings, we explicitly introduce dependence of the
discretized residuals on a collection of parameters. That is, we
redefine the standard DG residual in (\ref{eqn:disc_dg}) as
\begin{equation}
 \func{\rbm}{\Rbb^{N_\ubm}\times\Rbb^{N_\xbm}\times\Dcal}{\Rbb^{N_\ubm}}, \qquad
 \rbm : (\ubm,\xbm;\mubold) \mapsto \rbm(\ubm,\xbm;\mubold)
\end{equation}
and the enriched DG residual in (\ref{eqn:disc_dg1}) as
\begin{equation}
 \func{\Rbm}{\Rbb^{N_\ubm}\times\Rbb^{N_\xbm}\times\Dcal}{\Rbb^{N_\ubm'}},\qquad
 \Rbm : (\ubm,\xbm;\mubold) \mapsto \Rbm(\ubm,\xbm;\mubold),
\end{equation}
where $\Dcal\subset\Rbb^{N_\mubold}$ is a space of admissible parameter
configurations and $\mubold\in\Dcal$ is a vector of parameters that either
directly or indirectly defines relevant problem data. 
For this work, we assume that for any
$\mubold\in\Dcal$, the flow contains a lead shock that separates the farfield
from shock-downstream portion of the domain.  Furthermore, we assume
topologically equivalent or similar lead shocks for all $\mubold\in\Dcal$.

From these definitions, the HOIST solution in (\ref{eqn:pde-opt}) becomes
parameter dependent, and we let
$(\ubm_\mubold^\star,\ybm_\mubold^\star)\subset\Rbb^{N_\ubm}\times\Rbb^{N_\ybm}$
denote the solutions of (\ref{eqn:pde-opt}) with $\rbm$ and $\Rbm$ replaced by
$\rbm(\,\cdot\,,\,\cdot\,;\mubold)$ and $\Rbm(\,\cdot\,,\,\cdot\,;\mubold)$,
respectively. Because the optimization problem in (\ref{eqn:pde-opt}) is
non-convex, it will have multiple local minima. However, from a given
starting point $(\bar\ubm,\bar\ybm)$, the HOIST solver \cite{huang2022robust}
will return a single solution. We let
$\func{\Upsilon}{\Rbb^{N_\ubm}\times\Rbb^{N_\ybm}\times\Dcal}{\Rbb^{N_\ubm}\times\Rbb^{N_\ybm}}$
be an operator that maps the starting point $(\bar\ubm,\bar\ybm)$
and parameter configuration $\mubold$ to the element of the set
$(\ubm_\mubold^\star,\ybm_\mubold^\star)$ returned by the HOIST solver, i.e.,
\begin{equation}
 \Upsilon : (\bar\ubm,\bar\ybm,\mubold) \mapsto \Upsilon(\bar\ubm,\bar\ybm,\mubold) \in (\ubm_\mubold^\star,\ybm_\mubold^\star).
\end{equation}

\subsection{Many-query analysis}
\label{sec:bowshk:many}
In the remainder, we consider an abstract many-query setting where
our goal is to determine the solution of the conservation
law for all parameter configurations in an ordered subset $\Xi \subset \Dcal$
with $\Xi = \{\mubold_1,\dots,\mubold_N\}$. Depending on the many-query
setting, the set may be known \textit{a priori} (e.g., parameter sweeps,
static continuation, Monte Carlo approaches to uncertainty quantification),
or built adaptively (e.g., optimization, adaptive continuation, ``what-if''
scenarios). The solutions of the conservation law will be computed sequentially
according to the ordering of the parameters, where the initial guess
for the flow and mesh DoFs will come from the solution of the previous
parameter. Optimal ordering of the parameters will depend on the many-query
setting and the problem under consideration.

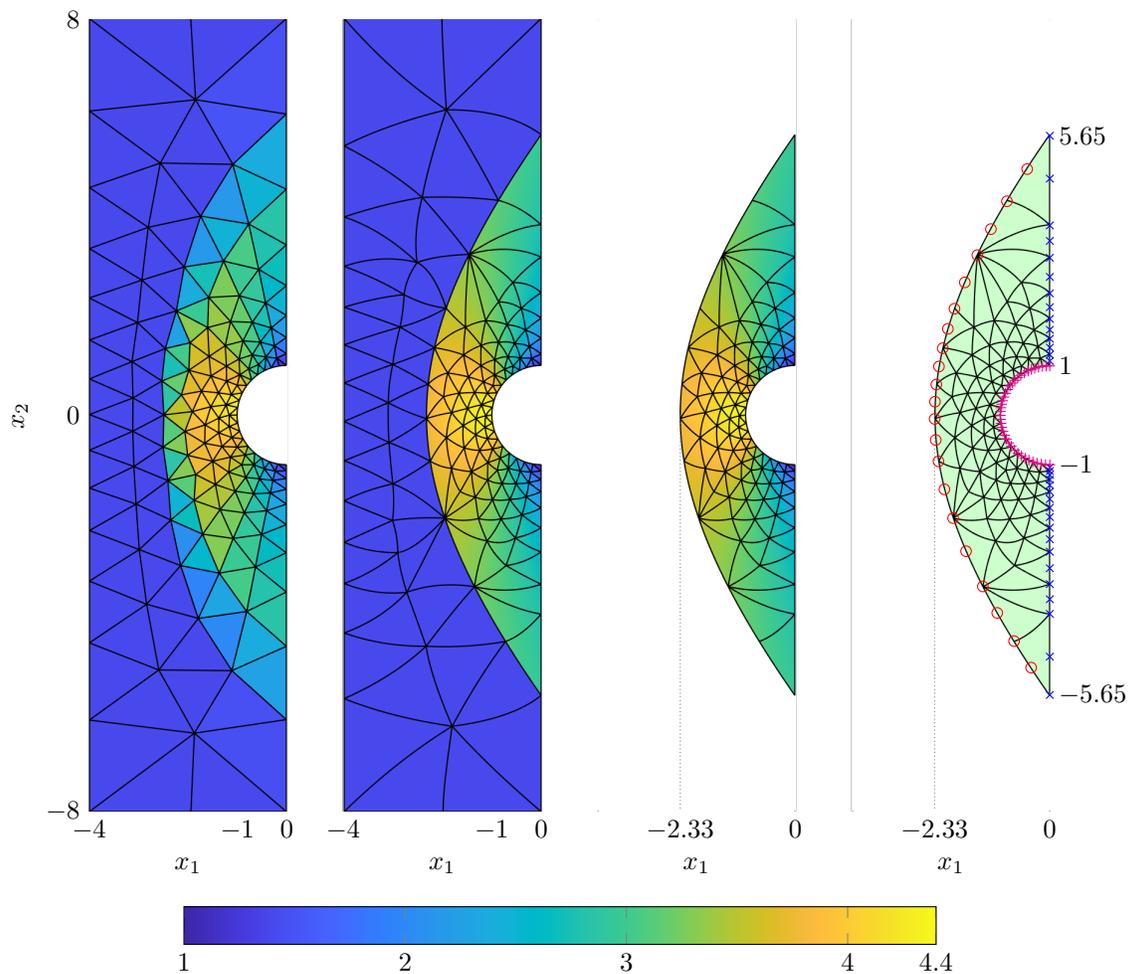
\begin{figure}
\centering
\ifbool{fastcompile}{}{

\begin{tikzpicture}
\begin{groupplot} [
group style={group size = 4 by 1, horizontal sep = 0.75cm, vertical sep = 0.5cm}]
\nextgroupplot[axis equal image, xtick={-4, -1, 0}, xlabel={$x_1$}, ytick={-8, 0, 8}, ylabel={$x_2$}, yticklabel pos=left, xmin=-4, xmax=0, ymin=-8, ymax=8, width=0.85\textwidth]
\addplot []
graphics [xmin=-4,xmax=0,ymin=-8,ymax=8] { 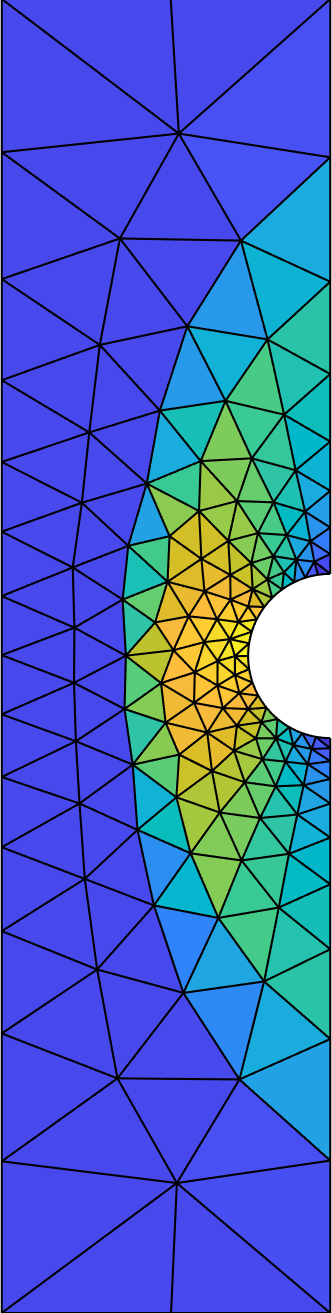};

\nextgroupplot[axis equal image, xtick={-4, -1, 0}, xlabel={$x_1$}, yticklabels={,,}, ylabel={}, xmin=-4, xmax=0, ymin=-8, ymax=8, width=0.85\textwidth]
\addplot []
graphics [xmin=-4,xmax=0,ymin=-8,ymax=8] { 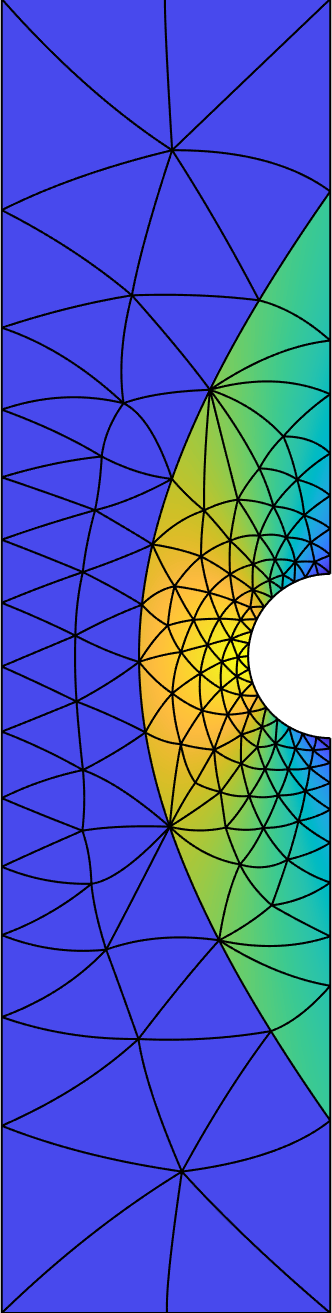};

\nextgroupplot[axis equal image, xtick={-2.33, 0}, xlabel={$x_1$}, yticklabels={,,}, ylabel={}, xmin=-4, xmax=0, ymin=-8, ymax=8, width=0.85\textwidth]
\addplot []
graphics [xmin=-4,xmax=0,ymin=-8,ymax=8] { 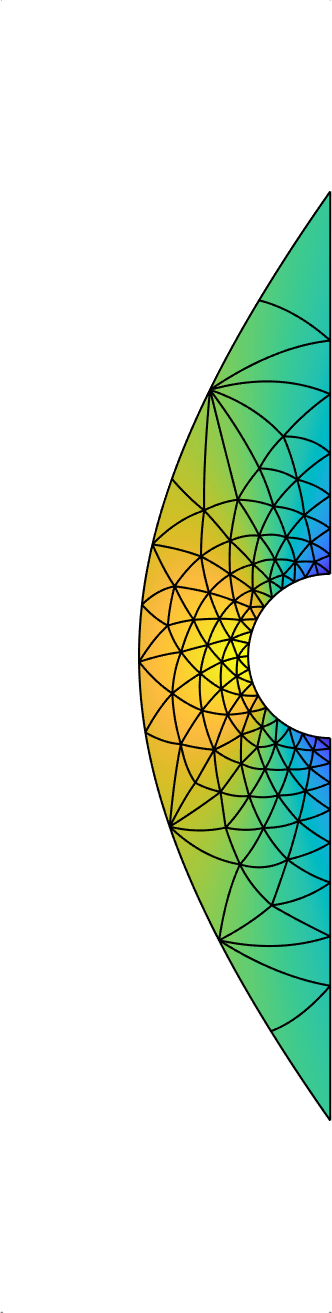};

\addplot [gray, densely dotted, thin, forget plot]
coordinates {
(-2.33000000e+00,  0.00000000e+00)
(-2.33000000e+00, -8.00000000e+00)};

\nextgroupplot[axis equal image, xtick={-2.33, 0}, xlabel={$x_1$}, ylabel={}, ytick={-5.65, -1, 1, 5.65}, yticklabel pos=right, xmin=-4, xmax=0, ymin=-8, ymax=8, width=0.85\textwidth]
\addplot []
graphics [xmin=-4,xmax=0,ymin=-8,ymax=8] { 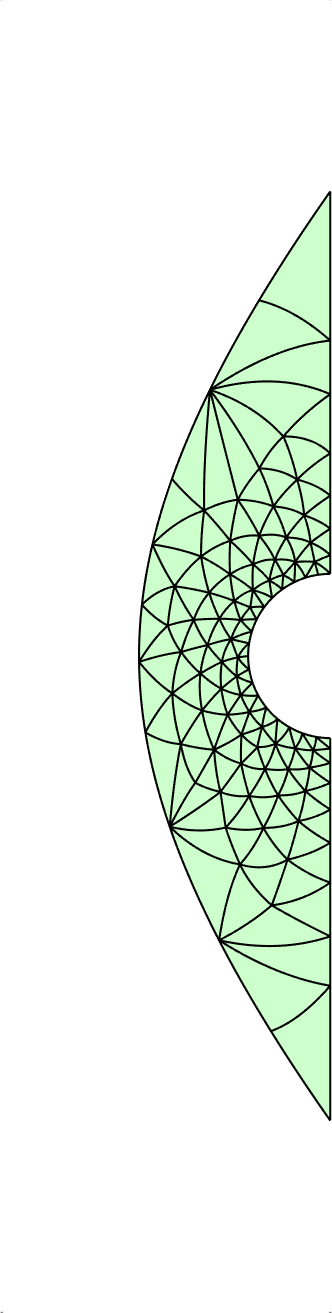};

\addplot [gray, densely dotted, thin, forget plot]
coordinates {
(-2.33000000e+00,  0.00000000e+00)
(-2.33000000e+00, -8.00000000e+00)};

\addplot [mark=o, mark size=2, mark options={thin}, red, only marks]
coordinates {
(-2.29338763e+00,  6.24074448e-01)
(-2.32154601e+00,  2.75556342e-01)
(-2.25106053e+00, -9.30050713e-01)
(-2.30598895e+00, -4.99031493e-01)
(-2.24116385e+00,  9.88242673e-01)
(-2.32792438e+00, -7.19542912e-02)
(-2.16506711e+00,  1.35833999e+00)
(-2.13113545e+00, -1.49577687e+00)
(-1.95478993e+00, -2.07899180e+00)
(-2.06054847e+00,  1.74995593e+00)
(-1.92768796e+00,  2.15635003e+00)
(-1.69291991e+00, -2.74898573e+00)
(-1.71846530e+00,  2.68982176e+00)
(-1.46468462e+00,  3.24116134e+00)
(-1.19102359e+00,  3.76780238e+00)
(-7.22133952e-01, -4.56770948e+00)
(-8.66836166e-01,  4.33075893e+00)
(-1.06302916e+00, -3.99667506e+00)
(-1.35357745e+00, -3.46194232e+00)
(-3.77515962e-01, -5.10408277e+00)
(-4.58231860e-01,  4.98141757e+00)};\label{line:inlet_nodes}

\addplot [mark=x, mark size=2, mark options={thin}, blue, only marks]
coordinates {
(-1.40133104e-18, -4.01567022e+00)
(-1.33887776e-18,  3.84099142e+00)
(-9.72703223e-19, -3.72036543e+00)
(-9.14516230e-19,  3.52887830e+00)
(-1.16138009e-18, -4.88223480e+00)
(-1.08476267e-18,  4.78967721e+00)
(-8.91797457e-19, -3.41383236e+00)
( 7.01814928e-18,  3.18796575e+00)
(-1.02521930e-18, -5.65663721e+00)
(-6.22113037e-18, -3.09342058e+00)
( 6.34867010e-18,  2.80153040e+00)
(-8.52493503e-19,  5.65680583e+00)
(-2.08441858e-18, -2.75866975e+00)
(-6.09729457e-19,  2.46701317e+00)
(-5.68799467e-19, -2.50646095e+00)
(-5.25706424e-19,  2.18790525e+00)
(-4.54674680e-19,  1.95195182e+00)
(-5.65940383e-19, -2.27368546e+00)
(-1.74566008e-19,  1.73117916e+00)
(-2.17974356e-19, -2.05690949e+00)
( 3.59343450e-20, -1.05431205e+00)
( 6.34632130e-20,  1.09950011e+00)
( 8.38557696e-20, -1.12674166e+00)
( 8.84499656e-20, -1.20526862e+00)
( 9.98320970e-20, -1.30212534e+00)
( 1.23438742e-19, -1.41161351e+00)
( 1.51949463e-19, -1.54384690e+00)
( 1.59871127e-19,  1.22576591e+00)
( 1.79382146e-19, -1.69278826e+00)
( 9.26840597e-20,  1.37032208e+00)
( 8.38604157e-20, -1.86887239e+00)
( 9.19689309e-21,  1.54994869e+00)};\label{line:outlet_nodes}

\addplot [mark=+, mark size=2, mark options={thin}, magenta, only marks]
coordinates {
(-1.00000000e+00, -2.79025549e-21)
(-9.97012089e-01,  7.72456754e-02)
(-9.96715901e-01, -8.09778563e-02)
(-9.88066211e-01,  1.54029744e-01)
(-9.86885174e-01, -1.61423834e-01)
(-9.74332297e-01,  2.25114584e-01)
(-9.69884311e-01, -2.43566053e-01)
(-9.55491249e-01,  2.95019446e-01)
(-9.46058981e-01, -3.23994451e-01)
(-9.28555006e-01,  3.71194829e-01)
(-9.14809150e-01, -4.03886395e-01)
(-8.95556924e-01,  4.44946959e-01)
(-8.76826986e-01, -4.80806028e-01)
(-8.52734360e-01,  5.22344821e-01)
(-8.30692985e-01, -5.56730783e-01)
(-8.03239831e-01,  5.95655751e-01)
(-7.78002394e-01, -6.28261310e-01)
(-7.57098403e-01,  6.53300856e-01)
(-7.13359161e-01, -7.00798621e-01)
(-7.06829283e-01,  7.07384171e-01)
(-6.48380645e-01,  7.61316320e-01)
(-6.41981526e-01, -7.66720106e-01)
(-5.85831051e-01,  8.10433205e-01)
(-5.73146886e-01, -8.19452651e-01)
(-5.14206258e-01,  8.57666558e-01)
(-5.00002795e-01, -8.66023790e-01)
(-4.38796340e-01,  8.98586541e-01)
(-4.22031579e-01, -9.06581131e-01)
(-3.72349626e-01,  9.28092536e-01)
(-3.40800418e-01, -9.40135668e-01)
(-3.03934756e-01,  9.52692849e-01)
(-2.57909944e-01, -9.66168961e-01)
(-2.27060186e-01,  9.73880728e-01)
(-1.73072620e-01, -9.84909066e-01)
(-1.48741824e-01,  9.88876064e-01)
(-8.68646478e-02, -9.96220123e-01)
(-7.45786031e-02,  9.97215138e-01)
( 0.00000000e+00, -1.00000000e+00)
( 0.00000000e+00,  1.00000000e+00)};\label{line:cyl_nodes}

\end{groupplot}\end{tikzpicture}

\colorbarMatlabParula{1}{2}{3}{4}{4.4}
}
 \caption{Initial guess (density) on a shock-agnostic mesh of $M_\infty=2$
  flow over cylinder (\textit{left}), the shock-aligned mesh and corresponding
  solution obtained using HOIST (\textit{middle-left}), the corresponding
  solution and mesh extracted from downstream the bow shock
  (\textit{middle-right}), and a schematic visualizing the
  mesh parametrization (\textit{right}). Legend: shock boundary
  nodes that move freely (\ref{line:inlet_nodes}), outlet boundary
  nodes that are constrained to slide along original boundary
 (\ref{line:outlet_nodes}), cylinder boundary nodes that are fixed
 (\ref{line:cyl_nodes}), and all remaining nodes (not shown for clarity)
 are determined from PDE-based smoothing.}
\label{fig:2d:mach2}
\end{figure}

For the first parameter $\mubold_1$, we follow the standard HOIST
approach (no PDE-based smoothing) \cite{huang2022robust} where the
initial mesh DoFs ($\bar\ybm$) come directly from a mesh generator
and the initial flow solution ($\bar\ubm$) is the first-order finite
volume solution of the parametrized PDE at $\mubold=\mubold_1$ (or
any other reasonable initial guess). From
this initial guess, we compute the HOIST solution and make the observation
that the solution in all elements upstream of the shock is constant and
equal to the farfield boundary condition, making this a costly waste of
degrees of freedom. Therefore, we remove all elements upstream of the lead
shock and directly apply the farfield boundary condition to the
lead shock itself (referred to as the \textit{shock boundary}
in the remainder). This procedure is illustrated in Figure~\ref{fig:2d:mach2}.

For all subsequent parameters, only elements downstream of the lead
shock remain so changes to the lead shock position and shape come from
deformations to the shock boundary. For problems where the lead
shock is the \textit{only} shock for all $\mubold\in\Dcal$, we take
the optimized mesh DoFs ($\ybm$) to be all DoFs on the shock
boundary and the smoothed mesh DoFs ($\xbm_\mathrm{s}$) to be all
remaining unconstrained mesh DoFs (Figure~\ref{fig:2d:mach2}).
In this case, there is no need to allow topological changes to the
mesh (e.g., element collapses as in \cite{huang2022robust}) so the
mesh parametrization $\phibold$ will be fixed for all parameters
$\mubold_2,\dots,\mubold_N$. For problems with secondary shocks,
all unconstrained mesh DoFs are taken as optimized mesh DoFs.
For either case, nodes that lie on the intersection of the lead shock
boundary and other boundaries should respect the constraints of the
fixed boundaries (Figure~\ref{fig:2d:mach2}). In this work, we 
consider problems with only lead shocks as well as problems 
with both lead and secondary shocks.

With the initial HOIST solve, element removal, and mesh parametrization
settled, we state the many-query analysis algorithm. Let
$(\ubm_1^\star,\ybm_1^\star)$ denote the HOIST solution
after the elements upstream of the shock have been removed
and the nodal coordinates are parametrized on the new mesh.
Then, for $k=2,\dots,N$, we define the HOIST solution at
parameter $\mubold_k$ as
\begin{equation}
 (\ubm_k^\star,\ybm_k^\star) \coloneqq \Upsilon(\ubm_{k-1}^\star,\ybm_{k-1}^\star, \mubold_k).
\end{equation}
That is, the HOIST solution for $\mubold_k$ is initialized from the HOIST
solution at $\mubold_{k-1}$. Because the nodal coordinates are parametrized
with optimized mesh DoFs on the lead shock boundary, the HOIST method naturally
deforms that boundary to align with the new lead shock location
(Section~\ref{sec:numexp}). The complete algorithm is summarized in
Algorithm~\ref{alg:bndmov1}.

\begin{remark}\label{remark:newton}
The number of unconstrained mesh DoFs have been significantly
reduced by removing elements upstream of the shock and only
including nodes on the lead shock as unconstrained mesh DoFs
relative to the standard HOIST setting. This allows the
lead shock to be tracked with high accuracy, but can stall
deep convergence of the HOIST solver to tight optimality tolerances.
Rapid, deep convergence can be restored by terminating the HOIST
solver once the mesh has converged and updating $\ubm_k^\star$ with
the solution of
\begin{equation}
 \rbm(\,\cdot\,,\phibold(\ybm_k^\star);\mubold_k) = \zerobold
\end{equation}
starting from the HOIST output $\ubm_k^\star$, e.g., using
Newton-Raphson iterations. That is, the HOIST flow solution is
updated with a fixed-mesh DG solve.
Algorithms~\ref{alg:bndmov1}-\ref{alg:bndmov2} both include this
(optional) fixed mesh solve.
\end{remark}

\begin{algorithm}
 \caption{HOIST method for many-query analysis with parametrized lead shock}
 \label{alg:bndmov1}
 \begin{algorithmic}[1]
  \REQUIRE Reference mesh of entire domain $\bar\Ecal_h$, parameter set $\Xi=\{\mubold_1,\dots,\mubold_N\}$
  \ENSURE HOIST solution over $\Xi$: $\{(\ubm_1^\star,\ybm_1^\star),\dots,(\ubm_N^\star,\ybm_N^\star)\}$
  \STATE \textbf{Initial HOIST solve}: Compute HOIST solution $(\bar\ubm_1^\star,\bar\ybm_1^\star)$ at $\mubold_1$ over $\bar\Ecal_h$
  \STATE \textbf{Reduce mesh}: Create reduced mesh $\Ecal_h$ as
   $\bar\Ecal_h$ without elements upstream of the lead shock
  \STATE \textbf{Transfer solution}: Transfer solution $\bar\ubm_1^\star$
 to reduced mesh $\ubm_1^\star$
  \STATE \textbf{Parametrize reduced mesh}: Determine parametrization of
 reduced mesh $\phibold$ (Section~\ref{sec:bowshk:many})
  \STATE \textbf{Transfer mesh DoFs}: Transfer mesh DoFs $\bar\ybm_1^\star$ to
 reduced mesh $\ybm_1^\star$
  \FOR{$k=2,\dots,N$}
   \STATE \textbf{HOIST solve}: $(\ubm_{k}^\star, \ybm_{k}^\star) = \Upsilon(\ubm_{k-1}^\star, \ybm_{k-1}^\star, \mubold_k)$
   \STATE \textbf{Fixed mesh solve}: Solve $\rbm(\ubm,\phibold(\ybm_k^\star);\mubold_k) = \zerobold$ for $\ubm$ with initial guess $\ubm_k^\star$, replace $\ubm_k^\star \leftarrow \ubm$
  \ENDFOR
\end{algorithmic}
\end{algorithm}

\subsection{Application: Parameter continuation with early termination}
\label{sec:bowshk:contin}
Algorithm~\ref{alg:bndmov1} is well-suited for many-query settings that
require an accurate solution for every $\mubold\in\Xi$ such as parameter
sweeps and ``what-if'' scenarios. However, applications such as optimization
and continuation only need the solution at $\mubold_N$; solutions at all
previous parameters either aid in finding $\mubold_N$ or robustly
computing the solution at $\mubold_N$.
As such, the solution at intermediate parameters
$\mubold_1,\dots,\mubold_{N-1}$ can be computed approximately to improve
computational efficiency. We focus on the continuation setting because
the tolerances used for intermediate solutions in an optimization setting
are intimately tied to global convergence theory \cite{kouri2014inexact}
and beyond the scope of this work.

To this end, we introduce two convergence criteria based on (1) the number of
iterations and (2) relative reduction of the DG residual. Let
$(\ubm_k^i,\ybm_k^i)$ denote the $i$th iteration of the HOIST SQP
solver at $\mubold_k$ with
$(\ubm_k^0,\ybm_k^0) = (\ubm_{k-1}^\star,\ybm_{k-1}^\star)$.
Then, we define the modified HOIST operator
$\func{\Upsilon_{n,\xi}}{\Rbb^{N_\ubm}\times\Rbb^{N_\ybm}\times\Dcal}{\Rbb^{N_\ubm}\times\Rbb^{N_\ybm}}$
for $n > 1$ and $\xi < 1$ to be
\begin{equation}
  \Upsilon_{n,\xi} : (\ubm_k^0,\ybm_k^0,\mubold_k) \mapsto (\tilde\ubm_k^\star,\tilde\ybm_k^\star) \coloneqq (\ubm_k^I,\ybm_k^I),
\end{equation}
where $I \leq n$ is the smallest number such that
\begin{equation}
 \norm{\rbm(\ubm_k^I,\phibold(\ybm_k^I); \mubold_k)} \leq
 \xi \norm{\rbm(\ubm_k^0,\phibold(\ybm_k^0); \mubold_k)}.
\end{equation}
That is, $(\tilde\ubm_k^\star,\tilde\ybm_k^\star)$ is the output
of the HOIST SQP sover at parameter $\mubold_k$ after $n$ iterations
or when the residual converges to a relative tolerance of $\xi$,
whichever occurs first.
Then, we introduce iteration limits $n_1 \leq n_2$ and tolerances
$\xi_2 \leq \xi_1$, and use $\Upsilon_{n_1,\xi_1}$ for intermediate
parameters $k=1,\dots,N-1$ and $\Upsilon_{n_2,\xi_2}$ for the
last parameter $k=N$. The complete algorithm is summarized in
Algorithm~\ref{alg:bndmov2}.

\begin{algorithm}
 \caption{HOIST method for parameter continuation with parametrized lead shock}
 \label{alg:bndmov2}
 \begin{algorithmic}[1]
  \REQUIRE Reference mesh of entire domain $\bar\Ecal_h$, continuation set $\Xi=\{\mubold_1,\dots,\mubold_N\}$, tolerances ($\xi_1$, $\xi_2$), iteration limits ($n_1$, $n_2$)
  \ENSURE HOIST solution at $\mubold_N$: $(\ubm_N^\star,\ybm_N^\star)$
  \STATE \textbf{Initial HOIST solve}: Compute HOIST solution $(\bar\ubm_1^\star,\bar\ybm_1^\star)$ at $\mubold_1$ over $\bar\Ecal_h$
  \STATE \textbf{Reduce mesh}: Create reduced mesh $\Ecal_h$ as
   $\bar\Ecal_h$ without elements upstream of the lead shock
  \STATE \textbf{Transfer solution}: Transfer solution $\bar\ubm_1^\star$
 to reduced mesh $\tilde\ubm_1^\star$
  \STATE \textbf{Parametrize reduced mesh}: Determine parametrization of
 reduced mesh $\phibold$ (Section~\ref{sec:bowshk:many})
  \STATE \textbf{Transfer mesh DoFs}: Transfer mesh DoFs $\bar\ybm_1^\star$ to
 reduced mesh $\tilde\ybm_1^\star$
  \FOR{$k=2,\dots,N-1$}
   \STATE \textbf{Approximate HOIST solve}: $(\tilde\ubm_{k}^\star, \tilde\ybm_{k}^\star) = \Upsilon_{n_1,\xi_1}(\tilde\ubm_{k-1}^\star, \tilde\ybm_{k-1}^\star, \mubold_k)$
  \ENDFOR
   \STATE \textbf{Final parameter}: $(\ubm_{N}^\star, \ybm_{N}^\star) = \Upsilon_{n_2,\xi_2}(\tilde\ubm_{N-1}^\star, \tilde\ybm_{N-1}^\star, \mubold_N)$
   \STATE \textbf{Fixed mesh solve}: Solve $\rbm(\ubm,\phibold(\ybm_N^\star);\mubold_N) = \zerobold$ for $\ubm$ with initial guess $\ubm_N^\star$, replace $\ubm_N^\star \leftarrow \ubm$
\end{algorithmic}
\end{algorithm}

\section{Numerical experiments}
\label{sec:numexp}
In this section, we demonstrate the proposed HOIST framework for
shock-dominated flow problems parametrized by the farfield Mach number
(Section~\ref{sec:numexp:mach}) and initial condition
(Section~\ref{sec:numexp:ic}). We examine the robustness of the framework
under different solver settings to show the framework is able to accurately
track the lead shock across parameter configurations and provide
high-quality solutions on extremely coarse meshes. 

\subsection{Mach continuation}
\label{sec:numexp:mach}
We begin by considering a series of compressible, inviscid flows with a
bow shock that separates the boundary state from the downstream flow. We consider
several simple bluff body flows
(Section~\ref{sec:numexp:mach2to3}-\ref{sec:numexp:3dmach2to3}) that only
possess a bow shock to study the method and close with flow over a double
wedge geometry (Section~\ref{sec:numexp:dblwdg}) that has complex shock-shock
interactions downstream of the bow shock. For these problems, the flow through
the domain $\Omega\subset\Rbb^d$ is modeled using the steady Euler equations
\begin{equation} \label{eqn:euler}
\begin{split}
\pder{}{x_j}\left(\rho(x) v_j(x)\right) &= 0, \\
\pder{}{x_j}\left(\rho(x) v_i(x)v_j(x)+P(x)\delta_{ij}\right) &= 0, \\
\pder{}{x_j}\left(\left[\rho(x)E(x)+P(x)\right]v_j(x)\right) &= 0,
\end{split}
\end{equation}
for all $x\in\Omega$, where $i=1,\dots,d$ and summation is implied over
the repeated index $j=1,\dots,d$. The density 
$\func{\rho}{\Omega}{\Rbb_{>0}}$, velocity 
$\func{v_i}{\Omega}{\Rbb}$ in the $x_i$
direction for $i=1,\dots,d$, and total energy 
$\func{E}{\Omega}{\Rbb_{>0}}$ of the fluid are
implicitly defined as the solution of (\ref{eqn:euler}). 
We assume the fluid is an ideal gas that follows the 
thermal and caloric state equations 
\begin{equation} \label{eqn:idealgas}
P = \rho R T, \qquad
e = \frac{1}{\gamma-1}\frac{P}{\rho},
\end{equation}
where $\func{P}{\Omega}{\Rbb_{>0}}$ and $\func{T}{\Omega}{\Rbb_{>0}}$ 
are the pressure and temperature of the fluid, $\func{e}{\Omega}{\Rbb_{>0}}$ 
is the internal energy of the fluid, $R$ is the specific gas constant, and 
$\gamma$ is the ratio of specific heats. The frozen sound speed 
$\func{c}{\Omega}{\Rbb_{>0}}$ and Mach number $\func{M}{\Omega}{\Rbb_{>0}}$ 
are defined as
\begin{equation} \label{eqn:mach}
c^2 := \frac{\gamma P}{\rho}, \qquad
M := \frac{\sqrt{v_i v_i}}{c}.
\end{equation}
From the definition of the specific total energy $\func{E}{\Omega}{\Rbb_{>0}}$ 
and total enthalpy $\func{H}{\Omega}{\Rbb_{>0}}$
\begin{equation} \label{eqn:EandH}
E := e + \frac{v_i v_i}{2}, \qquad
H := \frac{\rho E + P}{\rho},
\end{equation}
the pressure of the fluid, $\func{P}{\Omega}{\Rbb_{>0}}$, is directly 
related to the conservative variables for an ideal gas as 
\begin{equation} \label{eqn:EtoP}
 P = (\gamma-1)\left(\rho E - \frac{\rho v_i v_i}{2}\right).
\end{equation}
The relationships between the downstream stagnation pressure $p_{02}$
and the upstream pressure $P_\infty$, and between the stagnation temperature 
$T_0$ and the upstream temperature $T_\infty$, are defined through the upstream
Mach number $M_\infty$ and $\gamma$ as 
\begin{equation}\label{eqn:tot_prestemp}
\frac{p_{02}}{P_\infty} = \frac{
\left[\frac{(\gamma+1)M_\infty^2}{2}\right]^{\gamma/(\gamma-1)}
}{
\left[
\left(\frac{2\gamma M_\infty^2}{\gamma+1}\right) - \left(\frac{\gamma-1}{\gamma+1}\right)
\right]^{1/(\gamma-1)}
}, \qquad
\frac{T_0}{T_\infty} = 1 + \frac{\gamma-1}{2} M_\infty^2.
\end{equation}
Finally, the Euler equations can be written as a general system of conservation laws
(\ref{eqn:claw-phys}) as
\begin{equation*}
 U(x) = \begin{bmatrix} \rho(x) \\ \rho(x) v(x) \\ \rho(x) E(x) \end{bmatrix}, \qquad
 F(U(x)) = \begin{bmatrix} \rho(x)v(x)^T \\ \rho v(x)v(x)^T + P(x) I \\ (\rho(x)E(x)+P(x)) v(x)^T\end{bmatrix}, \qquad
 S(U(x)) = 0.
\end{equation*}

\begin{remark}
We found that it is necessary to use boundary conditions on the shock boundary
that directly evaluate the physical flux at the boundary state, as opposed to
boundary conditions based on Riemann solvers. We found that Riemann solver boundary
conditions can fail to enforce the boundary state on the shock boundary and derail
the implicit tracking process.
\end{remark}

\subsubsection{Inviscid flow over cylinder, Mach continuation to $M_\infty = 3$}
\label{sec:numexp:mach2to3}
In this problem, we solve for steady $M_\infty=3$ flow over a cylinder
in $d=2$ dimensions using Mach continuation beginning at $M_\infty=2$.
At $M_\infty=3$, we have the following stagnation quantities, which will
be used to make quantitative assessment of the accuracy of the HOIST method:
$p_{02}=12.061$, $T_0=2$, and $H_0=7$. Our approach begins by
applying the HOIST method at $M_\infty = 2$ on a shock-agnostic mesh
of the entire domain consisting of $229$ elements and extracting the
solution and mesh downstream of the lead shock to initialize our
continuation strategy, resulting in a reduced mesh with $172$ elements
(Figure~\ref{fig:2d:mach2}).
%
%
%

We examine the performance of the framework with different number of 
continuation stages, $N = 2, 5, 10, 20$. The intermediate stages are
solved with at most $n_1=30$ iterations or to a tolerance of $\xi_1=10^{-4}$,
and the final stage is solved with at most $n_2=100$ iterations or to a
tolerance of $\xi_2=10^{-8}$. The HOIST solver parameters
used for all $\Upsilon$ evaluations are $\lambda=10^{-4}$, $\kappa_0=1$,
and $(\zeta, \upsilon)=(4, 0.25)$ \cite{huang2022robust}. 
For all four cases, the SQP solver drives the DG residual
$\norm{\rbm(\ubm_k^i,\phibold(\ybm_k^i),\mubold_k)}$ 
rapidly towards the early termination criteria for each
intermediate stage; in the final stage, the enriched DG residual
$\norm{\Rbm(\ubm_N^i,\phibold(\ybm_N^i),\mubold_N)}$ plateaus within
15 iterations (Figure~\ref{fig:mach2to3_study2_hist}). Even though the 
DG residual does not reach a tight tolerance for any of the case, all cases 
exhibits deep convergence, i.e.,
$\norm{\rbm(\ubm,\phibold(\ybm_N^I),\mubold_N)}\sim\Ocal(10^{-14})$,
within 3 Newton iterations using the HOIST output as its starting point.
Thus, the shock is indeed tracked and the fully discrete PDE is
satisfied on the final discontinuity-aligned mesh, which shows the
continuation framework is robust with as few as two stages.
\ifbool{fastcompile}{}{
\begin{figure}[!htbp]
\centering
\input{./_py/mach2to3_study2.tikz}
\caption{
 SQP convergence history of the DG residual
 $\norm{\rbm(\ubm_k^i,\phibold(\ybm_k^i);\mubold_k)}$ (\ref{line:R0})
 and enriched DG residual $\norm{\Rbm(\ubm_k^i,\phibold(\ybm_k^i);\mubold_k)}$ 
 (\ref{line:R1err}), and the DG residual throughout the fixed-mesh
 Newton iterations
 $\norm{\rbm(\ubm_i,\phibold(\ybm_k^I);\mubold_k)}$ (\ref{line:R0_fix})
 with $N = 2, 5, 10, 20$ (\emph{top to bottom}) continuation stages
 for $M_\infty=2$ to $M_\infty=3$ continuation (cylinder).
 }
\label{fig:mach2to3_study2_hist}
\end{figure}
}

\begin{figure}[!htbp]
\centering
\ifbool{fastcompile}{}{
\begin{tikzpicture}
\begin{groupplot}[
  group style={
      group size=6 by 2,
      horizontal sep=0.1cm,
      vertical sep=0.1cm
  },
  width=0.7\textwidth,
  axis equal image,
  xlabel={},
  ylabel={},
  xtick = {-1.70},
  ytick = {-3.44, 3.44},
  xmin=-2.3, xmax=0,
  ymin=-5.65, ymax=5.65,
  yticklabel pos=right,
  xticklabel style={/pgf/number format/.cd,fixed zerofill,precision=2}
]

\nextgroupplot[yticklabel pos=left, yticklabels={,,}, xticklabels={,,}]
\addplot graphics [xmin=-2.3, xmax=0, ymin=-5.65, ymax=5.65] {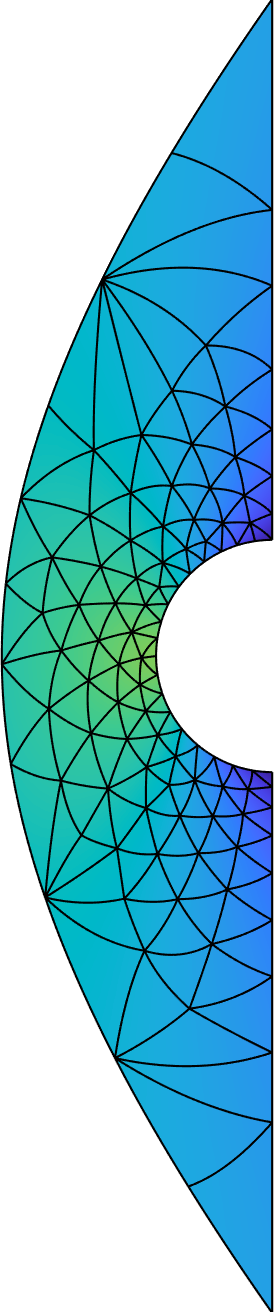};
\nextgroupplot[ylabel={}, yticklabels={,,}, xticklabels={,,}]
\addplot graphics [xmin=-2.3, xmax=0, ymin=-5.65, ymax=5.65] {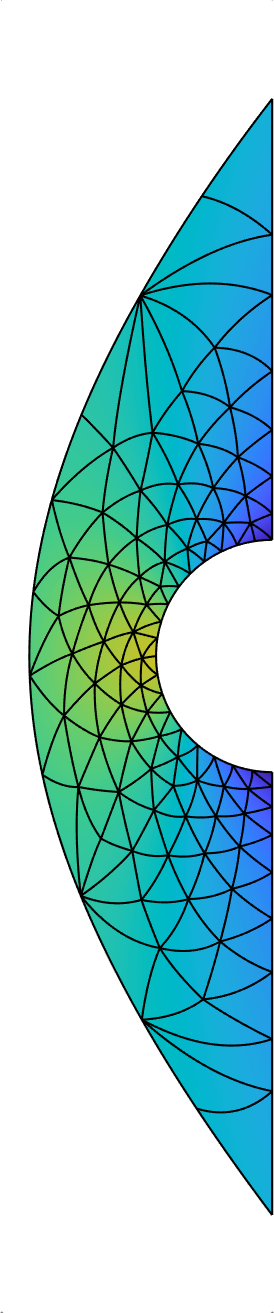};
\nextgroupplot[ylabel={}, yticklabels={,,}, xticklabels={,,}]
\addplot graphics [xmin=-2.3, xmax=0, ymin=-5.65, ymax=5.65] {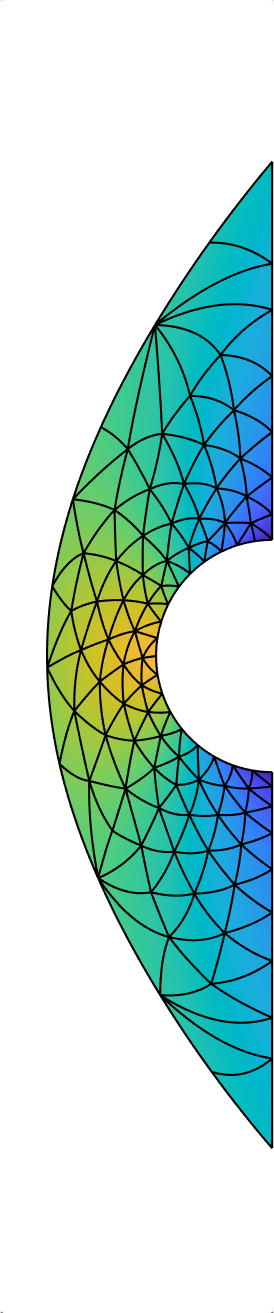};
\nextgroupplot[ylabel={}, yticklabels={,,}, xticklabels={,,}]
\addplot graphics [xmin=-2.3, xmax=0, ymin=-5.65, ymax=5.65] {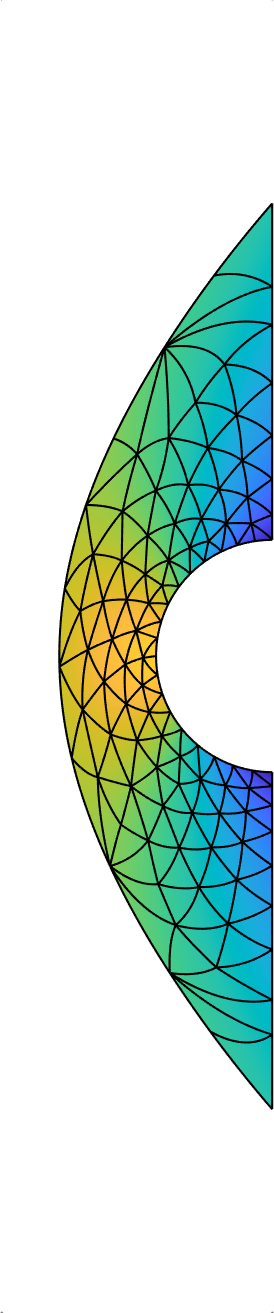};
\nextgroupplot[ylabel={}, yticklabels={,,}, xticklabels={,,}]
\addplot graphics [xmin=-2.3, xmax=0, ymin=-5.65, ymax=5.65] {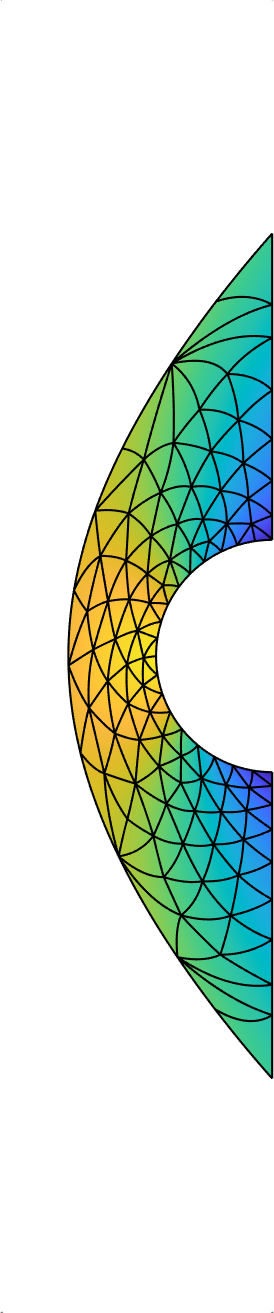};
\nextgroupplot[ylabel={}, xticklabel pos=top]
\addplot graphics [xmin=-2.3, xmax=0, ymin=-5.65, ymax=5.65] {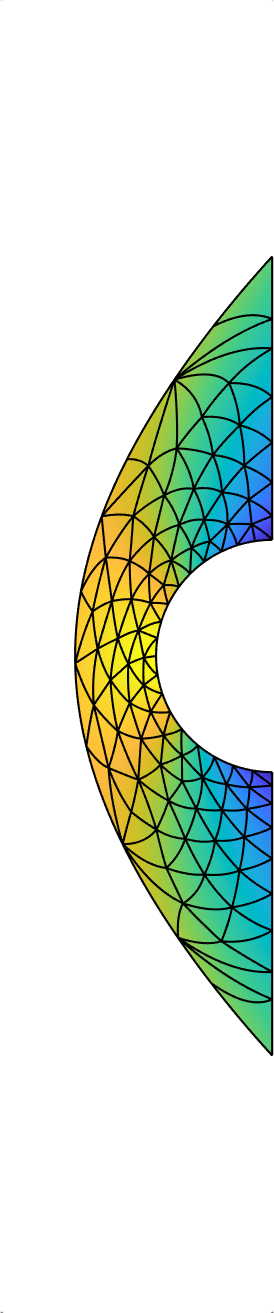};
\addplot [solid, gray, dashed]
coordinates {
( -1.7,  0)
( -1.7,  5.65)};

\nextgroupplot[yticklabel pos=left, yticklabels={,,}, xticklabels={,,}]
\addplot graphics [xmin=-2.3, xmax=0, ymin=-5.65, ymax=5.65] {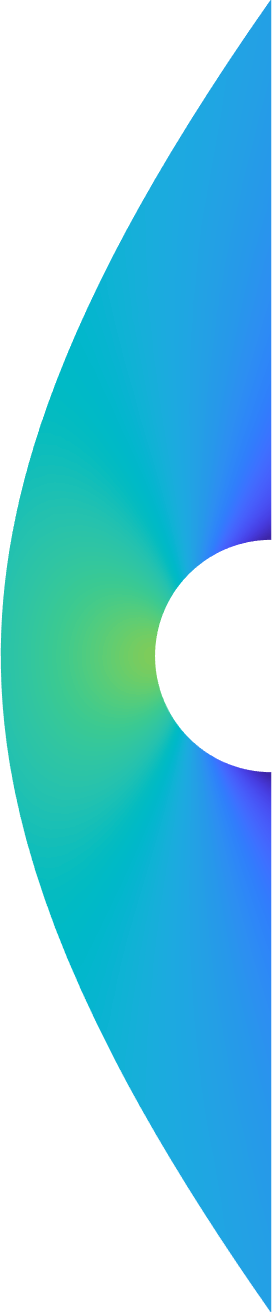};
\nextgroupplot[ylabel={}, yticklabels={,,}, xticklabels={,,}]
\addplot graphics [xmin=-2.3, xmax=0, ymin=-5.65, ymax=5.65] {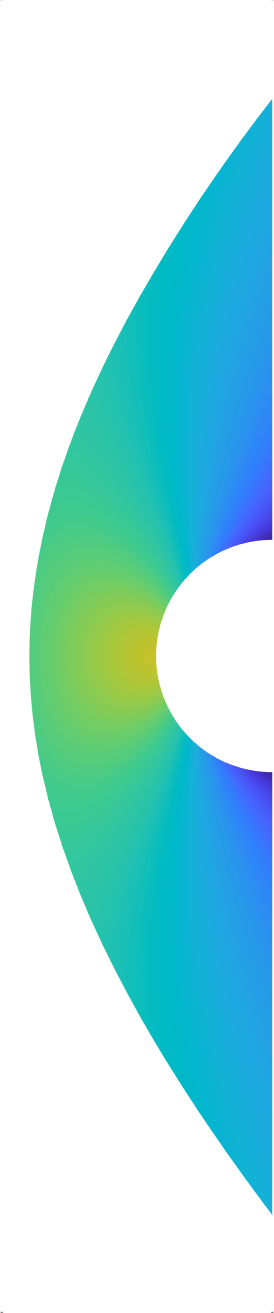};
\nextgroupplot[ylabel={}, yticklabels={,,}, xticklabels={,,}]
\addplot graphics [xmin=-2.3, xmax=0, ymin=-5.65, ymax=5.65] {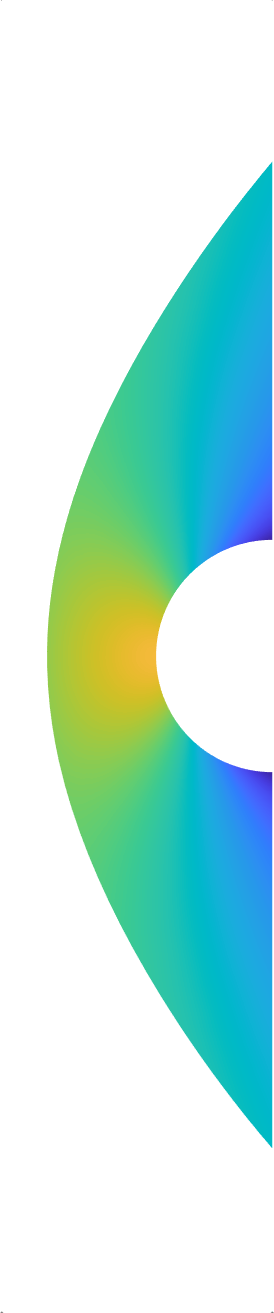};
\nextgroupplot[ylabel={}, yticklabels={,,}, xticklabels={,,}]
\addplot graphics [xmin=-2.3, xmax=0, ymin=-5.65, ymax=5.65] {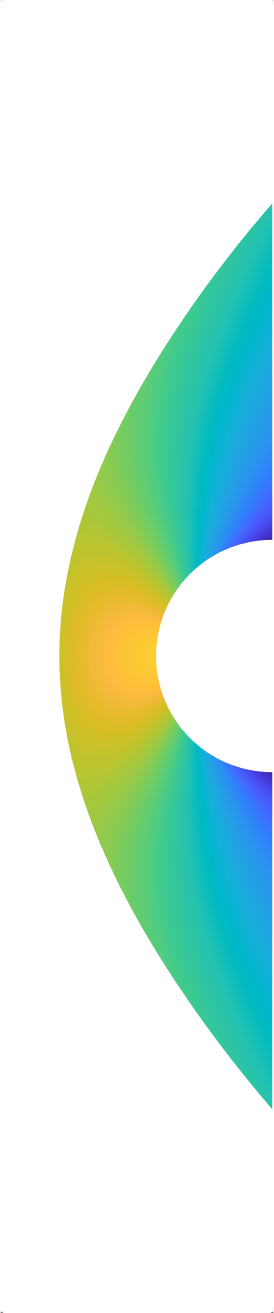};
\nextgroupplot[ylabel={}, yticklabels={,,}, xticklabels={,,}]
\addplot graphics [xmin=-2.3, xmax=0, ymin=-5.65, ymax=5.65] {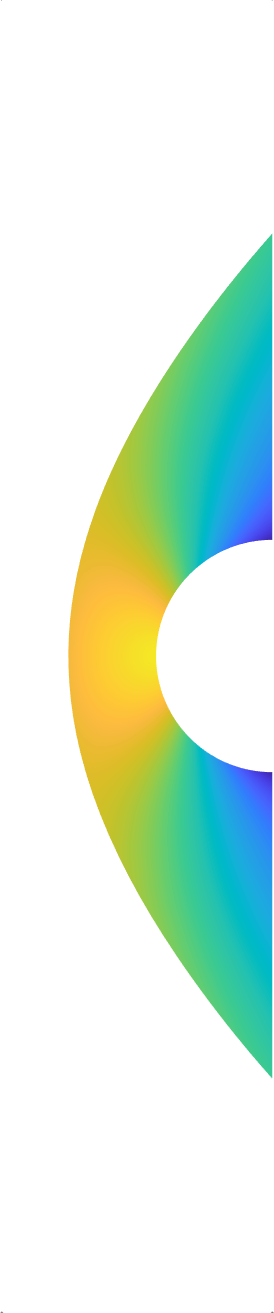};
\nextgroupplot[ylabel={}, xticklabels={,,}, yticklabels={,,}]
\addplot graphics [xmin=-2.3, xmax=0, ymin=-5.65, ymax=5.65] {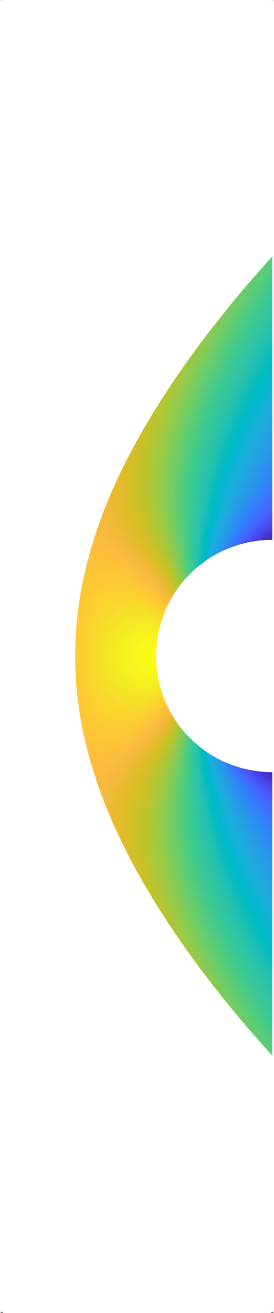};

\end{groupplot}
\end{tikzpicture}

\colorbarMatlabParula{1.14}{2}{3.5}{5}{6}
}
 \caption{Density distribution at all continuation stages
  ($M_\infty = 2, 2.2, 2.4, 2.6, 2.8, 3$) for $N=5$
  (\emph{left} to \emph{right}).}
\label{fig:2d:mach2to3}
\end{figure}

The proposed framework produces well-resolved, accurate shock profiles and
temperature and pressure distributions along the cylinder
(Figure~\ref{fig:mach2to3_study2_shkpos_temppres}). The HOIST solver
is robust in that it does not require parameter tuning across stages despite
the changing Mach number. The value and the relative error of stagnation
quantities evaluated at the stagnation point $(x_1, x_2) = (-1, 0)$ are
reported in Table~\ref{tab:stag_err}, which are all quite small, particularly
considering the coarse grid used (only $172$ quadratic elements).
For the $N=5$ case, Figure~\ref{fig:2d:mach2to3} shows the flow solution
and mesh at all Mach numbers considered, including intermediate (partially
converged) ones. In all cases, the shock is tracked with curved, high-order
elements and the solution is well-resolved throughout the domain. The overall 
mesh quality is also well-preserved as the shock boundary is 
compressed towards the cylinder.


\ifbool{fastcompile}{}{
\begin{figure}[!htbp]
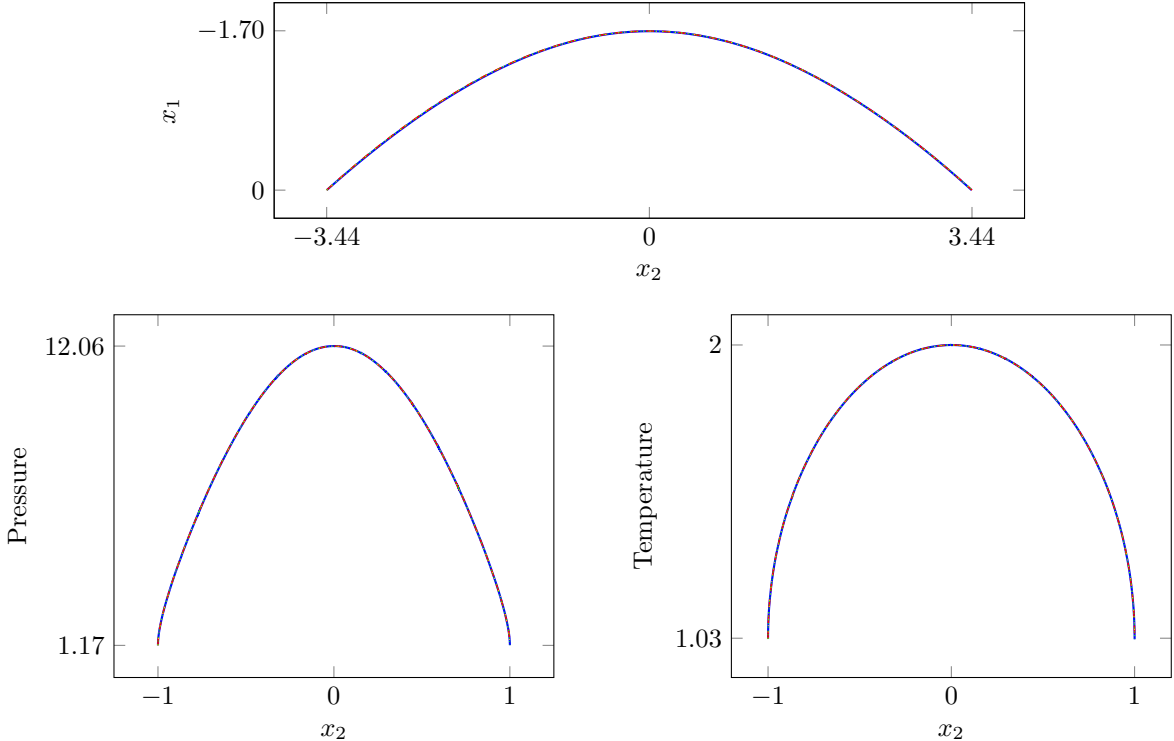

\centering
\input{./_py/mach2to3_study2_shkpos.tikz} \\\vspace{3mm}
\input{./_py/mach2to3_study2_pres.tikz} \qquad
\input{./_py/mach2to3_study2_temp.tikz}
\caption{
 Lead shock position (\textit{top}), and pressure (\textit{bottom left})
 and temperature (\textit{bottom right}) along the surface of the cylinder
 at $M_\infty=3$ using Mach continuation starting from $M_\infty=2$ with
 $N = 2$ (\ref{pnts:study2_2s_shk}), $N = 5$ (\ref{pnts:study2_5s_shk}),
 $N = 10$ (\ref{pnts:study2_10s_shk}), and $N = 20$ (\ref{pnts:study2_20s_shk})
 stages.
 }
\label{fig:mach2to3_study2_shkpos_temppres}
\end{figure}
}

\subsubsection{Inviscid flow over cylinder, Mach continuation to $M_\infty = 10$}
\label{sec:numexp:mach2to10}
Next, we increase the difficulty of the previous problem by increasing the
target Mach number to $M_\infty=10$ using Mach continuation beginning at
$M_\infty=2$. In this case, the stagnation quantities
are $p_{02}=129.217$, $T_0=15$, and $H_0=52.5$. We use the same initial
configuration as in Section~\ref{sec:numexp:mach2to3} (after applying HOIST
to the full domain at $M_\infty=2$ and removing elements upstream of the
bow shock), and we take the number of continuation stages to be $N = 40$
(5 stages per Mach number). The intermediate stages are
solved with at most $n_1=30$ iterations or to a tolerance of $\xi_1=10^{-4}$,
and the final stage is solved with at most $n_2=100$ iterations or to a
tolerance of $\xi_2=10^{-8}$. We use the same HOIST
solver parameters as in Section~\ref{sec:numexp:mach2to3} and they remain
the same throughout the optimization. 
In the final stage (from Mach 9.8 to Mach 10), the enriched DG residual
$\norm{\Rbm(\ubm_N^i,\phibold(\ybm_N^i),\mubold_N)}$ plateaus within 15
iterations (Figure~\ref{fig:mach2to10_hist}) and the DG residual
reaches a tolerance of $10^{-14}$ with only $4$ fixed-mesh Newton iterations.
\ifbool{fastcompile}{}{
\begin{figure}[!htbp]
\centering
\begin{tikzpicture}
\begin{groupplot} [
group style={group size = 3 by 1, horizontal sep = 1.4cm, vertical sep = 0.8cm}]
\nextgroupplot[xtick={0,50,100}, xlabel={Iteration ($i$)}, ytick={1e1, 1e-2}, xmin=-5, xmax=105, ymin=0.005, ymax=20.0, ymode=log, width=0.315\textwidth]
\addplot [black, thin, mark=*, mark size=0.75, mark options={solid, thin}, mark repeat=2]
coordinates {
( 1.00000000e+00,  9.83311785e+00)
( 2.00000000e+00,  4.39703806e+00)
( 3.00000000e+00,  2.26599546e+00)
( 4.00000000e+00,  1.28979406e+00)
( 5.00000000e+00,  4.62317185e-01)
( 6.00000000e+00,  3.88257923e-01)
( 7.00000000e+00,  3.15178624e-01)
( 8.00000000e+00,  2.32767461e-01)
( 9.00000000e+00,  1.51907419e-01)
( 1.00000000e+01,  1.28716574e-01)
( 1.10000000e+01,  1.19441694e-01)
( 1.20000000e+01,  1.19119464e-01)
( 1.30000000e+01,  1.04528708e-01)
( 1.40000000e+01,  1.16609753e-01)
( 1.50000000e+01,  1.18233625e-01)
( 1.60000000e+01,  1.03630723e-01)
( 1.70000000e+01,  1.64071377e-02)
( 1.80000000e+01,  3.01415629e-02)
( 1.90000000e+01,  2.64649321e-02)
( 2.00000000e+01,  2.60490747e-02)
( 2.10000000e+01,  2.19964937e-02)
( 2.20000000e+01,  1.93937613e-02)
( 2.30000000e+01,  1.81808273e-02)
( 2.40000000e+01,  1.07383841e-02)
( 2.50000000e+01,  1.46350622e-02)
( 2.60000000e+01,  1.44215189e-02)
( 2.70000000e+01,  1.41964126e-02)
( 2.80000000e+01,  1.11165054e-02)
( 2.90000000e+01,  1.14896348e-02)};\label{line:R0}

\addplot [black, thin, mark=*, mark size=0.75, mark options={solid, thin}, mark repeat=10]
coordinates {
( 3.00000000e+01,  9.54637524e-03)
( 3.10000000e+01,  9.30731049e-03)
( 3.20000000e+01,  9.33730850e-03)
( 3.30000000e+01,  9.30084848e-03)
( 3.40000000e+01,  9.39567770e-03)
( 3.50000000e+01,  9.56209280e-03)
( 3.60000000e+01,  9.63395572e-03)
( 3.70000000e+01,  9.68441478e-03)
( 3.80000000e+01,  1.00233679e-02)
( 3.90000000e+01,  1.00923928e-02)
( 4.00000000e+01,  1.00160975e-02)
( 4.10000000e+01,  1.00112070e-02)
( 4.20000000e+01,  9.96523197e-03)
( 4.30000000e+01,  9.99517040e-03)
( 4.40000000e+01,  9.94381357e-03)
( 4.50000000e+01,  1.00225000e-02)
( 4.60000000e+01,  1.02180876e-02)
( 4.70000000e+01,  1.02259195e-02)
( 4.80000000e+01,  1.01882345e-02)
( 4.90000000e+01,  1.00925646e-02)
( 5.00000000e+01,  1.01759201e-02)
( 5.10000000e+01,  1.02220047e-02)
( 5.20000000e+01,  1.01820871e-02)
( 5.30000000e+01,  1.00702420e-02)
( 5.40000000e+01,  1.01106197e-02)
( 5.50000000e+01,  1.01969727e-02)
( 5.60000000e+01,  1.01870154e-02)
( 5.70000000e+01,  1.01673329e-02)
( 5.80000000e+01,  1.02655609e-02)
( 5.90000000e+01,  1.02555366e-02)
( 6.00000000e+01,  1.02530329e-02)
( 6.10000000e+01,  1.02120793e-02)
( 6.20000000e+01,  1.02347270e-02)
( 6.30000000e+01,  1.02030969e-02)
( 6.40000000e+01,  1.02690279e-02)
( 6.50000000e+01,  1.02640139e-02)
( 6.60000000e+01,  1.02242061e-02)
( 6.70000000e+01,  1.02545168e-02)
( 6.80000000e+01,  1.02423076e-02)
( 6.90000000e+01,  1.02609726e-02)
( 7.00000000e+01,  1.02979459e-02)
( 7.10000000e+01,  1.02915413e-02)
( 7.20000000e+01,  1.02856194e-02)
( 7.30000000e+01,  1.02814498e-02)
( 7.40000000e+01,  1.02760390e-02)
( 7.50000000e+01,  1.02710216e-02)
( 7.60000000e+01,  1.02724843e-02)
( 7.70000000e+01,  1.02726721e-02)
( 7.80000000e+01,  1.02870385e-02)
( 7.90000000e+01,  1.02845270e-02)
( 8.00000000e+01,  1.02798865e-02)
( 8.10000000e+01,  1.02793692e-02)
( 8.20000000e+01,  1.02848964e-02)
( 8.30000000e+01,  1.02842785e-02)
( 8.40000000e+01,  1.02841740e-02)
( 8.50000000e+01,  1.02841348e-02)
( 8.60000000e+01,  1.02841250e-02)
( 8.70000000e+01,  1.02826342e-02)
( 8.80000000e+01,  1.02810669e-02)
( 8.90000000e+01,  1.02811108e-02)
( 9.00000000e+01,  1.02795356e-02)
( 9.10000000e+01,  1.02796823e-02)
( 9.20000000e+01,  1.02791965e-02)
( 9.30000000e+01,  1.02791752e-02)
( 9.40000000e+01,  1.02791556e-02)
( 9.50000000e+01,  1.02787840e-02)
( 9.60000000e+01,  1.02791059e-02)
( 9.70000000e+01,  1.02788204e-02)
( 9.80000000e+01,  1.02788734e-02)
( 9.90000000e+01,  1.02788239e-02)};\label{line:R0}

\nextgroupplot[xtick={0,50,100}, xlabel={Iteration ($i$)}, ytick={1e1, 1e0}, xmin=-5, xmax=105, ymin=1.0, ymax=10.0, ymode=log, width=0.315\textwidth]
\addplot [blue, thin, mark=square*, mark size=0.75, mark options={solid, thin}, mark repeat=3]
coordinates {
( 1.00000000e+00,  8.43457410e+00)
( 2.00000000e+00,  4.65086980e+00)
( 3.00000000e+00,  2.39976988e+00)
( 4.00000000e+00,  1.68013876e+00)
( 5.00000000e+00,  1.29720137e+00)
( 6.00000000e+00,  1.28943872e+00)
( 7.00000000e+00,  1.24858928e+00)
( 8.00000000e+00,  1.23391421e+00)
( 9.00000000e+00,  1.21292353e+00)
( 1.00000000e+01,  1.21269742e+00)
( 1.10000000e+01,  1.21157393e+00)
( 1.20000000e+01,  1.21044145e+00)
( 1.30000000e+01,  1.20914828e+00)
( 1.40000000e+01,  1.21436915e+00)};\label{line:R1err}

\addplot [blue, thin, mark=square*, mark size=0.75, mark options={solid, thin}, mark repeat=10]
coordinates {
( 1.50000000e+01,  1.21322440e+00)
( 1.60000000e+01,  1.21866774e+00)
( 1.70000000e+01,  1.21529042e+00)
( 1.80000000e+01,  1.21497595e+00)
( 1.90000000e+01,  1.21512804e+00)
( 2.00000000e+01,  1.21561293e+00)
( 2.10000000e+01,  1.21501345e+00)
( 2.20000000e+01,  1.21538829e+00)
( 2.30000000e+01,  1.21588056e+00)
( 2.40000000e+01,  1.21447318e+00)
( 2.50000000e+01,  1.21414339e+00)
( 2.60000000e+01,  1.21418199e+00)
( 2.70000000e+01,  1.21458867e+00)
( 2.80000000e+01,  1.21429168e+00)
( 2.90000000e+01,  1.21418247e+00)
( 3.00000000e+01,  1.21429154e+00)
( 3.10000000e+01,  1.21418507e+00)
( 3.20000000e+01,  1.21411836e+00)
( 3.30000000e+01,  1.21413654e+00)
( 3.40000000e+01,  1.21417112e+00)
( 3.50000000e+01,  1.21410064e+00)
( 3.60000000e+01,  1.21411393e+00)
( 3.70000000e+01,  1.21419495e+00)
( 3.80000000e+01,  1.21404953e+00)
( 3.90000000e+01,  1.21400047e+00)
( 4.00000000e+01,  1.21399108e+00)
( 4.10000000e+01,  1.21399196e+00)
( 4.20000000e+01,  1.21406188e+00)
( 4.30000000e+01,  1.21401432e+00)
( 4.40000000e+01,  1.21404375e+00)
( 4.50000000e+01,  1.21406792e+00)
( 4.60000000e+01,  1.21399900e+00)
( 4.70000000e+01,  1.21397816e+00)
( 4.80000000e+01,  1.21397733e+00)
( 4.90000000e+01,  1.21404566e+00)
( 5.00000000e+01,  1.21399875e+00)
( 5.10000000e+01,  1.21398317e+00)
( 5.20000000e+01,  1.21398473e+00)
( 5.30000000e+01,  1.21401495e+00)
( 5.40000000e+01,  1.21399792e+00)
( 5.50000000e+01,  1.21397732e+00)
( 5.60000000e+01,  1.21398104e+00)
( 5.70000000e+01,  1.21399200e+00)
( 5.80000000e+01,  1.21396095e+00)
( 5.90000000e+01,  1.21395983e+00)
( 6.00000000e+01,  1.21395967e+00)
( 6.10000000e+01,  1.21399630e+00)
( 6.20000000e+01,  1.21397288e+00)
( 6.30000000e+01,  1.21398468e+00)
( 6.40000000e+01,  1.21395930e+00)
( 6.50000000e+01,  1.21396046e+00)
( 6.60000000e+01,  1.21397715e+00)
( 6.70000000e+01,  1.21396558e+00)
( 6.80000000e+01,  1.21397270e+00)
( 6.90000000e+01,  1.21397957e+00)
( 7.00000000e+01,  1.21396265e+00)
( 7.10000000e+01,  1.21396041e+00)
( 7.20000000e+01,  1.21395791e+00)
( 7.30000000e+01,  1.21396251e+00)
( 7.40000000e+01,  1.21395991e+00)
( 7.50000000e+01,  1.21396119e+00)
( 7.60000000e+01,  1.21396856e+00)
( 7.70000000e+01,  1.21396281e+00)
( 7.80000000e+01,  1.21395861e+00)
( 7.90000000e+01,  1.21395941e+00)
( 8.00000000e+01,  1.21396294e+00)
( 8.10000000e+01,  1.21396034e+00)
( 8.20000000e+01,  1.21395835e+00)
( 8.30000000e+01,  1.21395827e+00)
( 8.40000000e+01,  1.21395825e+00)
( 8.50000000e+01,  1.21395824e+00)
( 8.60000000e+01,  1.21395825e+00)
( 8.70000000e+01,  1.21395926e+00)
( 8.80000000e+01,  1.21395898e+00)
( 8.90000000e+01,  1.21395918e+00)
( 9.00000000e+01,  1.21395982e+00)
( 9.10000000e+01,  1.21395917e+00)
( 9.20000000e+01,  1.21395909e+00)
( 9.30000000e+01,  1.21395906e+00)
( 9.40000000e+01,  1.21395906e+00)
( 9.50000000e+01,  1.21395918e+00)
( 9.60000000e+01,  1.21395906e+00)
( 9.70000000e+01,  1.21395909e+00)
( 9.80000000e+01,  1.21395909e+00)
( 9.90000000e+01,  1.21395905e+00)};\label{line:R1err}

\nextgroupplot[xtick={0,2,4}, xlabel={Iteration ($i$)}, xmin=-0.25, xmax=4.25, ymin=1e-14, ymax=0.051394119583632776, ymode=log, width=0.315\textwidth]
\addplot [red, thin, mark=triangle*, mark size=1, mark options={solid, thin}, mark repeat=0]
coordinates {
( 0.00000000e+00,  1.02788239e-02)
( 1.00000000e+00,  8.69413055e-03)
( 2.00000000e+00,  1.12945576e-05)
( 3.00000000e+00,  5.78216189e-12)
( 4.00000000e+00,  9.26127042e-13)};\label{line:R0_fix}

\end{groupplot}\end{tikzpicture}
\caption{
 SQP convergence history of the DG residual
 $\norm{\rbm(\ubm_k^i,\phibold(\ybm_k^i),\mubold_k)}$ (\ref{line:R0})
 and enriched DG residual $\norm{\Rbm(\ubm_k^i,\phibold(\ybm_k^i),\mubold_k)}$
 (\ref{line:R1err}), and the DG residual throughout the fixed-mesh
 Newton iterations
 $\norm{\rbm(\ubm_i,\phibold(\ybm_k^I),\mubold_k)}$ for the
 final stage from Mach 9.8 to 10 (\ref{line:R0_fix}) for $M_\infty=2$
 to $M_\infty=10$ continuation (cylinder).
 }
\label{fig:mach2to10_hist}
\end{figure}
}

\begin{figure}[!htbp]
\centering
\ifbool{fastcompile}{}{
\begin{tikzpicture}
\begin{groupplot}[
  group style={
      group size=9 by 2,
      horizontal sep=0.1cm,
      vertical sep=0.1cm
  },
  width=0.65\textwidth,
  axis equal image,
  xlabel={},
  ylabel={},
  xtick = {-1.40},
  ytick = {-2.40, 2.40},
  xmin=-2.3, xmax=0,
  ymin=-5.65, ymax=5.65,
  yticklabel pos=right,
  xticklabel style={/pgf/number format/.cd,fixed zerofill,precision=2},
  yticklabel style={/pgf/number format/.cd,fixed zerofill,precision=2}
]

\nextgroupplot[yticklabel pos=left, yticklabels={,,}, xticklabels={,,}]
\addplot graphics [xmin=-2.3, xmax=0, ymin=-5.65, ymax=5.65] {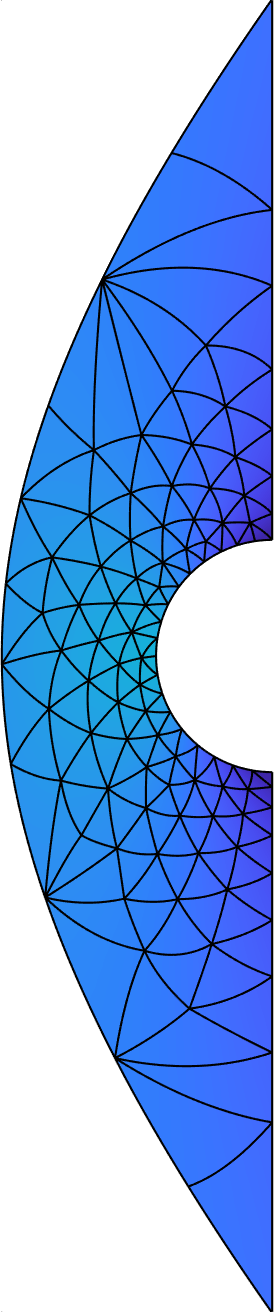};
\nextgroupplot[ylabel={}, yticklabels={,,}, xticklabels={,,}]
\addplot graphics [xmin=-2.3, xmax=0, ymin=-5.65, ymax=5.65] {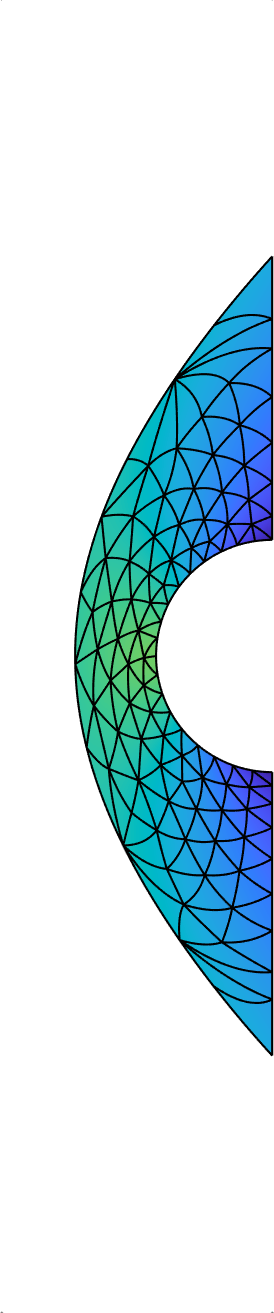};
\nextgroupplot[ylabel={}, yticklabels={,,}, xticklabels={,,}]
\addplot graphics [xmin=-2.3, xmax=0, ymin=-5.65, ymax=5.65] {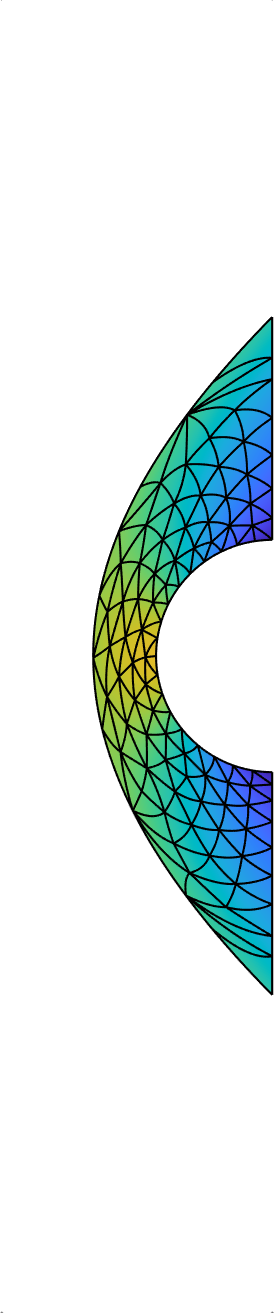};
\nextgroupplot[ylabel={}, yticklabels={,,}, xticklabels={,,}]
\addplot graphics [xmin=-2.3, xmax=0, ymin=-5.65, ymax=5.65] {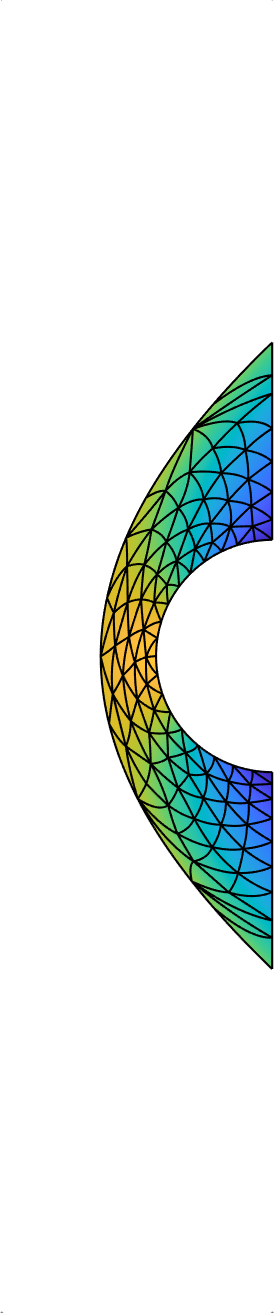};
\nextgroupplot[ylabel={}, yticklabels={,,}, xticklabels={,,}]
\addplot graphics [xmin=-2.3, xmax=0, ymin=-5.65, ymax=5.65] {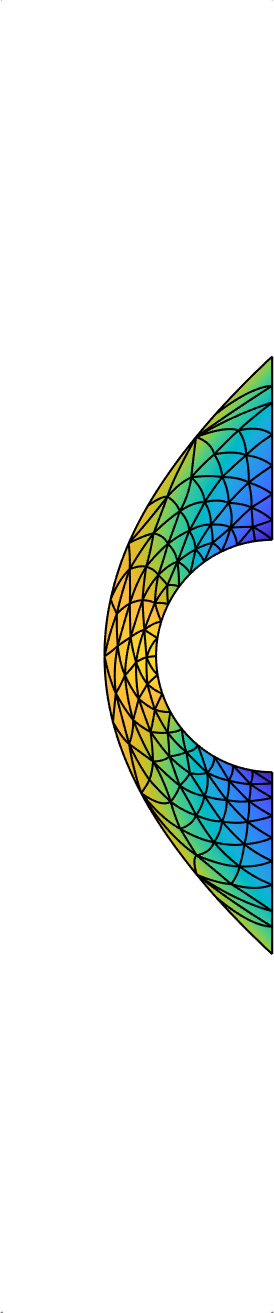};
\nextgroupplot[ylabel={}, yticklabels={,,}, xticklabels={,,}]
\addplot graphics [xmin=-2.3, xmax=0, ymin=-5.65, ymax=5.65] {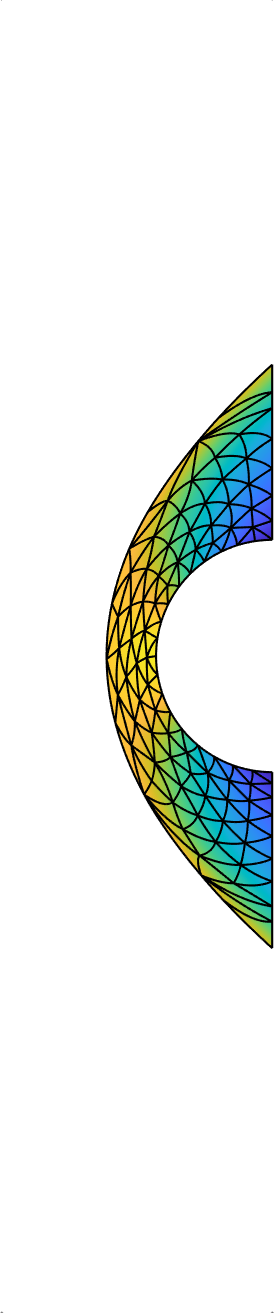};
\nextgroupplot[ylabel={}, yticklabels={,,}, xticklabels={,,}]
\addplot graphics [xmin=-2.3, xmax=0, ymin=-5.65, ymax=5.65] {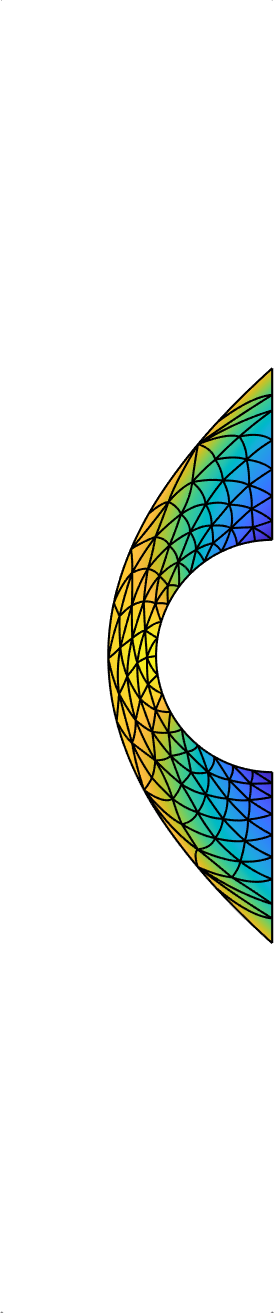};
\nextgroupplot[ylabel={}, yticklabels={,,}, xticklabels={,,}]
\addplot graphics [xmin=-2.3, xmax=0, ymin=-5.65, ymax=5.65] {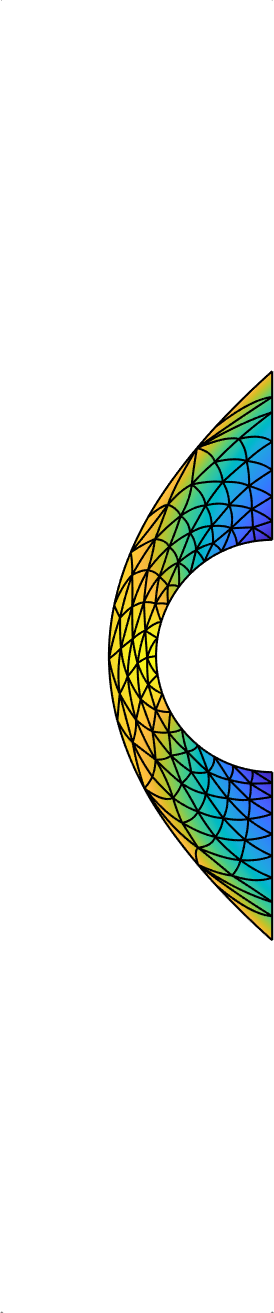};
\nextgroupplot[ylabel={}, xticklabel pos=top]
\addplot graphics [xmin=-2.3, xmax=0, ymin=-5.65, ymax=5.65] {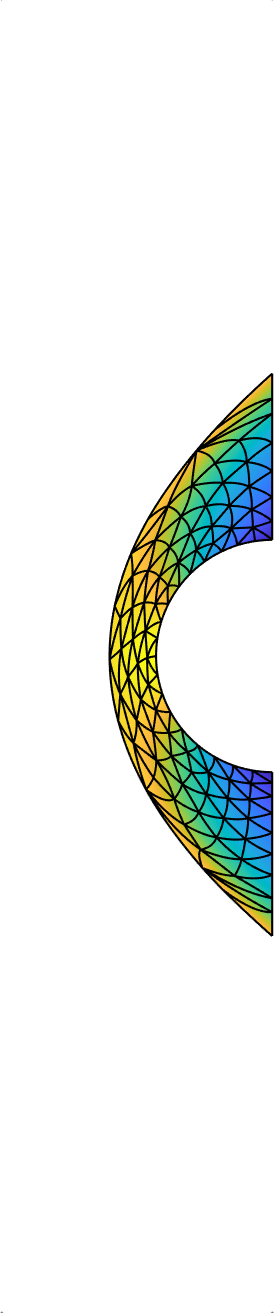};
\addplot [solid, gray, dashed]
coordinates {
( -1.4,  0)
( -1.4,  5.65)};

\nextgroupplot[yticklabel pos=left, yticklabels={,,}, xticklabels={,,}]
\addplot graphics [xmin=-2.3, xmax=0, ymin=-5.65, ymax=5.65] {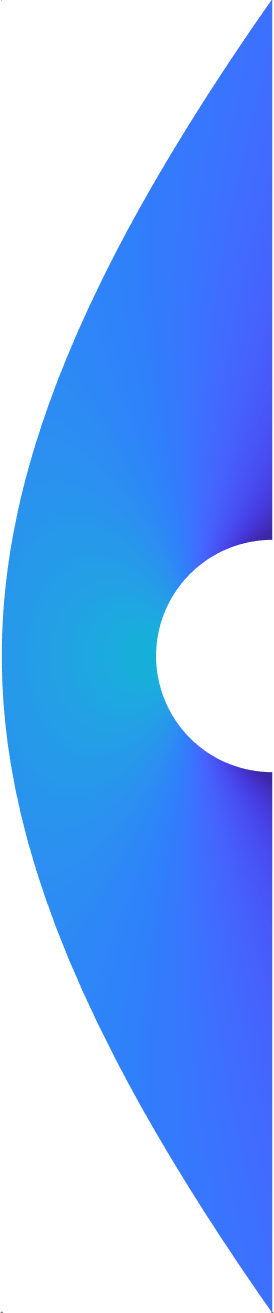};
\nextgroupplot[ylabel={}, yticklabels={,,}, xticklabels={,,}]
\addplot graphics [xmin=-2.3, xmax=0, ymin=-5.65, ymax=5.65] {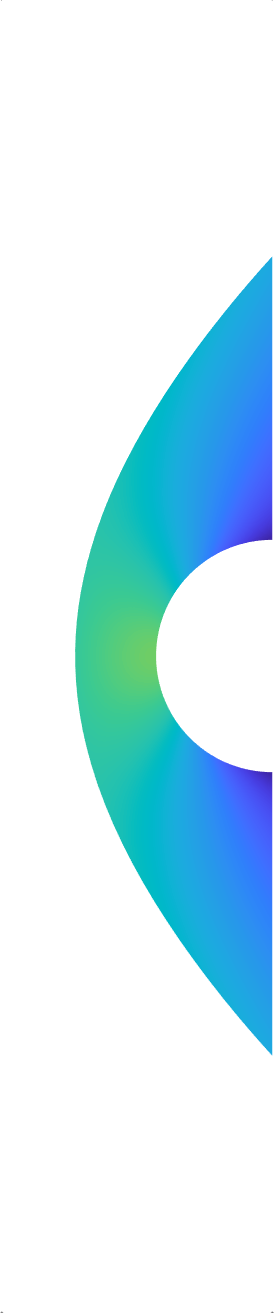};
\nextgroupplot[ylabel={}, yticklabels={,,}, xticklabels={,,}]
\addplot graphics [xmin=-2.3, xmax=0, ymin=-5.65, ymax=5.65] {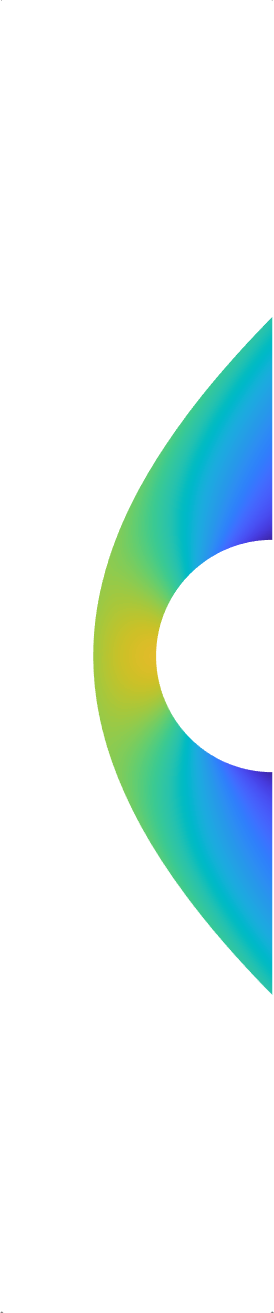};
\nextgroupplot[ylabel={}, yticklabels={,,}, xticklabels={,,}]
\addplot graphics [xmin=-2.3, xmax=0, ymin=-5.65, ymax=5.65] {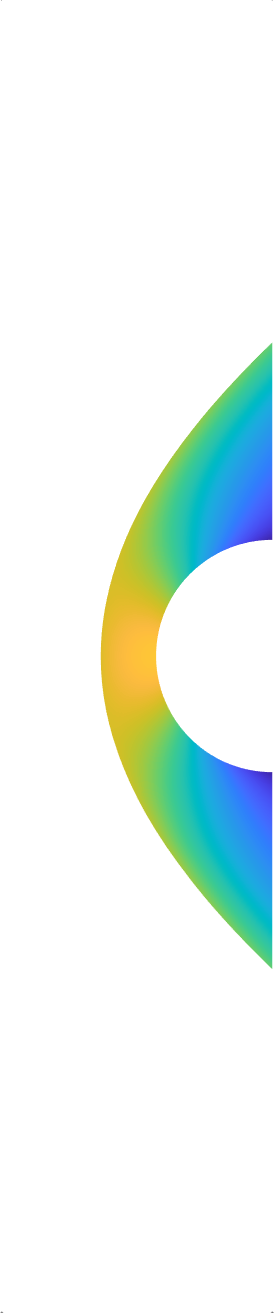};
\nextgroupplot[ylabel={}, yticklabels={,,}, xticklabels={,,}]
\addplot graphics [xmin=-2.3, xmax=0, ymin=-5.65, ymax=5.65] {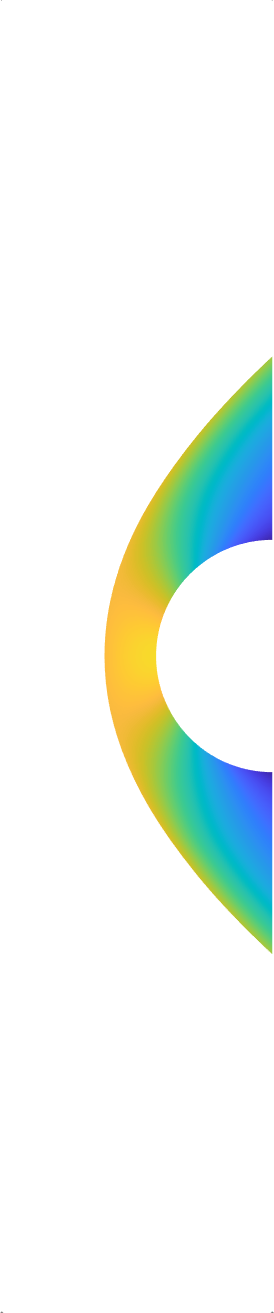};
\nextgroupplot[ylabel={}, yticklabels={,,}, xticklabels={,,}]
\addplot graphics [xmin=-2.3, xmax=0, ymin=-5.65, ymax=5.65] {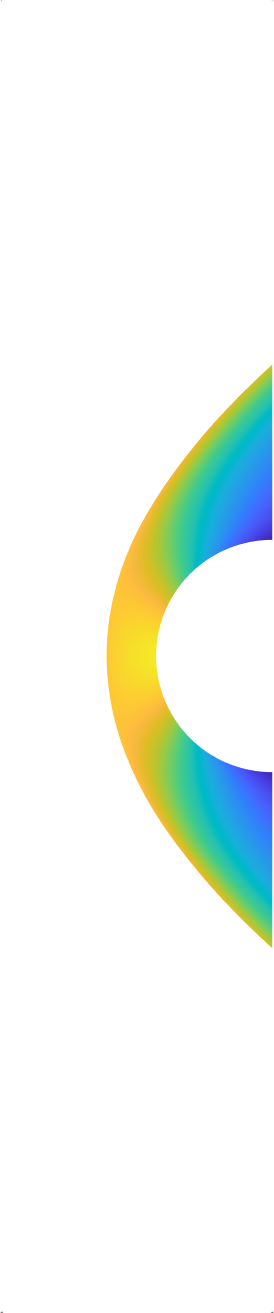};
\nextgroupplot[ylabel={}, yticklabels={,,}, xticklabels={,,}]
\addplot graphics [xmin=-2.3, xmax=0, ymin=-5.65, ymax=5.65] {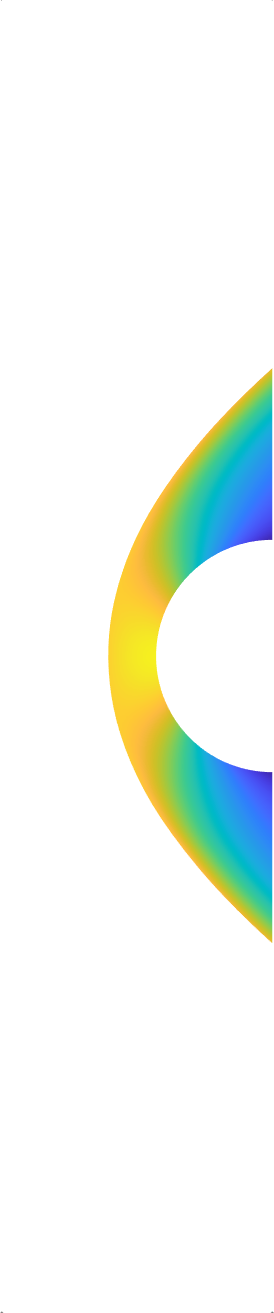};
\nextgroupplot[ylabel={}, yticklabels={,,}, xticklabels={,,}]
\addplot graphics [xmin=-2.3, xmax=0, ymin=-5.65, ymax=5.65] {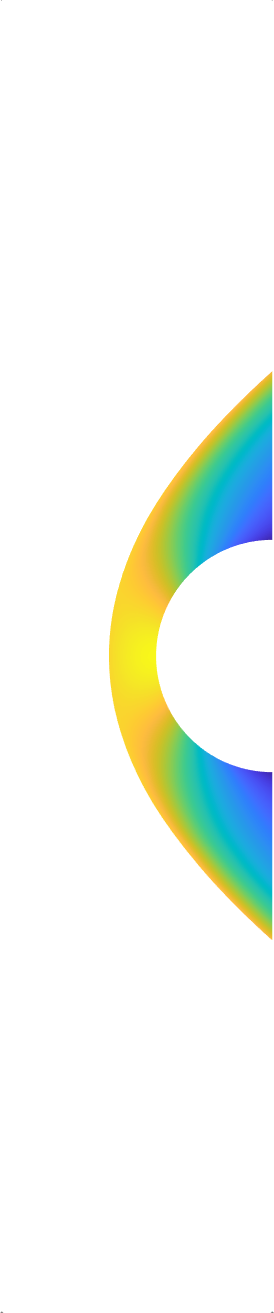};
\nextgroupplot[ylabel={}, xticklabels={,,}, yticklabels={,,}]
\addplot graphics [xmin=-2.3, xmax=0, ymin=-5.65, ymax=5.65] {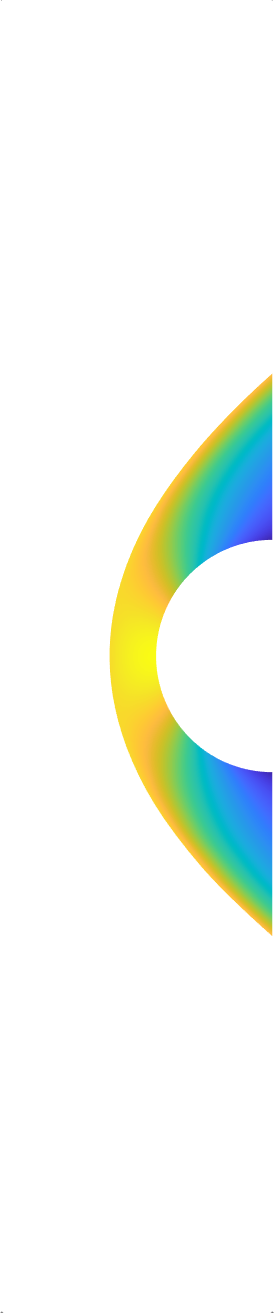};

\end{groupplot}
\end{tikzpicture}

\colorbarMatlabParula{1.2}{3}{5}{7}{8.6}
}
 \caption{Density distribution at selected continuation stages
  ($M_\infty = 2, 3, 4, \dots, 10$) (\emph{left} to \emph{right}).}
\label{fig:2d:mach2to10}
\end{figure}

The proposed framework accurately tracks the shock and the fully discrete PDE
is satisfied to a tight tolerance on the final discontinuity-aligned mesh. The
framework robustly handled the large Mach variation that causes the shock
to move significantly from its position at the initial Mach number $M_\infty=2$
to the target Mach number $M_\infty=10$
(Figures~\ref{fig:mach2to3_study2_shkpos_temppres}-\ref{fig:2d:mach2to10}).
The peak pressure and temperature along the surface are also significantly
larger than the original $M_\infty=2$ case, but are smooth and well-resolved
nonetheless. The value and relative error of stagnation quantities evaluated
at the stagnation point $(x_1, x_2) = (-1, 0)$ are reported in 
Table~\ref{tab:stag_err}, which are still quite small.
Figure~\ref{fig:2d:mach2to10} shows the flow solution
and mesh at all Mach numbers considered, including intermediate (partially
converged) ones. In all cases, the shock is tracked with curved, high-order
elements and the solution is well-resolved throughout the domain. The overall
mesh quality is also well-preserved as the shock boundary is substantially
compresses toward the cylinder.

\ifbool{fastcompile}{}{
\begin{figure}[!htbp]
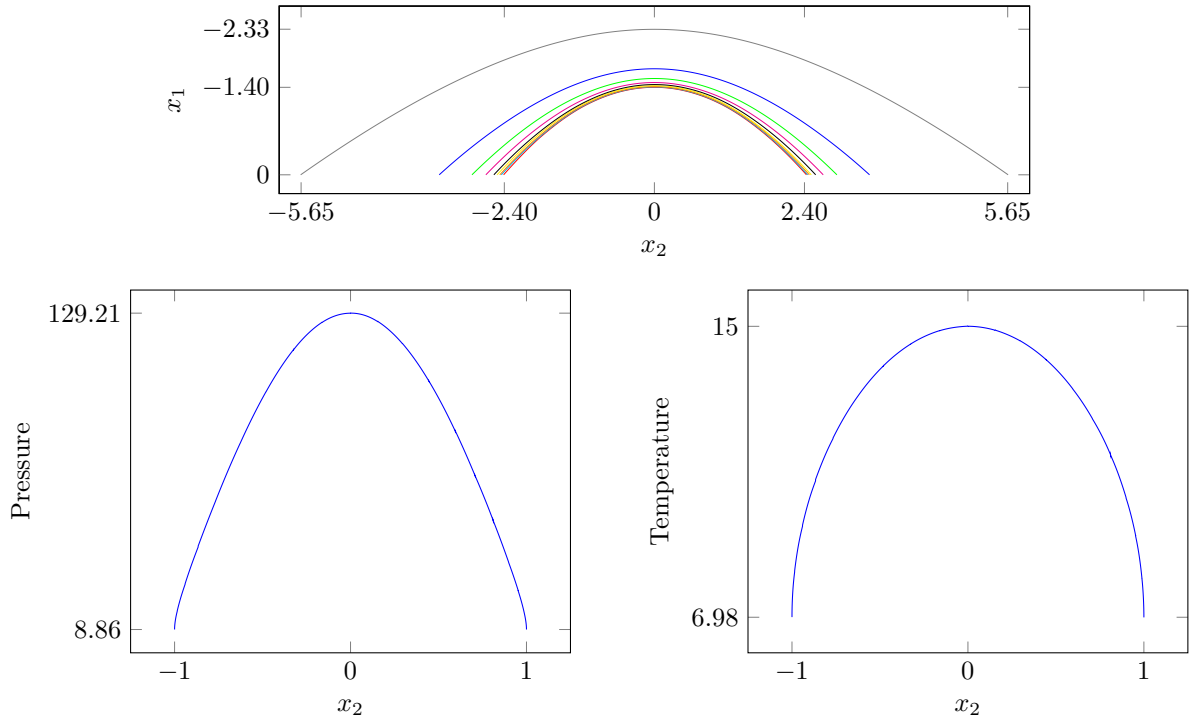

\centering
\input{./_py/mach2to10_shkpos.tikz} \\\vspace{3mm}
\input{./_py/mach2to10_pres.tikz} \qquad
\input{./_py/mach2to10_temp.tikz}
\caption{
 Lead shock positions (\textit{top}) for $M_\infty=2,3,4,\dots,10$,
 and pressure (\textit{bottom left}) and temperature (\textit{bottom right})
 along the surface of the cylinder at $M_\infty=10$ using Mach continuation
 starting from $M_\infty=2$ with $N=40$ stages.
 }
\label{fig:mach2to10_shkpos_temppres}
\end{figure}
}

\subsubsection{Inviscid flow over sphere, Mach continuation to $M_\infty = 3$}
\label{sec:numexp:3dmach2to3}
\begin{figure}[!htbp]
\centering
\ifbool{fastcompile}{}{
\raisebox{-0.5\height}{\input{./_py/sph1_geom.tikz}} \qquad
\raisebox{-0.4\height}{\includegraphics[width=0.21\textwidth]{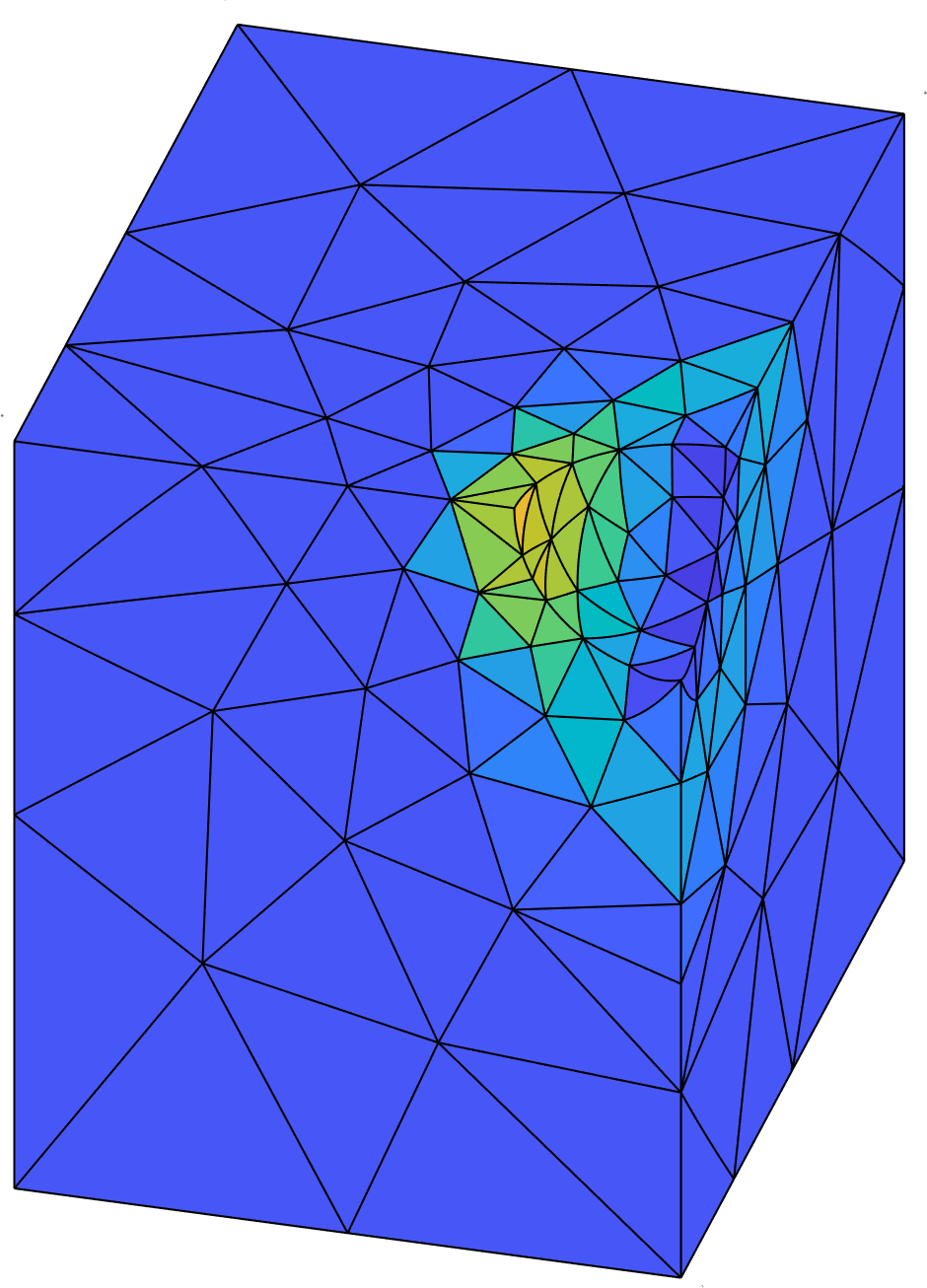}}\qquad
\raisebox{-0.4\height}{\includegraphics[width=0.21\textwidth]{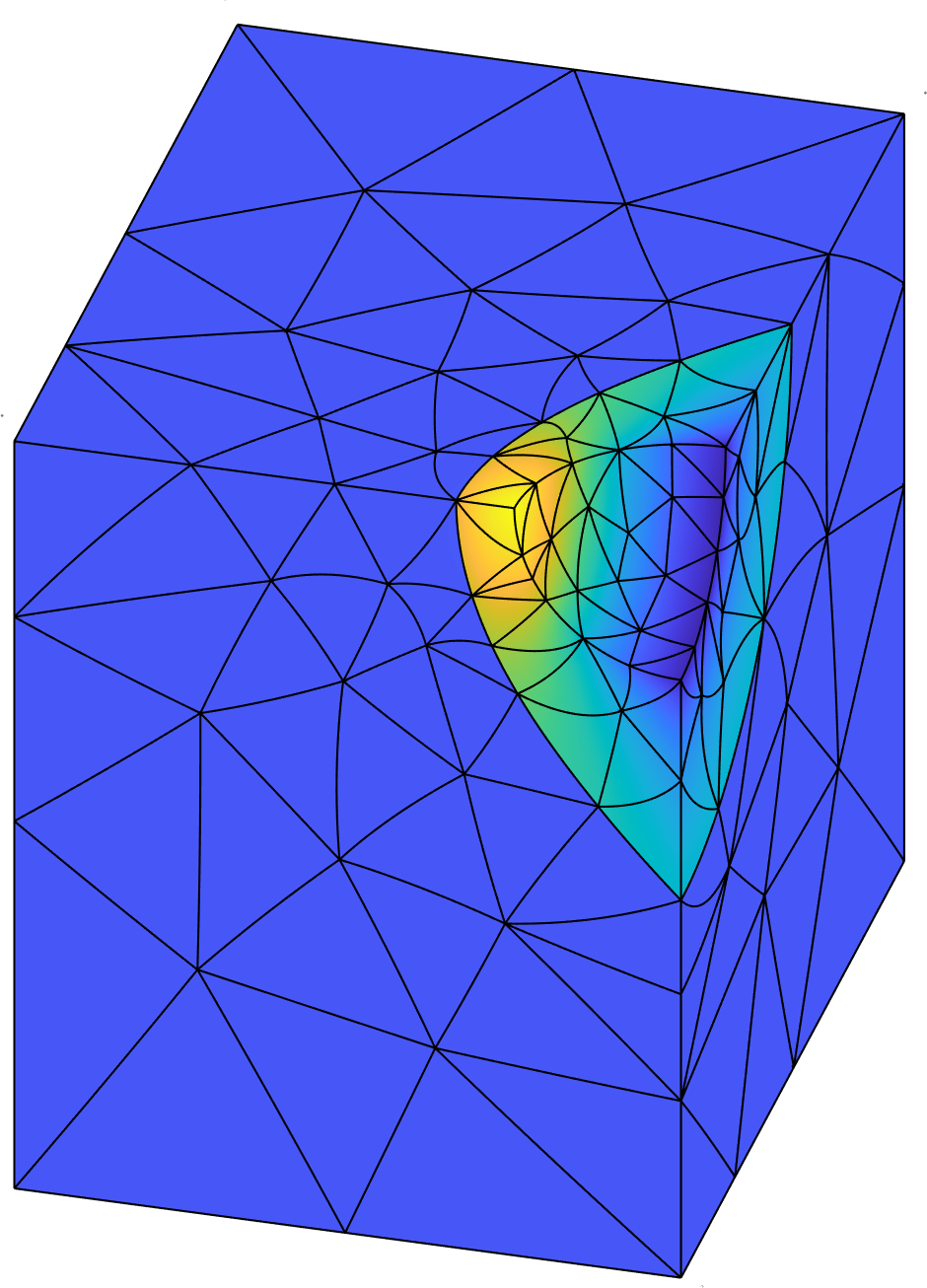}}\qquad
\raisebox{-0.4\height}{\includegraphics[width=0.088\textwidth]{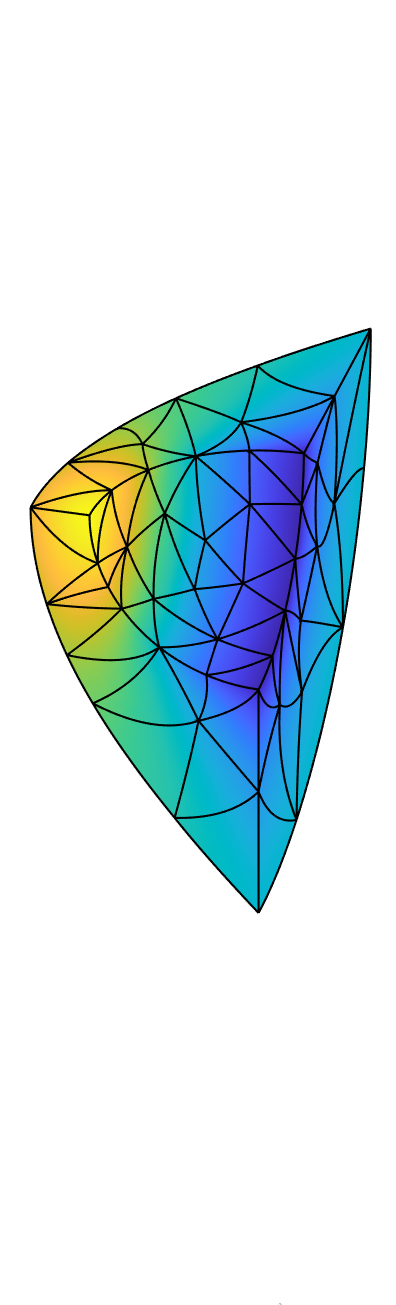}}
}
\colorbarMatlabParula{1}{2}{3}{4}{4.4}
 \caption{Flow domain (\textit{left}) and the $M_\infty=2$ initial guess
 (density) on a shock-agnostic mesh (\textit{middle-left}), the shock-aligned
 mesh and corresponding solution obtained using HOIST (\textit{middle-right}),
 and the corresponding solution and mesh extracted from downstream the bow
 shock (\textit{right}).}
\label{fig:3d:mach2}
\end{figure}

Next, we consider steady $M_\infty=3$ flow over a
sphere using Mach continuation beginning at $M_\infty=2$. At $M_\infty=3$,
we have the following stagnation quantities, which will be used to make
quantitative assessment of the accuracy of the HOIST method:
$p_{02}=12.061$, $T_0=2$, and $H_0=7$ (same as the 2D $M_\infty=3$
case in Section~\ref{sec:numexp:mach2to3}). To reduce the computational
cost, we model only an eighth of the geometry and use symmetry boundary
conditions. Our approach begins by applying the HOIST method at $M_\infty = 2$
on a shock-agnostic mesh of the entire domain consisting of $491$ elements
and extracting the solution and mesh downstream of the bow shock to initialize
our continuation strategy, resulting in a reduced mesh with $159$ elements
(Figure~\ref{fig:3d:mach2}).

\begin{figure}[!htbp]
\centering
\ifbool{fastcompile}{}{
\includegraphics[width=0.35\textwidth]{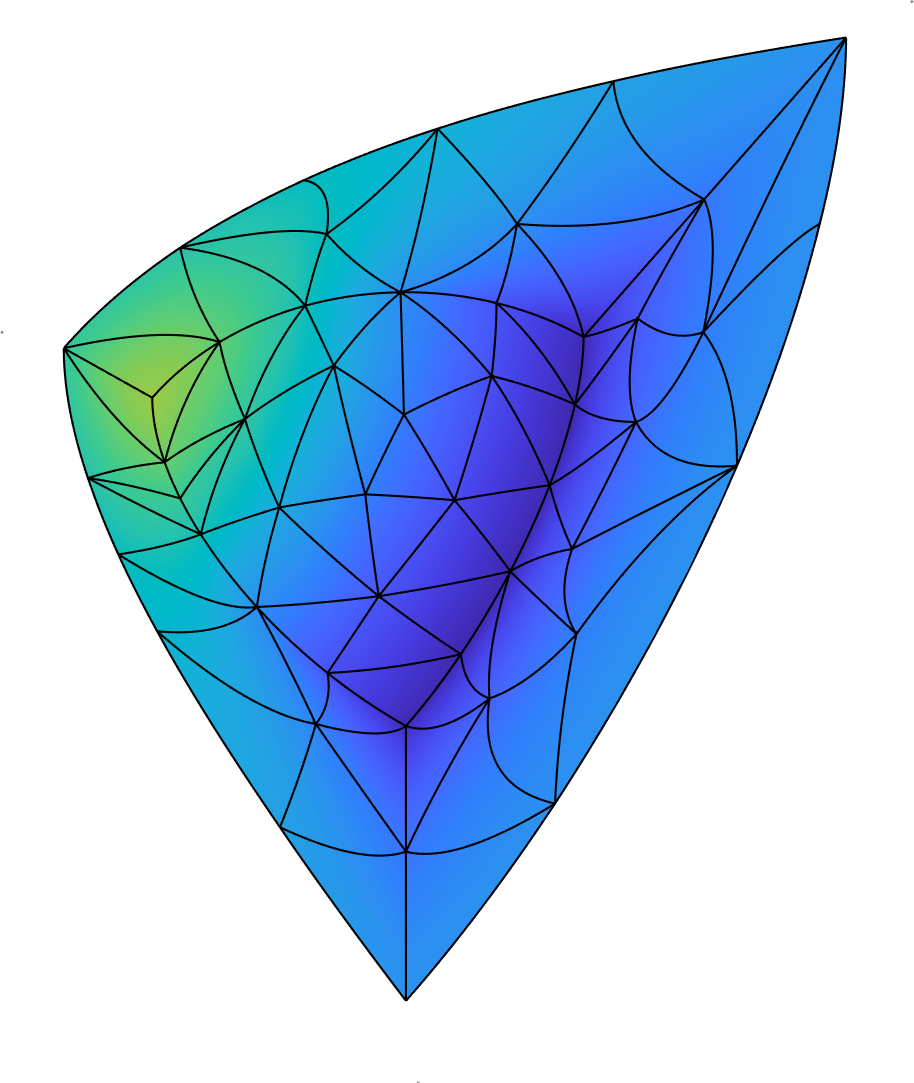} 
\includegraphics[width=0.35\textwidth]{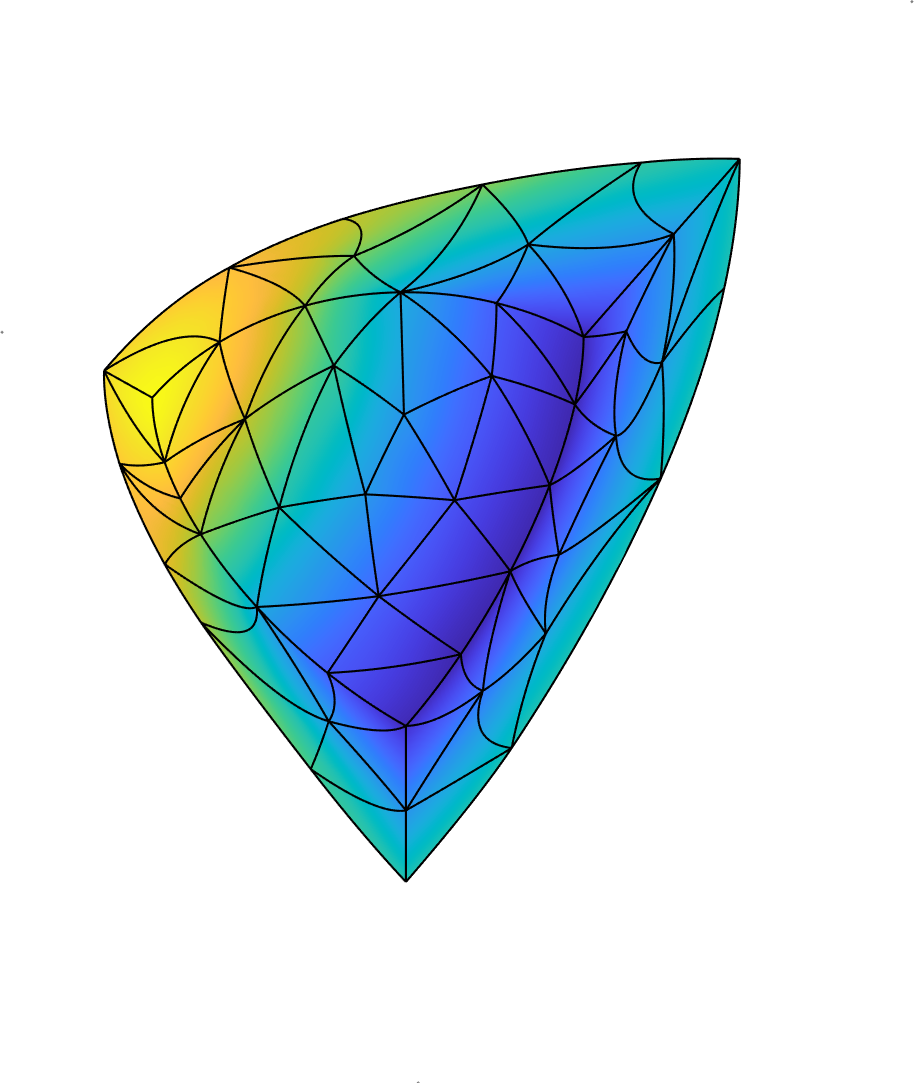} \\
\includegraphics[width=0.35\textwidth]{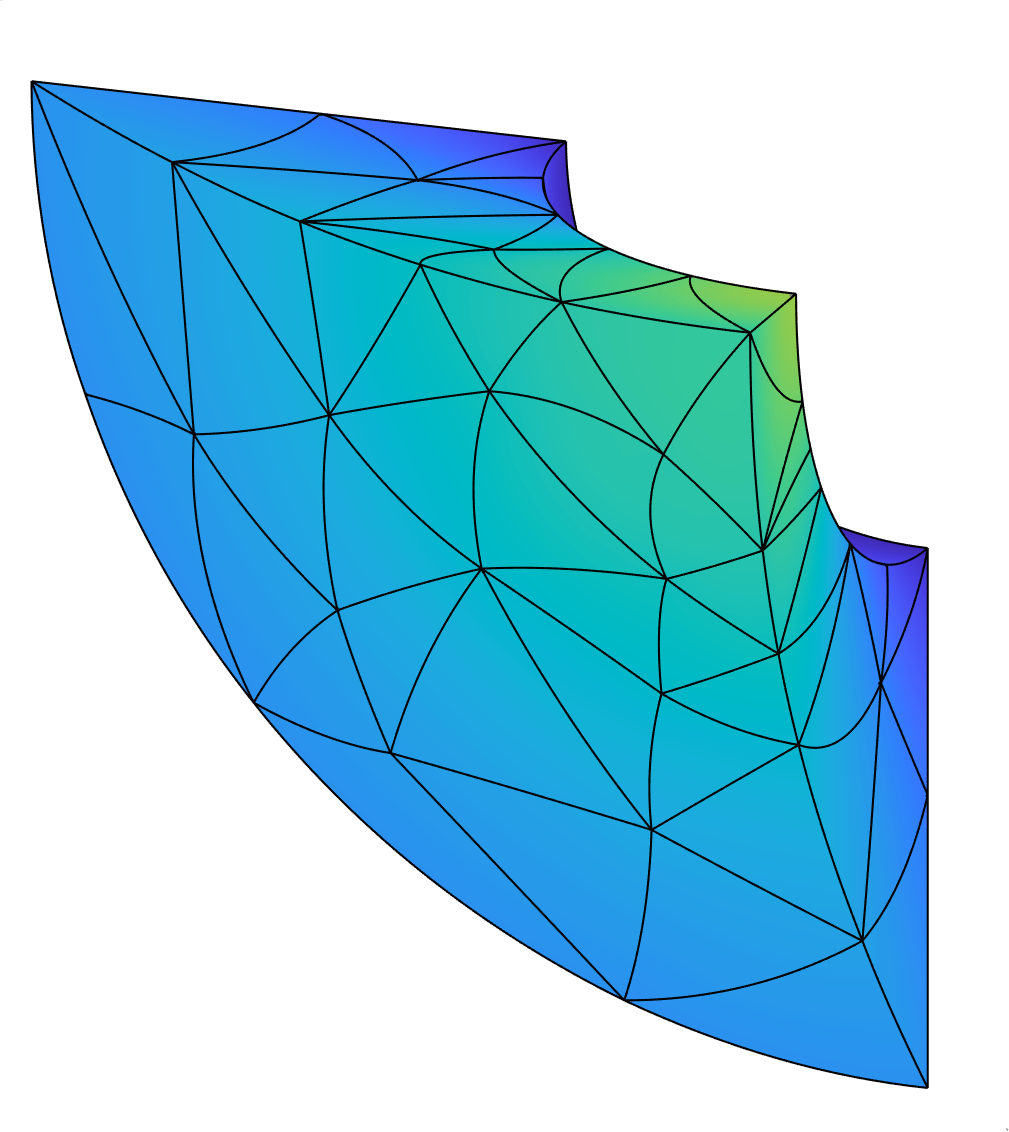} 
\includegraphics[width=0.35\textwidth]{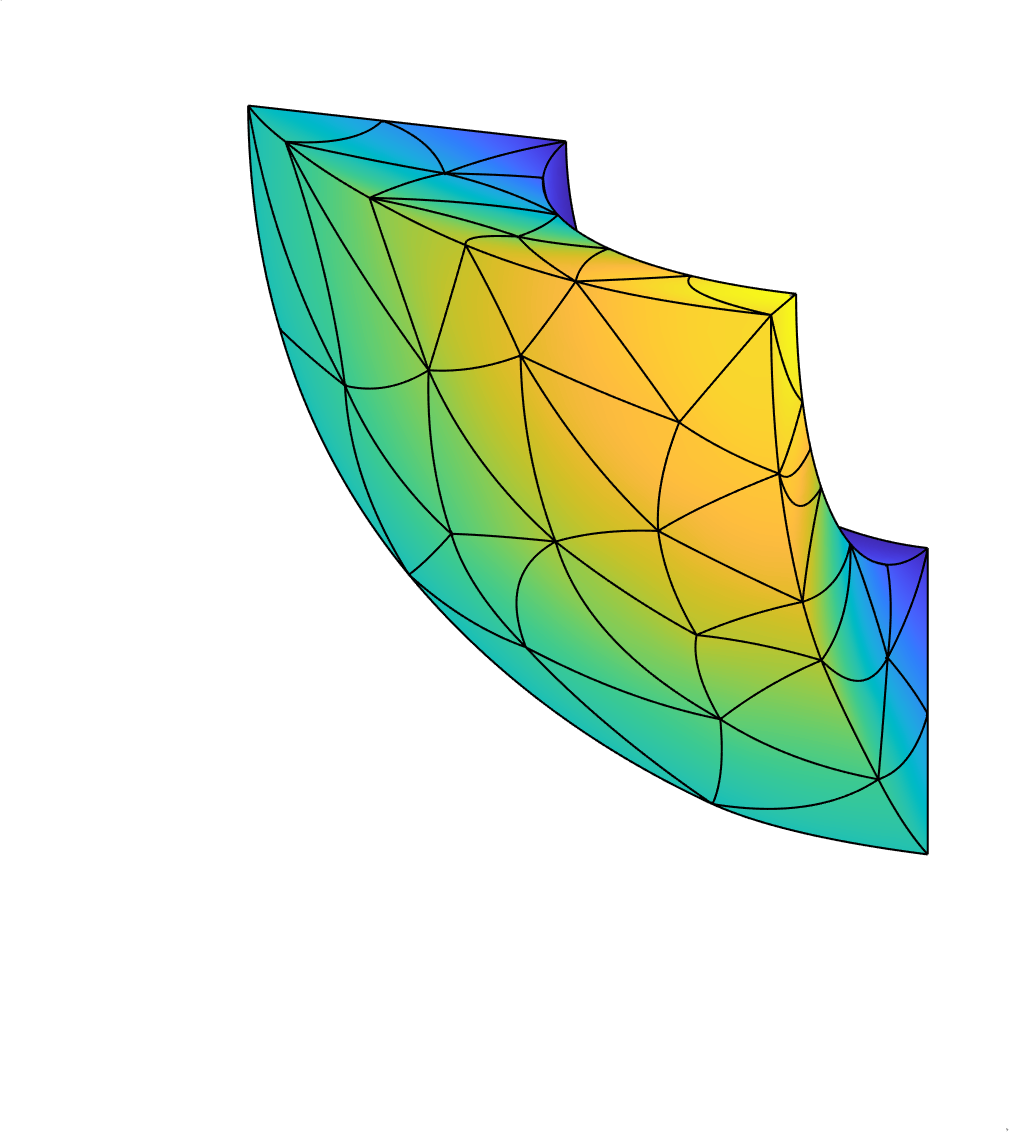} 
}
\colorbarMatlabParula{0.8}{2}{3}{4}{4.45}
 \caption{Density distribution at the initial $M_\infty=2$ (\textit{left})
 and final $M_\infty=3$ (\textit{right}) continuation stage with
 two views (\textit{top to bottom}).}
\label{fig:3d:mach2to3}
\end{figure}

We apply our continuation strategy with $N=10$ stages. The intermediate stages
are solved with at most $n_1=15$ iterations or to a tolerance of $\xi_1=10^{-4}$,
and the final stage is solved with at most $n_2=100$ iterations or to a
tolerance of $\xi_2=10^{-8}$. We use the same HOIST
solver parameters as in Section~\ref{sec:numexp:mach2to3} for all
intermediate stages, and we set $\lambda=0.5$, $\kappa_0=10^{-12}$ 
in the final stage to allow the solver to further improve the solution.
Similarly rapid convergence of the HOIST solver and fixed-mesh solve in
the final stage as the two-dimensional problems are observed
(Figure~\ref{fig:mach2to3_sph_hist}).
Figure~\ref{fig:3d:mach2to3} shows the flow solution and mesh at the
initial ($M_\infty=2$) and final ($M_\infty=3$) Mach numbers. We
observe that the high-order elements conform to the curvature of the
bow shock, which leads to a highly accurate solution on an extremely
coarse mesh (only $159$ quadratic elements) when combined with the
high-order flow field approximation. The overall mesh quality is
well-preserved as the shock boundary and mesh elements compress towards
the sphere. The value and relative error of stagnation quantities evaluated
at  the stagnation point $(x_1, x_2, x_3) = (-1, 0, 0)$ are reported in 
Table~\ref{tab:stag_err}, which demonstrates the high accuracy per DoF
of the overall approach for three-dimensional problems.
\ifbool{fastcompile}{}{
\begin{figure}[!htbp]
\centering
\begin{tikzpicture}
\begin{groupplot} [
group style={group size = 3 by 1, horizontal sep = 1.4cm, vertical sep = 0.8cm}]
\nextgroupplot[xtick={0,50,100}, xlabel={Iteration ($i$)}, ytick={1e-1, 1e-4}, xmin=-5, xmax=105, ymin=4e-05, ymax=0.2, ymode=log, width=0.315\textwidth]
\addplot [black, thin, mark=*, mark size=0.75, mark options={solid, thin}, mark repeat=2]
coordinates {
( 1.00000000e+00,  7.89279184e-02)
( 2.00000000e+00,  3.07979938e-02)
( 3.00000000e+00,  1.95198650e-02)
( 4.00000000e+00,  2.12819507e-02)
( 5.00000000e+00,  1.68101577e-02)
( 6.00000000e+00,  1.49575762e-02)
( 7.00000000e+00,  1.30828472e-02)
( 8.00000000e+00,  1.00581340e-02)
( 9.00000000e+00,  9.38440828e-03)
( 1.00000000e+01,  6.40085004e-03)
( 1.10000000e+01,  7.10679997e-03)
( 1.20000000e+01,  4.45060074e-03)
( 1.30000000e+01,  3.41273276e-03)
( 1.40000000e+01,  1.24035147e-03)
( 1.50000000e+01,  1.03189798e-03)
( 1.60000000e+01,  9.78501835e-04)
( 1.70000000e+01,  1.07097941e-03)
( 1.80000000e+01,  9.65462740e-04)
( 1.90000000e+01,  4.37896584e-04)
( 2.00000000e+01,  4.27091644e-04)
( 2.10000000e+01,  4.57784269e-04)
( 2.20000000e+01,  3.14177377e-04)
( 2.30000000e+01,  3.22384685e-04)
( 2.40000000e+01,  3.99045185e-04)
( 2.50000000e+01,  3.64727968e-04)
( 2.60000000e+01,  2.59285974e-04)
( 2.70000000e+01,  2.17385126e-04)
( 2.80000000e+01,  2.08109935e-04)
( 2.90000000e+01,  1.99388553e-04)};\label{line:R0}

\addplot [black, thin, mark=*, mark size=0.75, mark options={solid, thin}, mark repeat=10]
coordinates {
( 3.00000000e+01,  1.99960328e-04)
( 3.10000000e+01,  2.06326058e-04)
( 3.20000000e+01,  1.16404191e-04)
( 3.30000000e+01,  1.20341093e-04)
( 3.40000000e+01,  1.15415283e-04)
( 3.50000000e+01,  1.00323694e-04)
( 3.60000000e+01,  1.11677268e-04)
( 3.70000000e+01,  1.08653412e-04)
( 3.80000000e+01,  1.06025723e-04)
( 3.90000000e+01,  1.06394647e-04)
( 4.00000000e+01,  1.06230962e-04)
( 4.10000000e+01,  1.06723352e-04)
( 4.20000000e+01,  1.06923016e-04)
( 4.30000000e+01,  1.11675416e-04)
( 4.40000000e+01,  1.09665386e-04)
( 4.50000000e+01,  1.05508599e-04)
( 4.60000000e+01,  9.60491742e-05)
( 4.70000000e+01,  8.57086666e-05)
( 4.80000000e+01,  8.59776475e-05)
( 4.90000000e+01,  8.68638534e-05)
( 5.00000000e+01,  8.40482604e-05)
( 5.10000000e+01,  8.15451572e-05)
( 5.20000000e+01,  7.80280807e-05)
( 5.30000000e+01,  7.88122399e-05)
( 5.40000000e+01,  7.95167244e-05)
( 5.50000000e+01,  7.97482524e-05)
( 5.60000000e+01,  8.00137857e-05)
( 5.70000000e+01,  8.08686587e-05)
( 5.80000000e+01,  8.08753405e-05)
( 5.90000000e+01,  7.90352568e-05)
( 6.00000000e+01,  7.82762735e-05)
( 6.10000000e+01,  7.90238455e-05)
( 6.20000000e+01,  7.98675432e-05)
( 6.30000000e+01,  8.00684422e-05)
( 6.40000000e+01,  7.70772698e-05)
( 6.50000000e+01,  7.77541280e-05)
( 6.60000000e+01,  7.69511938e-05)
( 6.70000000e+01,  7.55957702e-05)
( 6.80000000e+01,  7.85617250e-05)
( 6.90000000e+01,  7.67824265e-05)
( 7.00000000e+01,  7.71147736e-05)
( 7.10000000e+01,  7.72897824e-05)
( 7.20000000e+01,  7.64043040e-05)
( 7.30000000e+01,  7.56988279e-05)
( 7.40000000e+01,  7.33102165e-05)
( 7.50000000e+01,  7.34834158e-05)
( 7.60000000e+01,  7.31168482e-05)
( 7.70000000e+01,  7.34325658e-05)
( 7.80000000e+01,  7.36892905e-05)
( 7.90000000e+01,  7.37333374e-05)
( 8.00000000e+01,  7.39649247e-05)
( 8.10000000e+01,  7.35110215e-05)
( 8.20000000e+01,  7.36343768e-05)
( 8.30000000e+01,  7.36213636e-05)
( 8.40000000e+01,  7.29381016e-05)
( 8.50000000e+01,  7.29113325e-05)
( 8.60000000e+01,  7.28323049e-05)
( 8.70000000e+01,  7.38977253e-05)
( 8.80000000e+01,  7.31645253e-05)
( 8.90000000e+01,  7.33214699e-05)
( 9.00000000e+01,  7.31659015e-05)
( 9.10000000e+01,  7.34692820e-05)
( 9.20000000e+01,  7.32268980e-05)
( 9.30000000e+01,  7.32106294e-05)
( 9.40000000e+01,  7.33156210e-05)
( 9.50000000e+01,  7.33790130e-05)
( 9.60000000e+01,  7.33796870e-05)
( 9.70000000e+01,  7.39743630e-05)
( 9.80000000e+01,  7.28792897e-05)
( 9.90000000e+01,  7.27056538e-05)};\label{line:R0}

\nextgroupplot[xtick={0,50,100}, xlabel={Iteration ($i$)}, ytick={0.07, 0.13}, xmin=-5, xmax=105, ymin=0.06, ymax=0.14, yticklabel style={/pgf/number format/fixed}, width=0.315\textwidth]
\addplot [blue, thin, mark=square*, mark size=0.75, mark options={solid, thin}, mark repeat=3]
coordinates {
( 1.00000000e+00,  1.28215955e-01)
( 2.00000000e+00,  8.32751107e-02)
( 3.00000000e+00,  7.93032994e-02)
( 4.00000000e+00,  8.92427542e-02)
( 5.00000000e+00,  7.76571229e-02)
( 6.00000000e+00,  7.76982412e-02)
( 7.00000000e+00,  7.52736614e-02)
( 8.00000000e+00,  7.49803119e-02)
( 9.00000000e+00,  7.38362247e-02)
( 1.00000000e+01,  7.39200685e-02)
( 1.10000000e+01,  7.33205942e-02)
( 1.20000000e+01,  7.35201250e-02)
( 1.30000000e+01,  7.29423430e-02)
( 1.40000000e+01,  7.34250126e-02)};\label{line:R1err}

\addplot [blue, thin, mark=square*, mark size=0.75, mark options={solid, thin}, mark repeat=10]
coordinates {
( 1.50000000e+01,  7.29280092e-02)
( 1.60000000e+01,  7.28268533e-02)
( 1.70000000e+01,  7.25655540e-02)
( 1.80000000e+01,  7.25483399e-02)
( 1.90000000e+01,  7.24536966e-02)
( 2.00000000e+01,  7.23925851e-02)
( 2.10000000e+01,  7.23440129e-02)
( 2.20000000e+01,  7.22977923e-02)
( 2.30000000e+01,  7.22338681e-02)
( 2.40000000e+01,  7.21814258e-02)
( 2.50000000e+01,  7.21474237e-02)
( 2.60000000e+01,  7.20850372e-02)
( 2.70000000e+01,  7.20740819e-02)
( 2.80000000e+01,  7.20605556e-02)
( 2.90000000e+01,  7.20234424e-02)
( 3.00000000e+01,  7.19858040e-02)
( 3.10000000e+01,  7.19688219e-02)
( 3.20000000e+01,  7.19182319e-02)
( 3.30000000e+01,  7.19031456e-02)
( 3.40000000e+01,  7.18705361e-02)
( 3.50000000e+01,  7.18534906e-02)
( 3.60000000e+01,  7.18255583e-02)
( 3.70000000e+01,  7.18193577e-02)
( 3.80000000e+01,  7.18141493e-02)
( 3.90000000e+01,  7.18123291e-02)
( 4.00000000e+01,  7.18111779e-02)
( 4.10000000e+01,  7.18101905e-02)
( 4.20000000e+01,  7.18096197e-02)
( 4.30000000e+01,  7.18063135e-02)
( 4.40000000e+01,  7.18037716e-02)
( 4.50000000e+01,  7.17789187e-02)
( 4.60000000e+01,  7.17689871e-02)
( 4.70000000e+01,  7.17539190e-02)
( 4.80000000e+01,  7.17528006e-02)
( 4.90000000e+01,  7.17450326e-02)
( 5.00000000e+01,  7.17413387e-02)
( 5.10000000e+01,  7.17283028e-02)
( 5.20000000e+01,  7.17194045e-02)
( 5.30000000e+01,  7.17174092e-02)
( 5.40000000e+01,  7.17162083e-02)
( 5.50000000e+01,  7.17152801e-02)
( 5.60000000e+01,  7.17147052e-02)
( 5.70000000e+01,  7.17079575e-02)
( 5.80000000e+01,  7.17035943e-02)
( 5.90000000e+01,  7.16961848e-02)
( 6.00000000e+01,  7.16933436e-02)
( 6.10000000e+01,  7.16893633e-02)
( 6.20000000e+01,  7.16854304e-02)
( 6.30000000e+01,  7.16848156e-02)
( 6.40000000e+01,  7.16778602e-02)
( 6.50000000e+01,  7.16753351e-02)
( 6.60000000e+01,  7.16717214e-02)
( 6.70000000e+01,  7.16698350e-02)
( 6.80000000e+01,  7.16626608e-02)
( 6.90000000e+01,  7.16606122e-02)
( 7.00000000e+01,  7.16586914e-02)
( 7.10000000e+01,  7.16581223e-02)
( 7.20000000e+01,  7.16458973e-02)
( 7.30000000e+01,  7.16432598e-02)
( 7.40000000e+01,  7.16390475e-02)
( 7.50000000e+01,  7.16387600e-02)
( 7.60000000e+01,  7.16353292e-02)
( 7.70000000e+01,  7.16334912e-02)
( 7.80000000e+01,  7.16324897e-02)
( 7.90000000e+01,  7.16322358e-02)
( 8.00000000e+01,  7.16307471e-02)
( 8.10000000e+01,  7.16301261e-02)
( 8.20000000e+01,  7.16300082e-02)
( 8.30000000e+01,  7.16270879e-02)
( 8.40000000e+01,  7.16258138e-02)
( 8.50000000e+01,  7.16249990e-02)
( 8.60000000e+01,  7.16245564e-02)
( 8.70000000e+01,  7.16212805e-02)
( 8.80000000e+01,  7.16202337e-02)
( 8.90000000e+01,  7.16086476e-02)
( 9.00000000e+01,  7.16035319e-02)
( 9.10000000e+01,  7.16029708e-02)
( 9.20000000e+01,  7.16027217e-02)
( 9.30000000e+01,  7.16026624e-02)
( 9.40000000e+01,  7.16025952e-02)
( 9.50000000e+01,  7.16024816e-02)
( 9.60000000e+01,  7.16024153e-02)
( 9.70000000e+01,  7.15967372e-02)
( 9.80000000e+01,  7.15942027e-02)
( 9.90000000e+01,  7.15937286e-02)};\label{line:R1err}

\nextgroupplot[xtick={0,2,4,6}, xlabel={Iteration ($i$)}, xmin=-0.25, xmax=4.25, ymin=2e-15, ymax=0.0036352826901350204, ymode=log, width=0.315\textwidth]
\addplot [red, thin, mark=triangle*, mark size=1, mark options={solid, thin}, mark repeat=0]
coordinates {
( 0.00000000e+00,  7.27056538e-05)
( 1.00000000e+00,  2.67682621e-05)
( 2.00000000e+00,  1.73722609e-07)
( 3.00000000e+00,  3.81781891e-14)
( 4.00000000e+00,  1.83969646e-14)};\label{line:R0_fix}

\end{groupplot}\end{tikzpicture}
\caption{
 SQP convergence history of the DG residual
 $\norm{\rbm(\ubm_k^i,\phibold(\ybm_k^i),\mubold_k)}$ (\ref{line:R0})
 and enriched DG residual $\norm{\Rbm(\ubm_k^i,\phibold(\ybm_k^i),\mubold_k)}$
 (\ref{line:R1err}), and the DG residual throughout the fixed-mesh
 Newton iterations
 $\norm{\rbm(\ubm_i,\phibold(\ybm_k^I),\mubold_k)}$ for the
 final stage from Mach 2.9 to 3 (\ref{line:R0_fix}) for $M_\infty=2$
 to $M_\infty=3$ continuation (sphere).
 }
\label{fig:mach2to3_sph_hist}
\end{figure}
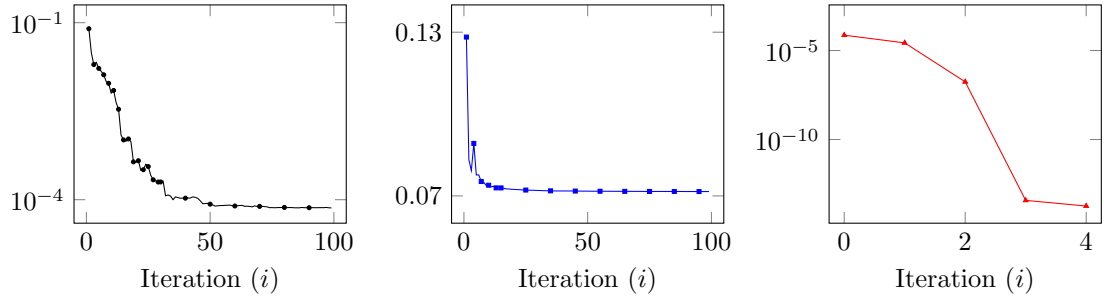
}

\begin{table}
\centering
\caption{Quantities of interest at the stagnation point produced by HOIST
 method at final continuation stage ($\check{p}_{02}$, $\check{T}_0$,
 $\check{H}_0$) and the corresponding relative error for all numerical
 experiments.}
\label{tab:mach2to3:stag_err}
\begin{tabular}{c|c|cc|cc|cc}
& $N$ & 
$\check{p}_{02}$ & $|p_{02}-\check{p}_{02}| / p_{02}$ & 
$\check{T}_0$ & $|T_0-\check{T}_0| / T_0$ & 
$\check{H}_0$ & $|H_0-\check{H}_0| / H_0$ \\\hline
$M=3$, 2D & 2 & 12.0610 & 6.8939e-06 & 1.9999 & 5.9976e-05 & 6.9996 & 5.8529e-05 \\
   & 5 & 12.0588 & 1.7742e-04 & 1.9998 & 9.4091e-05 & 6.9994 & 9.2435e-05 \\
   & 10 & 12.0592 & 1.4654e-04 & 1.9999 & 5.3712e-05 & 6.9996 & 5.1601e-05 \\
   & 20 & 12.0587 & 1.9162e-04 & 1.9998 & 8.3470e-05 & 6.9994 & 8.1685e-05 \\\hline

$M=10$, 2D & 40 & 129.2054 & 8.9226e-05 & 15.0029 & 1.9480e-04 & 52.5104 & 1.9905e-04 \\\hline

$M=3$, 3D & 10 & 11.9799 & 6.7254e-03 & 1.9990 & 5.1980e-04 & 6.9966 & 4.8529e-04

\end{tabular}
\label{tab:stag_err}
\end{table}

\subsubsection{Inviscid flow over double wedge, Mach continuation to $M_\infty = 6.8$}
\label{sec:numexp:dblwdg}
Next, we consider steady $M_\infty=6.8$ flow over a
double-wedge geometry with angles $15^\circ$ and $35^\circ$ 
using Mach continuation beginning at $M_\infty=2.8$. 
Our approach begins by applying the HOIST method at $M_\infty = 2.8$
on a shock-agnostic mesh of the entire domain consisting of $681$ elements
and extracting the solution and mesh downstream of the lead shock to initialize
our continuation strategy, resulting in a reduced mesh with $576$ elements
(Figure~\ref{fig:dw_ic}). In this experiment, the shock boundary consists 
of an oblique shock followed by a bow shock, which introduces highly non-uniform
mesh motion during the parameter sweep. In addition, this problem possesses 
complex shock-shock interactions and a supersonic jet downstream of the lead shock
as $M_\infty$ is increased.

We apply our continuation strategy with $N=121$ stages. The intermediate stages
are solved with $n_1=50$ iterations and the final stage is solved with $n_2=100$
iterations. The HOIST solver parameters
used for all $\Upsilon$ evaluations are $\lambda=10^{-1}$, $\kappa_0=10^{-2}$,
and $(\zeta, \upsilon)=(2, 0.5)$ \cite{huang2022robust}.
Figure~\ref{fig:dw:coarse_sweep} shows the flow solution and mesh at
selected Mach numbers during the parameter sweep. The density range
varies drastically from the initial Mach 2.8 flow to the final Mach 6.8 flow,
and complex shock-shock interactions emerge. As the inflow Mach increases, a
Type IV shock interaction \cite{candler1997dw} is observed and a thin layer of
elements conforms to the curved supersonic jet.

\ifbool{fastcompile}{}{
\begin{figure}[!htbp]
\centering
\input{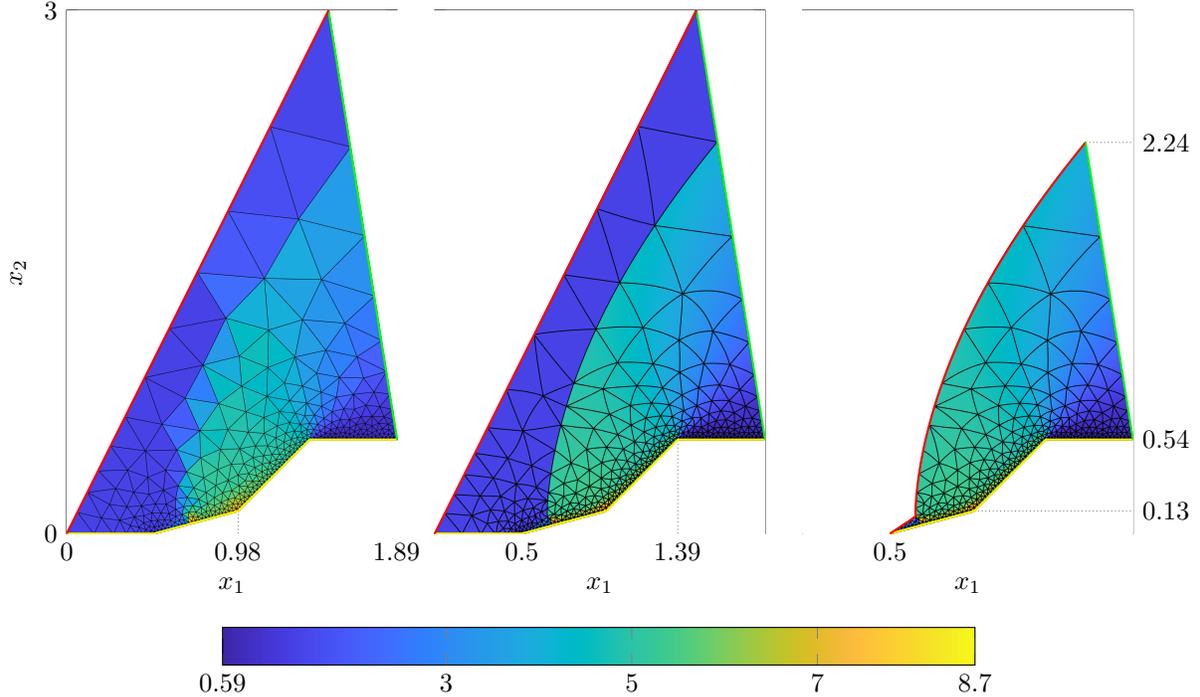}
\colorbarMatlabParula{0.59}{3}{5}{7}{8.70}
\caption{
 Initial guess (density) on a shock-agnostic mesh of $M_\infty = 2.8$
 flow over the double-wedge geometry (\textit{left}), the shock-aligned
 mesh and corresponding solution obtained using HOIST (\textit{middle}),
 and the corresponding solution and mesh extracted from downstream the lead
 shock (\textit{right}). The corresponding boundary conditions are inflow
 (\ref{line:dw:inlet}), slip walls (\ref{line:dw:wall}), and outflow
 (\ref{line:dw:outlet}).
}
\label{fig:dw_ic}
\end{figure}
}

\ifbool{fastcompile}{}{
\begin{landscape}
\begin{figure}[!htbp]
\centering
\begin{tikzpicture}
\begin{groupplot}[
  group style={
      group size=6 by 2,
      horizontal sep=0.1cm,
      vertical sep=0.1cm
  },
  width=0.525\textwidth,
  axis equal image,
  xlabel={},
  ylabel={},
  xtick = {1.72},
  ytick = {1.60},
  xmin=0.5, xmax=1.89254,
  ymin=0, ymax=2.242928,
  yticklabel pos=right,
  xticklabel style={/pgf/number format/.cd,fixed zerofill,precision=2},
  yticklabel style={/pgf/number format/.cd,fixed zerofill,precision=2}
]

\nextgroupplot[yticklabel pos=left, yticklabels={,,}, xticklabels={,,}]
\addplot graphics [xmin=0.5, xmax=1.89254, ymin=0, ymax=2.242928] {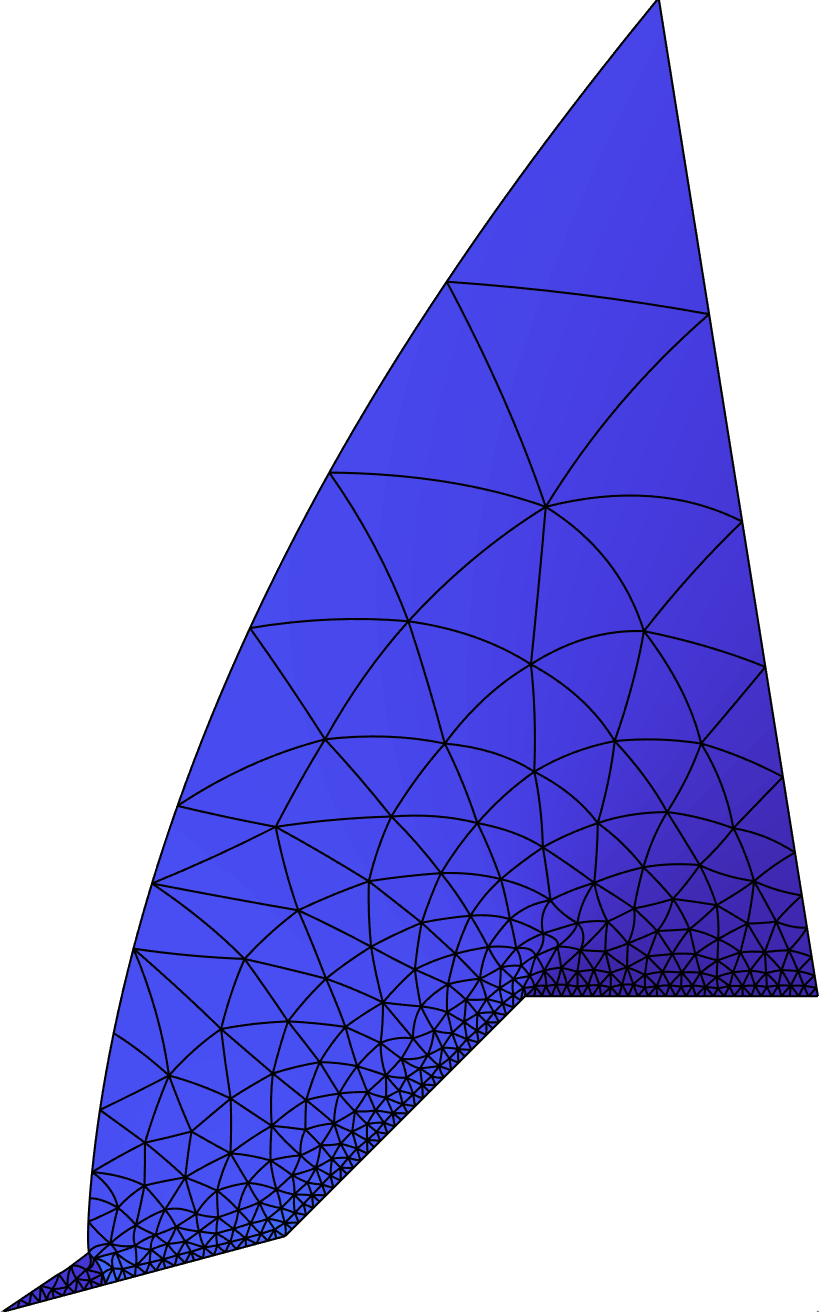};
\nextgroupplot[yticklabel pos=left, yticklabels={,,}, xticklabels={,,}]
\addplot graphics [xmin=0.5, xmax=1.89254, ymin=0, ymax=2.242928] {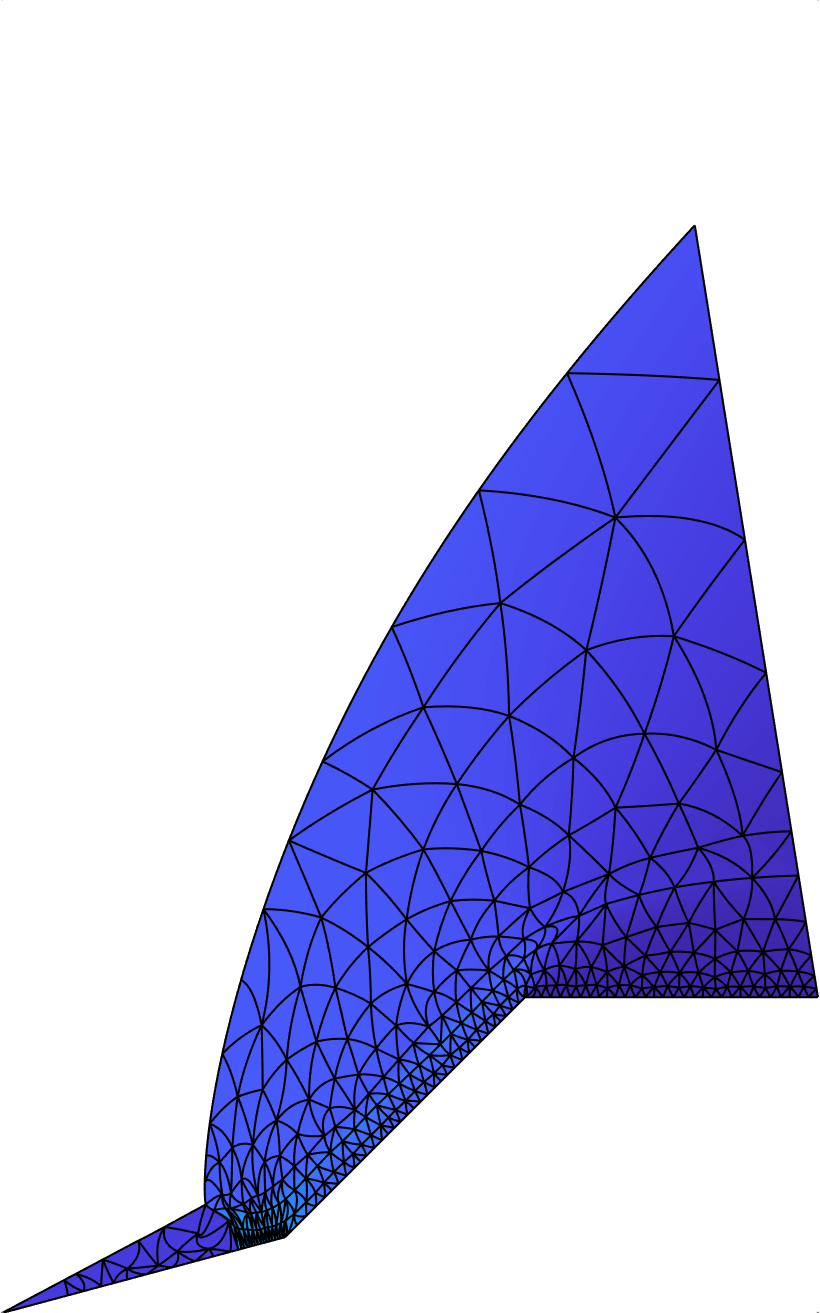};
\nextgroupplot[yticklabel pos=left, yticklabels={,,}, xticklabels={,,}]
\addplot graphics [xmin=0.5, xmax=1.89254, ymin=0, ymax=2.242928] {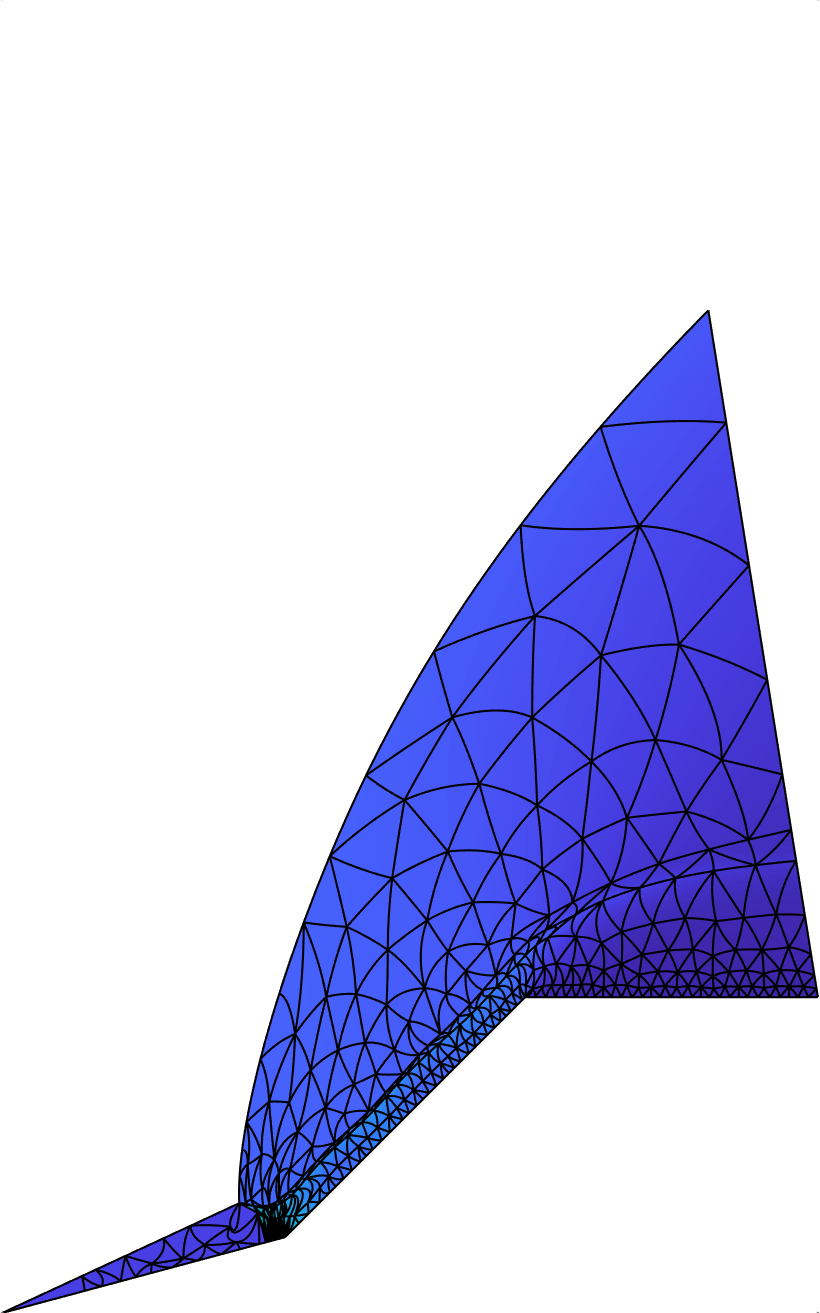};
\nextgroupplot[yticklabel pos=left, yticklabels={,,}, xticklabels={,,}]
\addplot graphics [xmin=0.5, xmax=1.89254, ymin=0, ymax=2.242928] {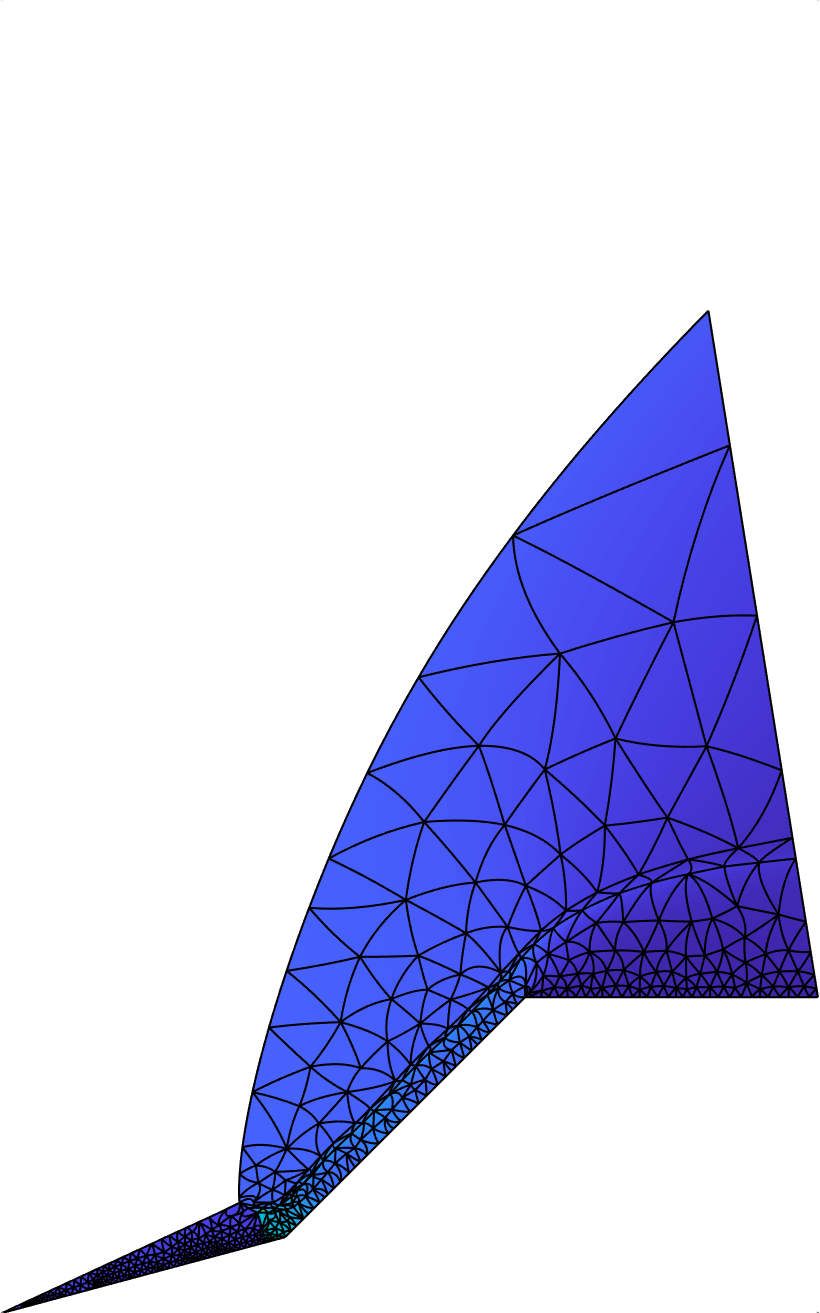};
\nextgroupplot[yticklabel pos=left, yticklabels={,,}, xticklabels={,,}]
\addplot graphics [xmin=0.5, xmax=1.89254, ymin=0, ymax=2.242928] {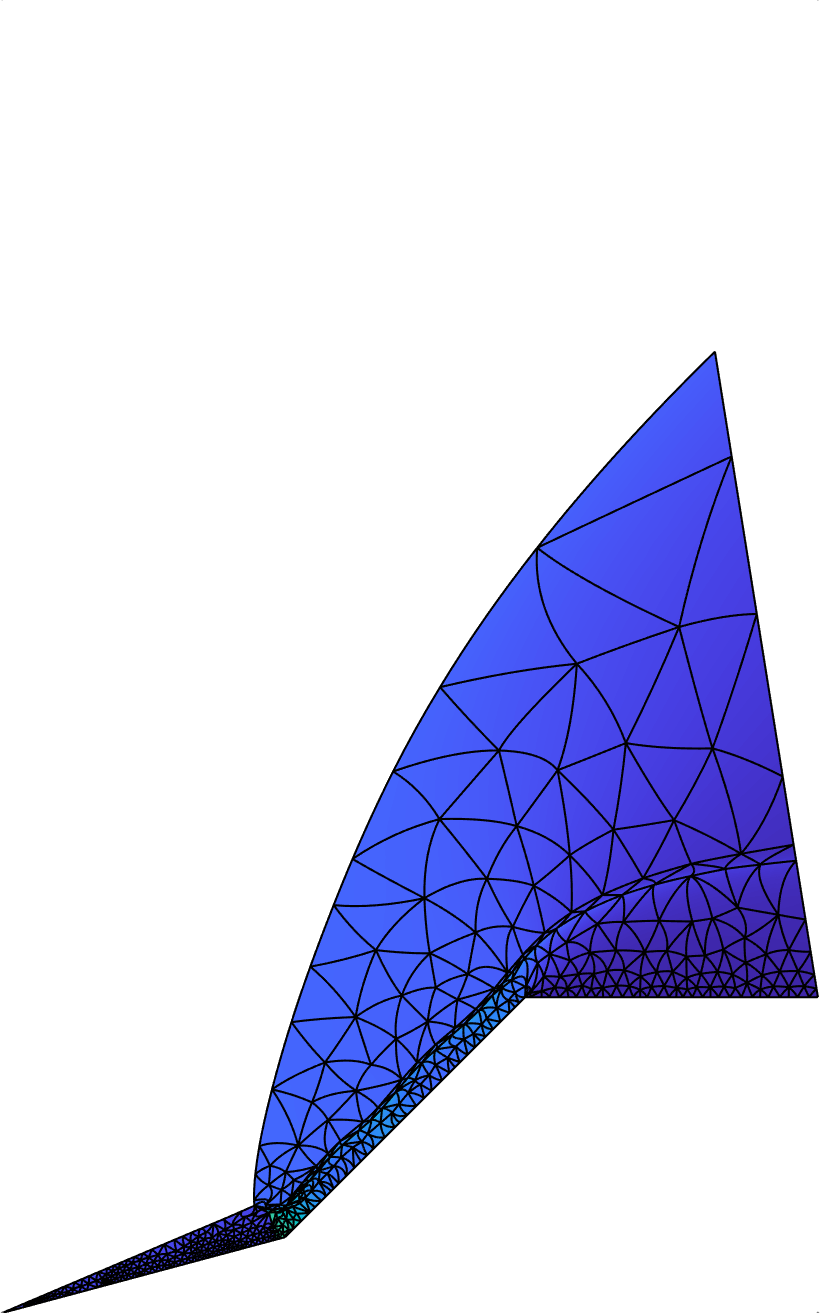};
\nextgroupplot[ylabel={}, xticklabel pos=top]
\addplot graphics [xmin=0.5, xmax=1.89254, ymin=0, ymax=2.242928] {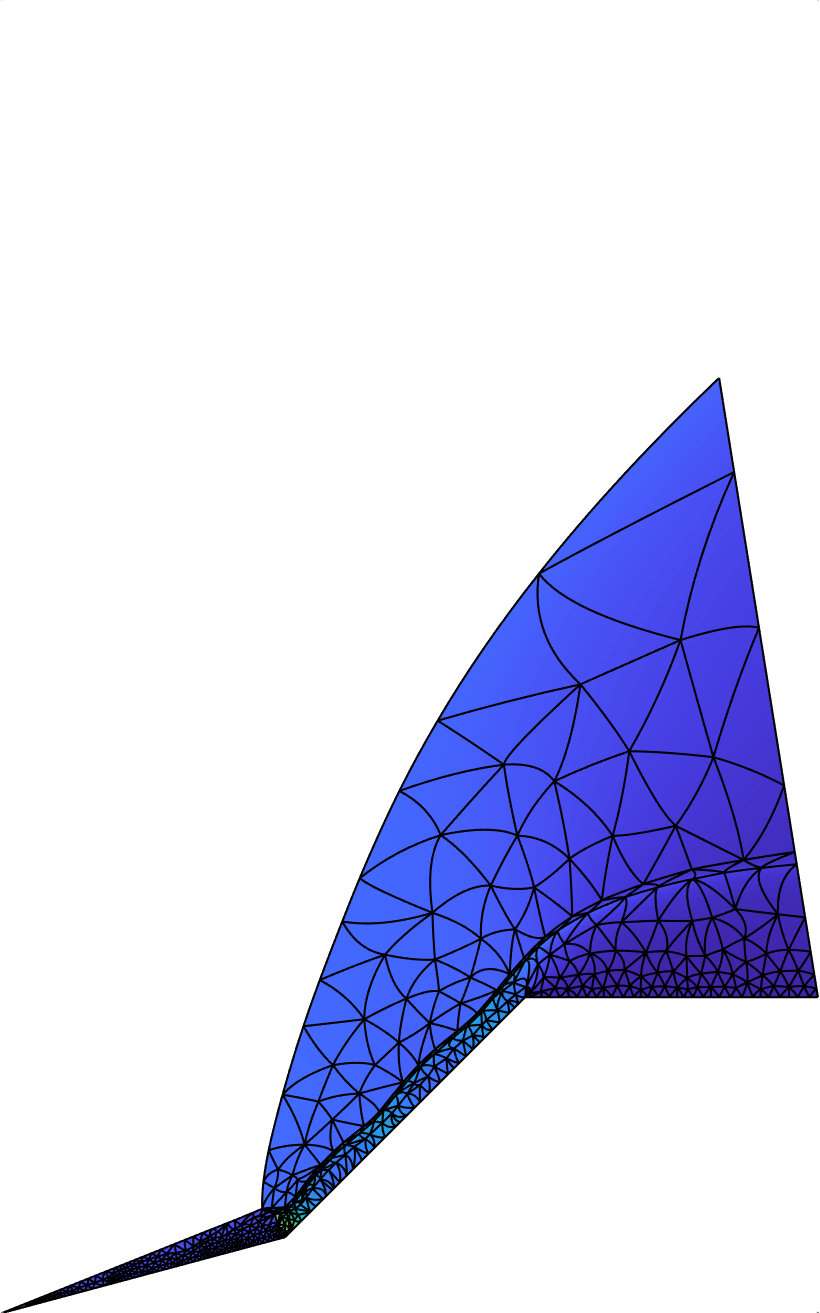};
\addplot [solid, gray, dashed]
coordinates {
( 1.724,  1.596)
( 1.89254,  1.596)};
\addplot [solid, gray, dashed]
coordinates {
( 1.724,  1.596)
( 1.724,  2.242928)};

\nextgroupplot[yticklabel pos=left, yticklabels={,,}, xticklabels={,,}]
\addplot graphics [xmin=0.5, xmax=1.89254, ymin=0, ymax=2.242928] {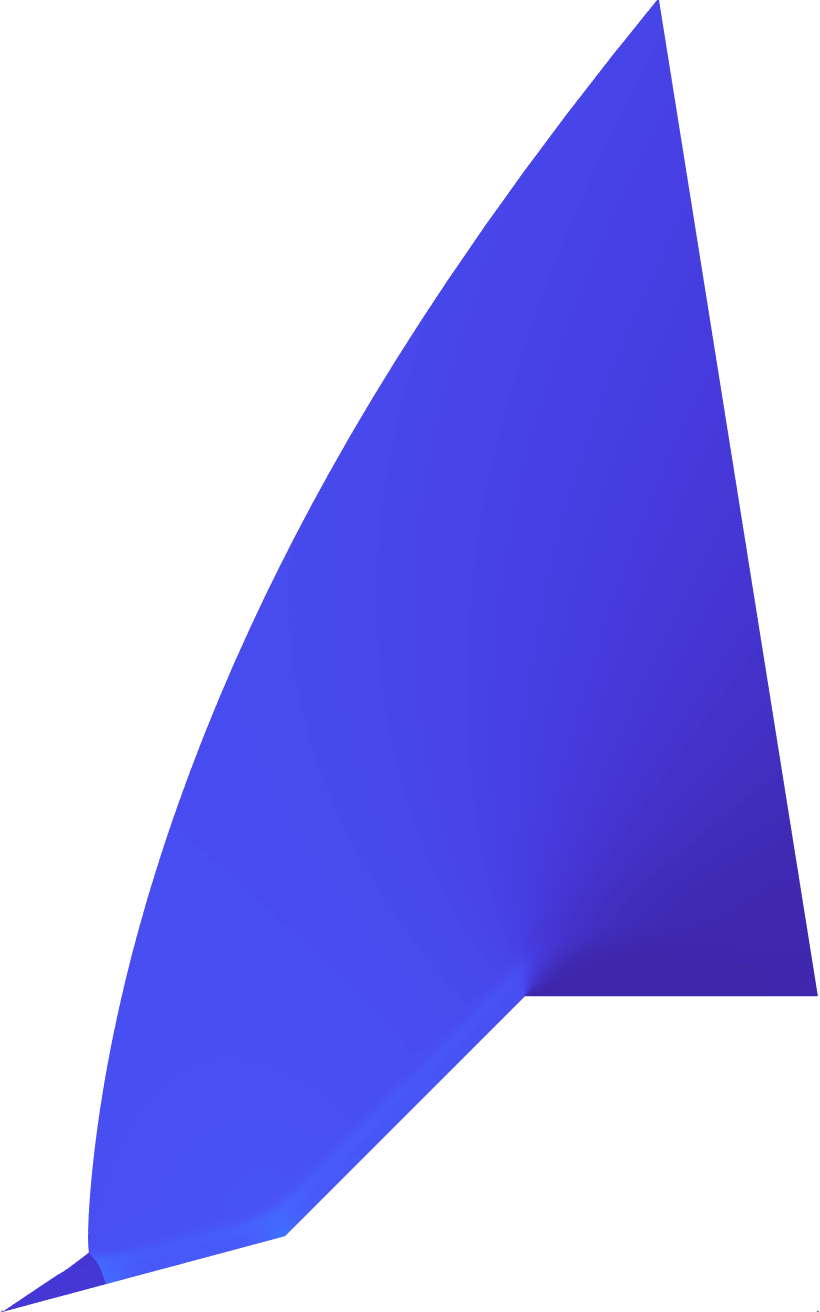};
\nextgroupplot[yticklabel pos=left, yticklabels={,,}, xticklabels={,,}]
\addplot graphics [xmin=0.5, xmax=1.89254, ymin=0, ymax=2.242928] {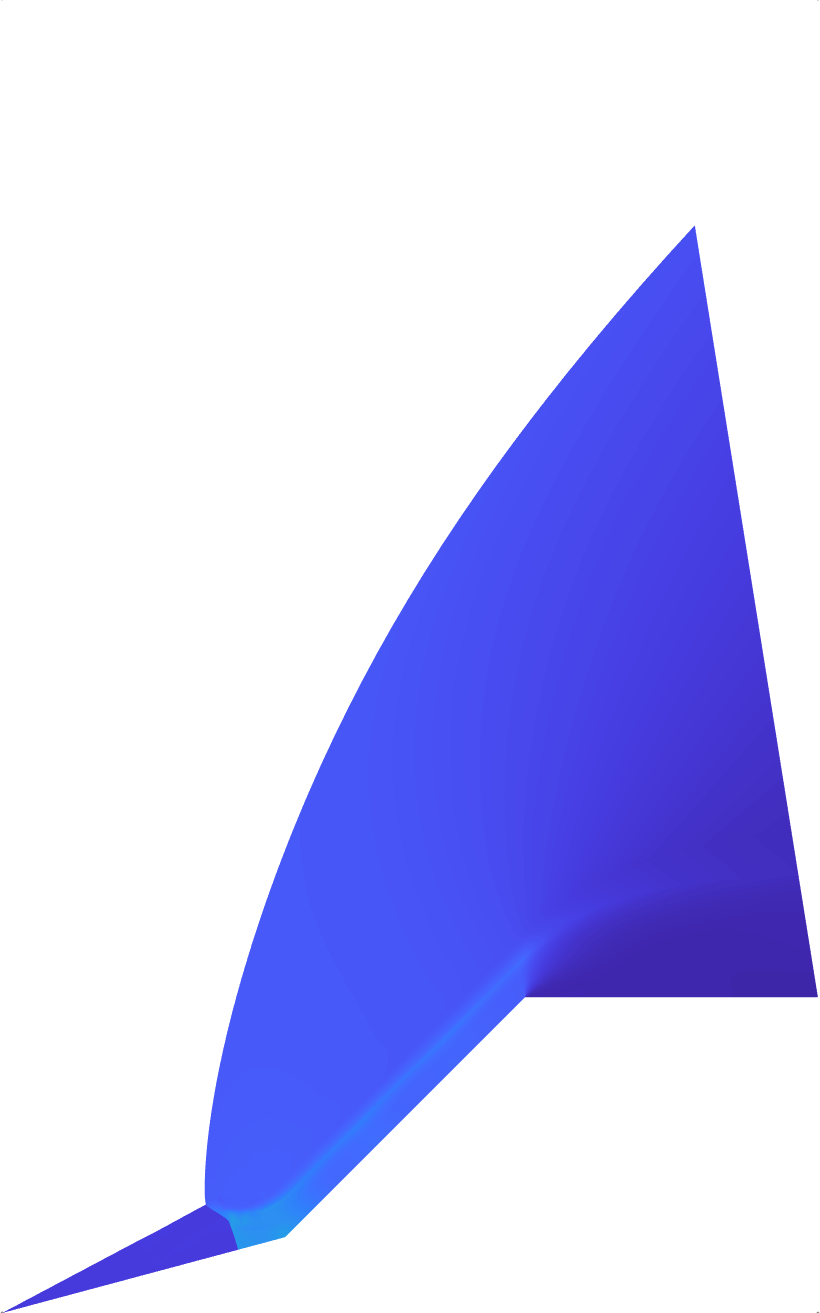};
\nextgroupplot[yticklabel pos=left, yticklabels={,,}, xticklabels={,,}]
\addplot graphics [xmin=0.5, xmax=1.89254, ymin=0, ymax=2.242928] {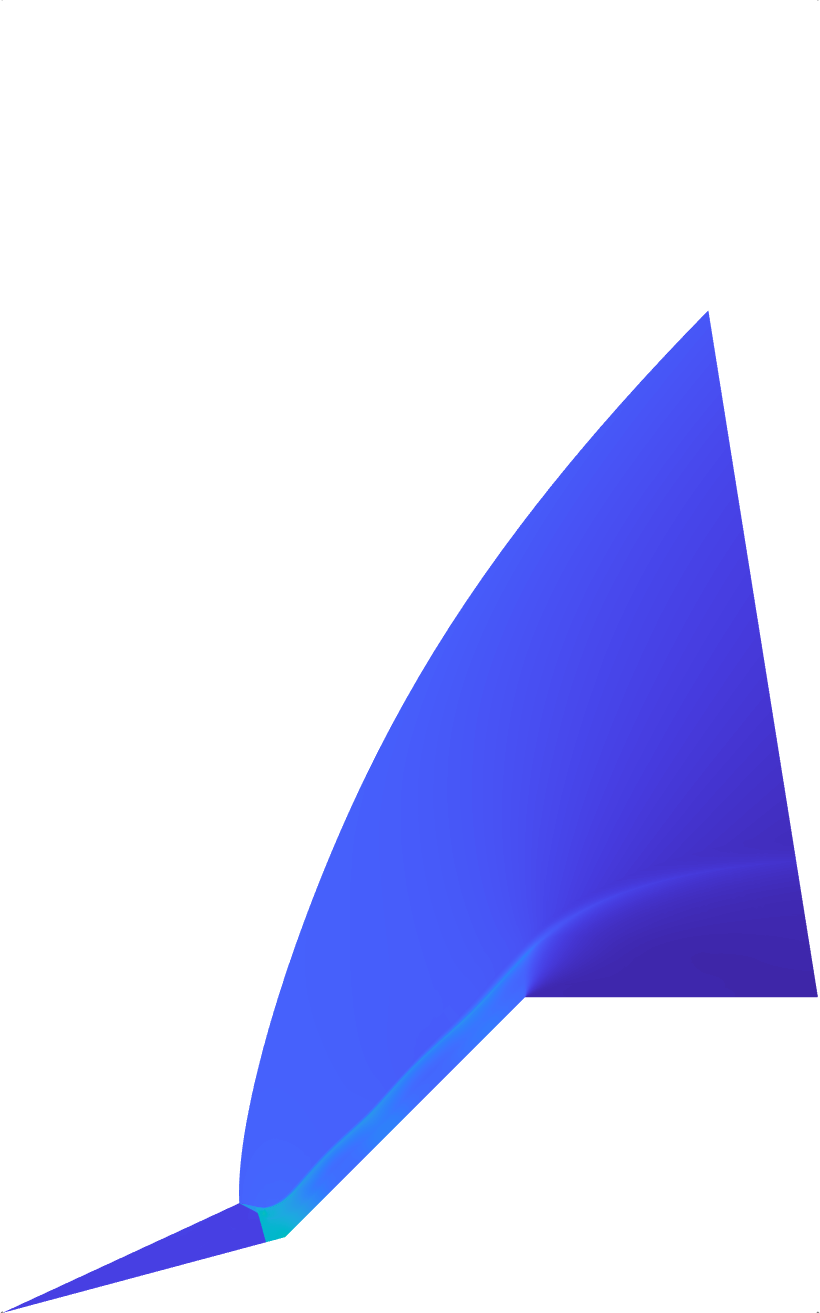};
\nextgroupplot[yticklabel pos=left, yticklabels={,,}, xticklabels={,,}]
\addplot graphics [xmin=0.5, xmax=1.89254, ymin=0, ymax=2.242928] {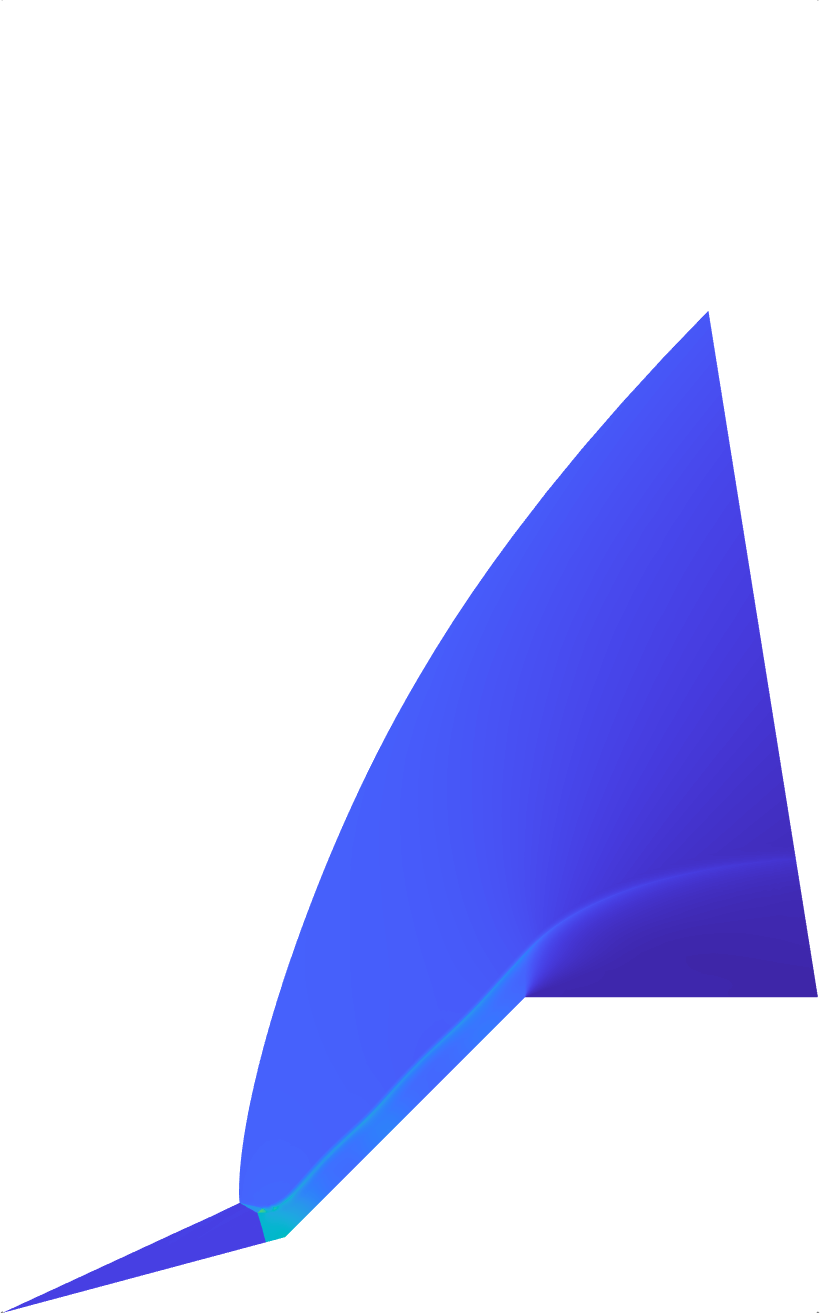};
\nextgroupplot[yticklabel pos=left, yticklabels={,,}, xticklabels={,,}]
\addplot graphics [xmin=0.5, xmax=1.89254, ymin=0, ymax=2.242928] {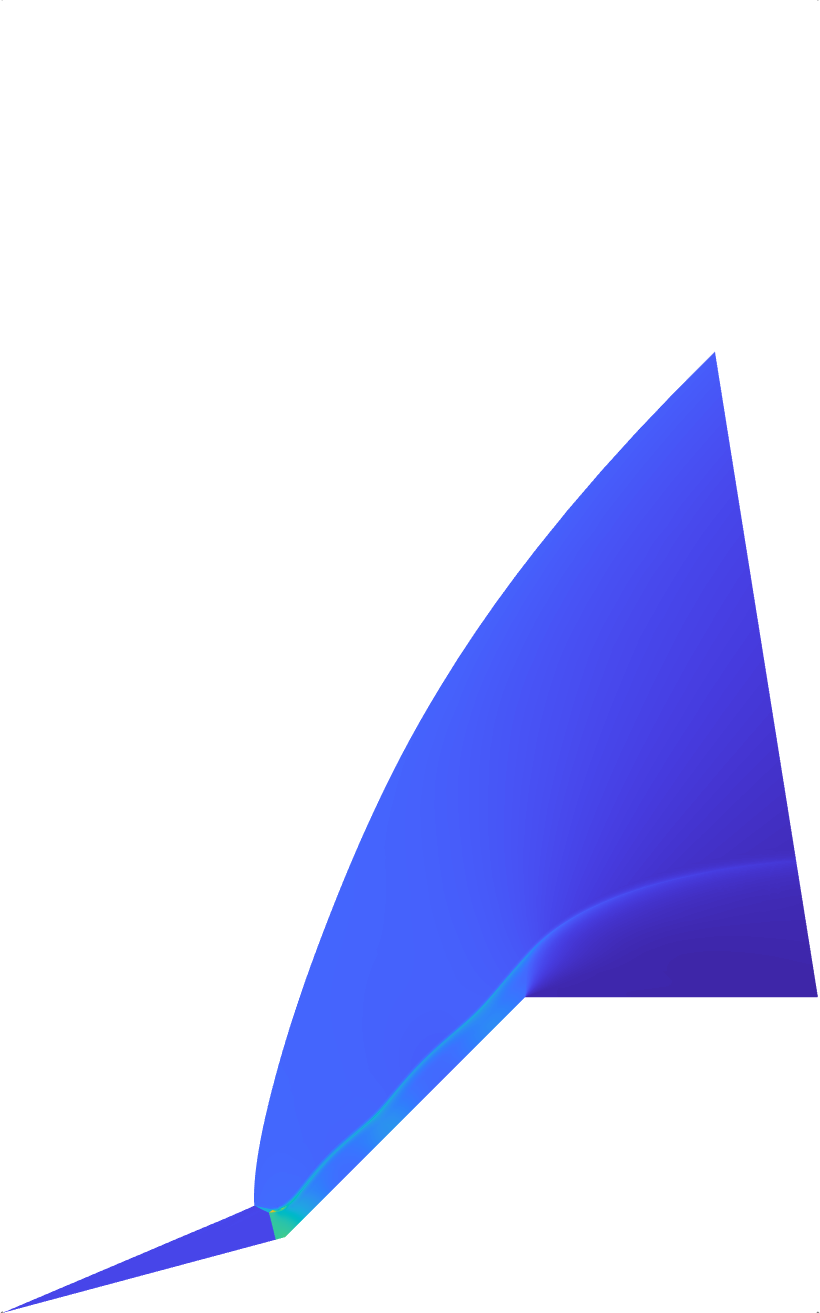};
\nextgroupplot[yticklabel pos=left, yticklabels={,,}, xticklabels={,,}]
\addplot graphics [xmin=0.5, xmax=1.89254, ymin=0, ymax=2.242928] {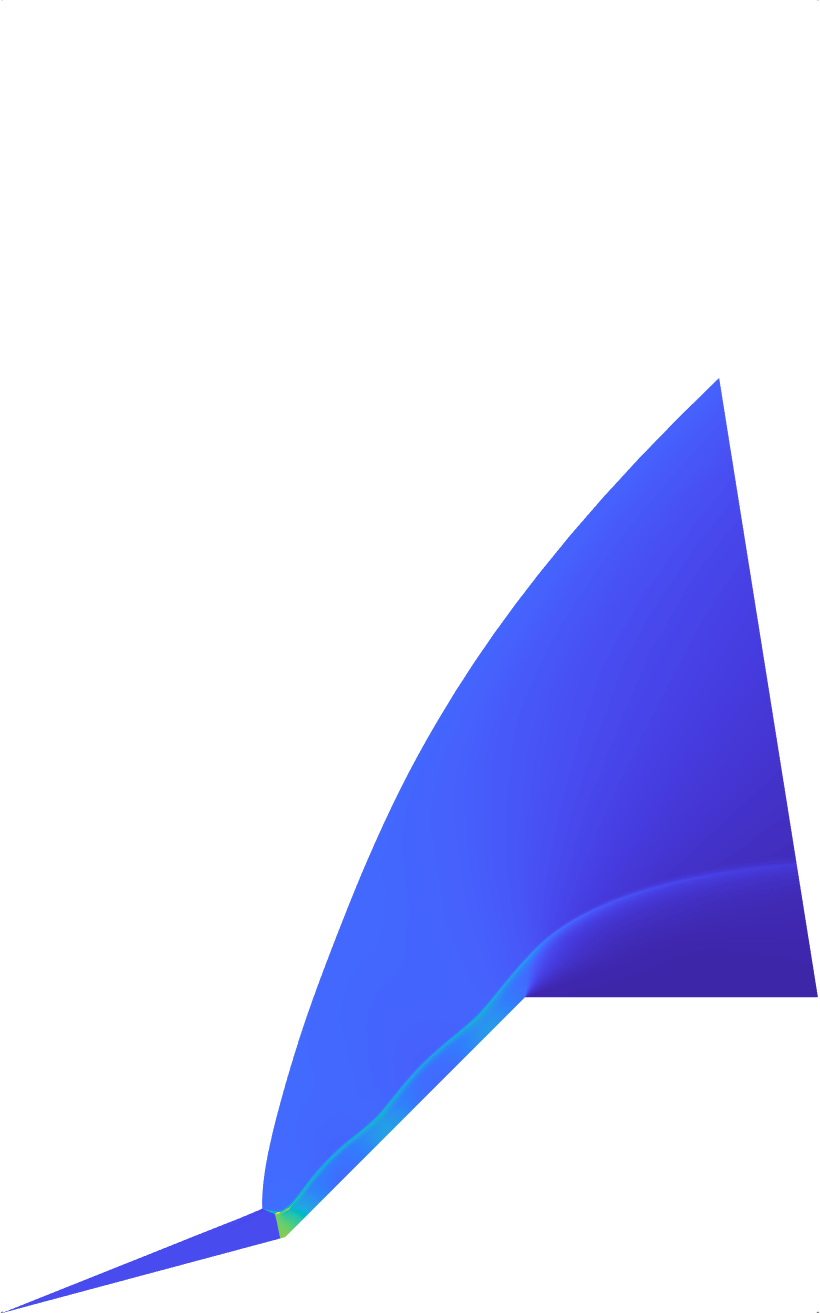};

\end{groupplot}
\end{tikzpicture}

\colorbarMatlabParula{0.66}{8}{18}{28}{33.61}
 \caption{
 Density distribution for double-wedge simulation (\textit{left to right}):
 $M_\infty = 2.8$,
 $M_\infty = 3.77$,
 $M_\infty = 4.77$ (stage 60, before re-mesh),
 $M_\infty = 4.77$ (stage 60, after re-mesh),
 $M_\infty = 5.77$, and
 $M_\infty = 6.8$.
  }
\label{fig:dw:coarse_sweep}
\end{figure}
\end{landscape}
}

\begin{figure}[!htbp]
\centering
\ifbool{fastcompile}{}{
\includegraphics[width=0.24\textwidth]{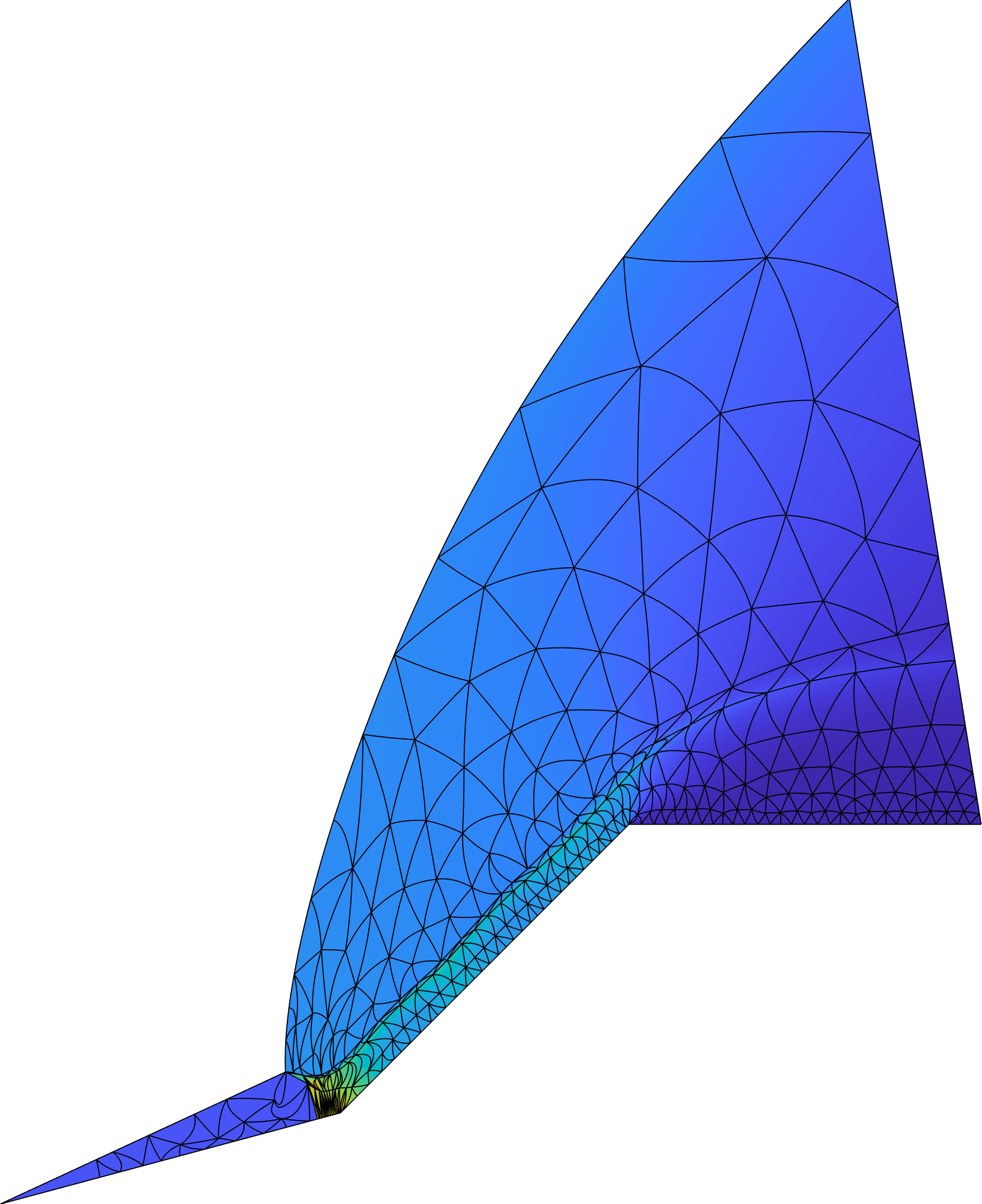} 
\includegraphics[width=0.24\textwidth]{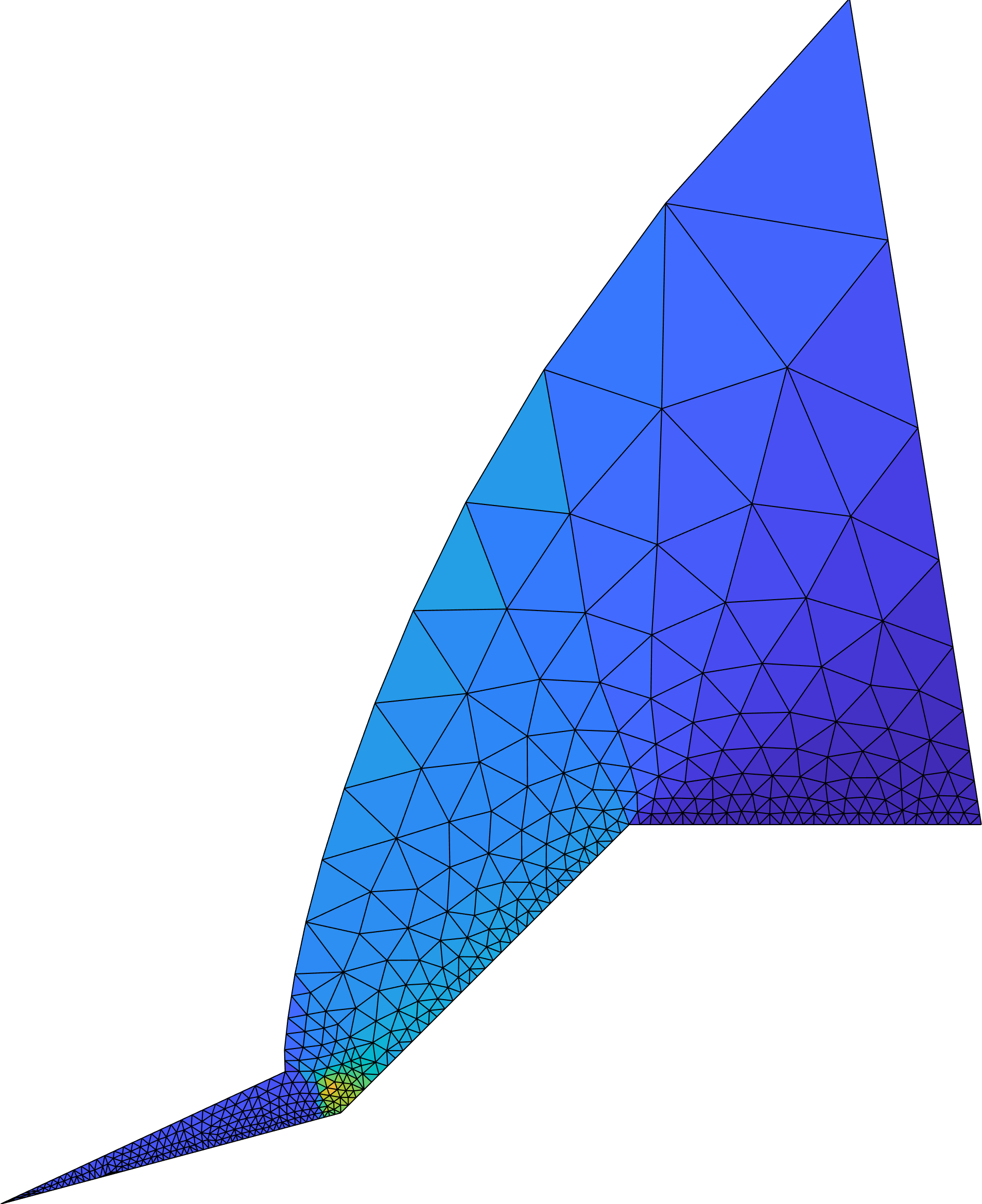} 
\includegraphics[width=0.24\textwidth]{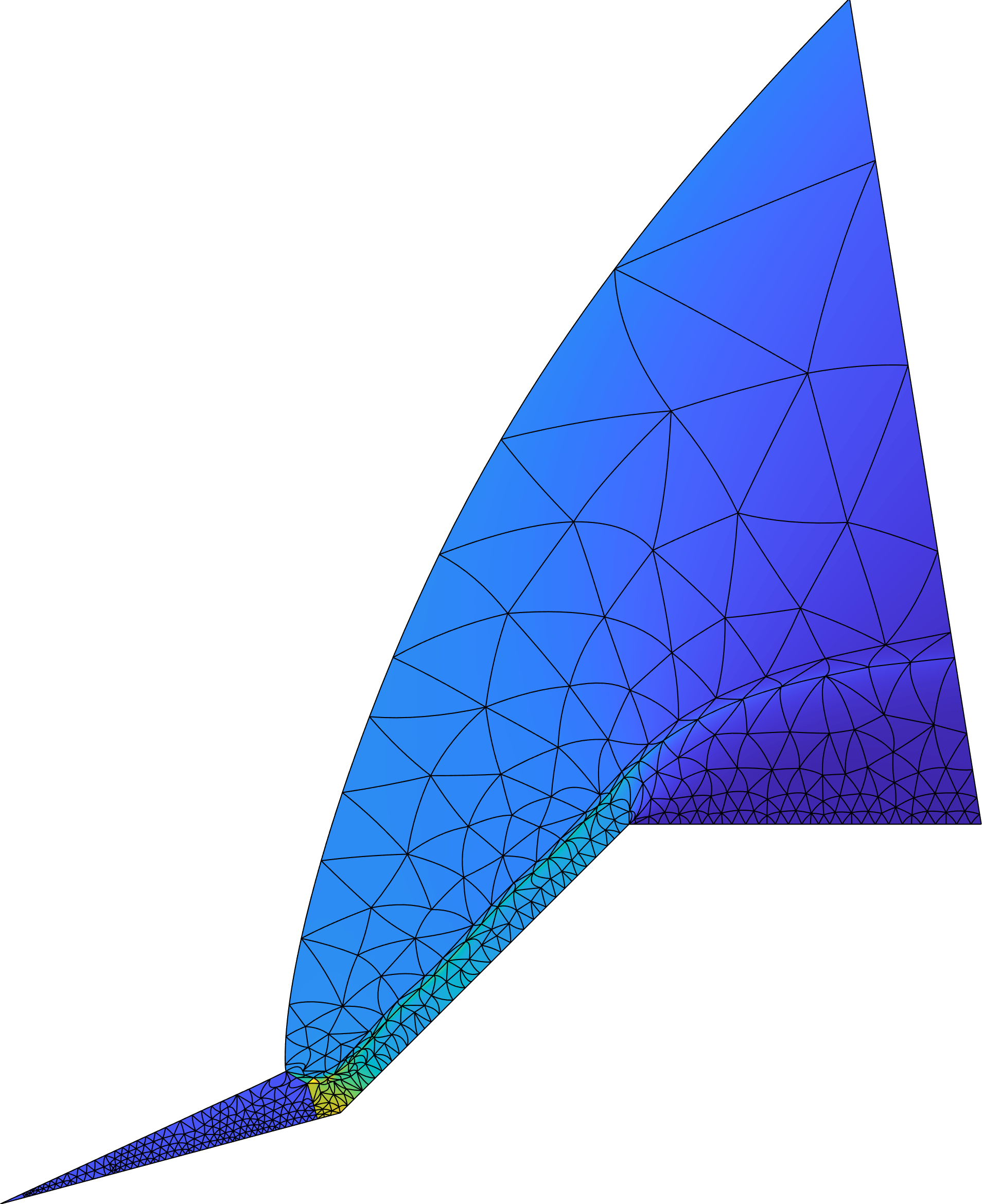}
\includegraphics[width=0.24\textwidth]{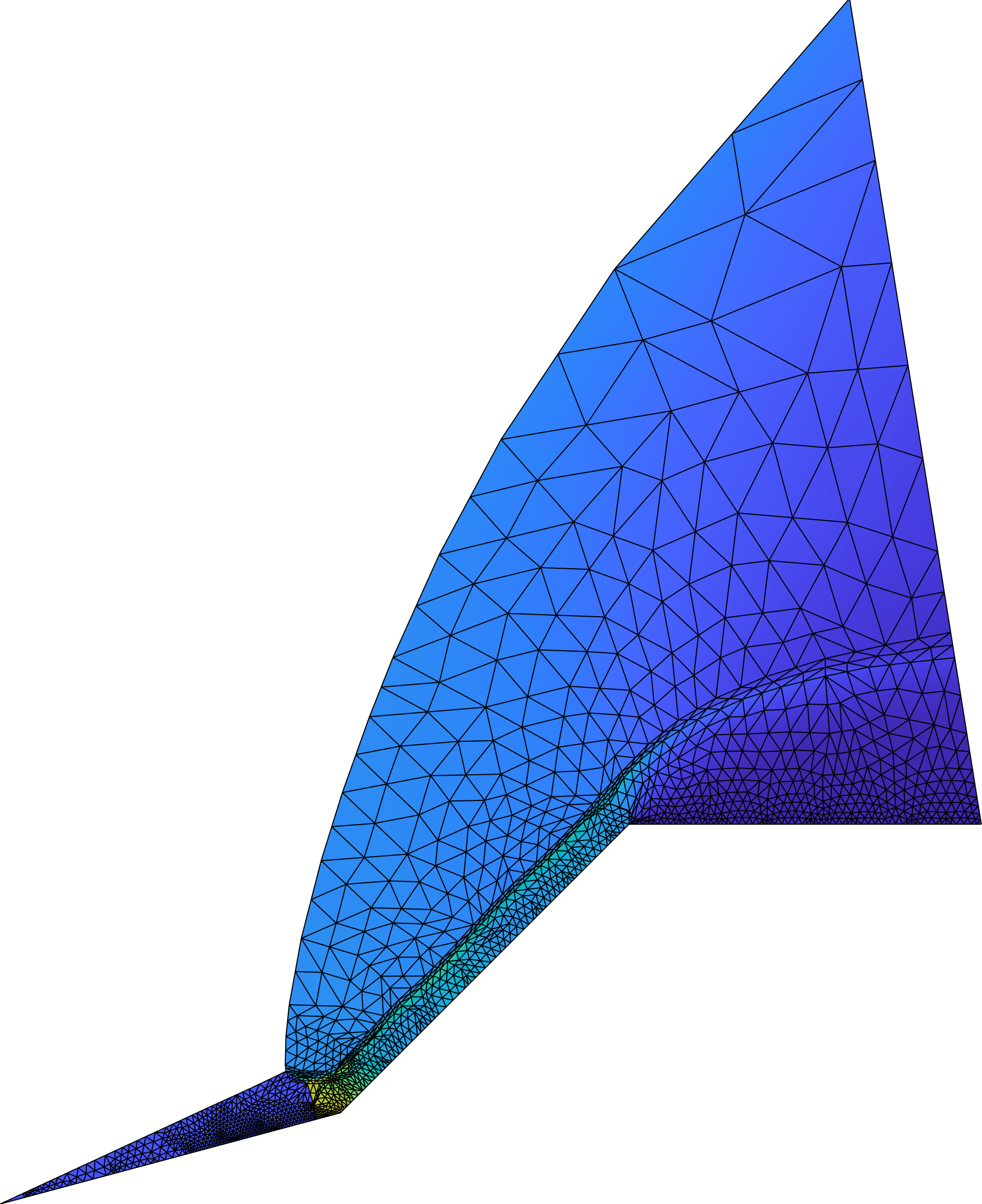}
}
\colorbarMatlabParula{0.67}{5}{10}{15}{20.85}
 \caption{
  Density distributions at $M_\infty = 4.77$ (stage 60) for double-wedge simulation:
  before re-mesh (\emph{left}), after re-mesh with a piecewise constant
  initial guess for the flow solution (\emph{middle-left}), the HOIST
  solution after re-mesh (\textit{middle-right}), and the HOIST solution
  after refinement of new mesh (\textit{right}).
 }
\label{fig:dw:remesh}
\end{figure}

Due to the nature of this problem, the elements condense near the intersection of
the two wedges at $(x_1, x_2) = (0.98, 0.13)$ forming long and skinny triangles as
the Mach number increases. These unavoidably crowded elements can degrade the solution
quality, which necessitates re-meshing at an intermediate parameter stage. We choose
to re-mesh at the middle of the parameter sweep without any attempt to preserve internal
shock structures or even the shock boundary. After the re-mesh, we use the HOIST method
with a free shock boundary to re-solve the stage to recover the appropriate lead shock
position and internal shock structure. As shown in Figure~\ref{fig:dw:remesh}, the initial
re-meshed configuration of 837 elements is shock-agnostic with respect to the shock
interactions downstream the lead shock, and the tracked configuration of 809 elements
resolves those solution features. 

\begin{figure}[!htbp]
\centering
\ifbool{fastcompile}{}{
\includegraphics[height=0.3\textwidth]{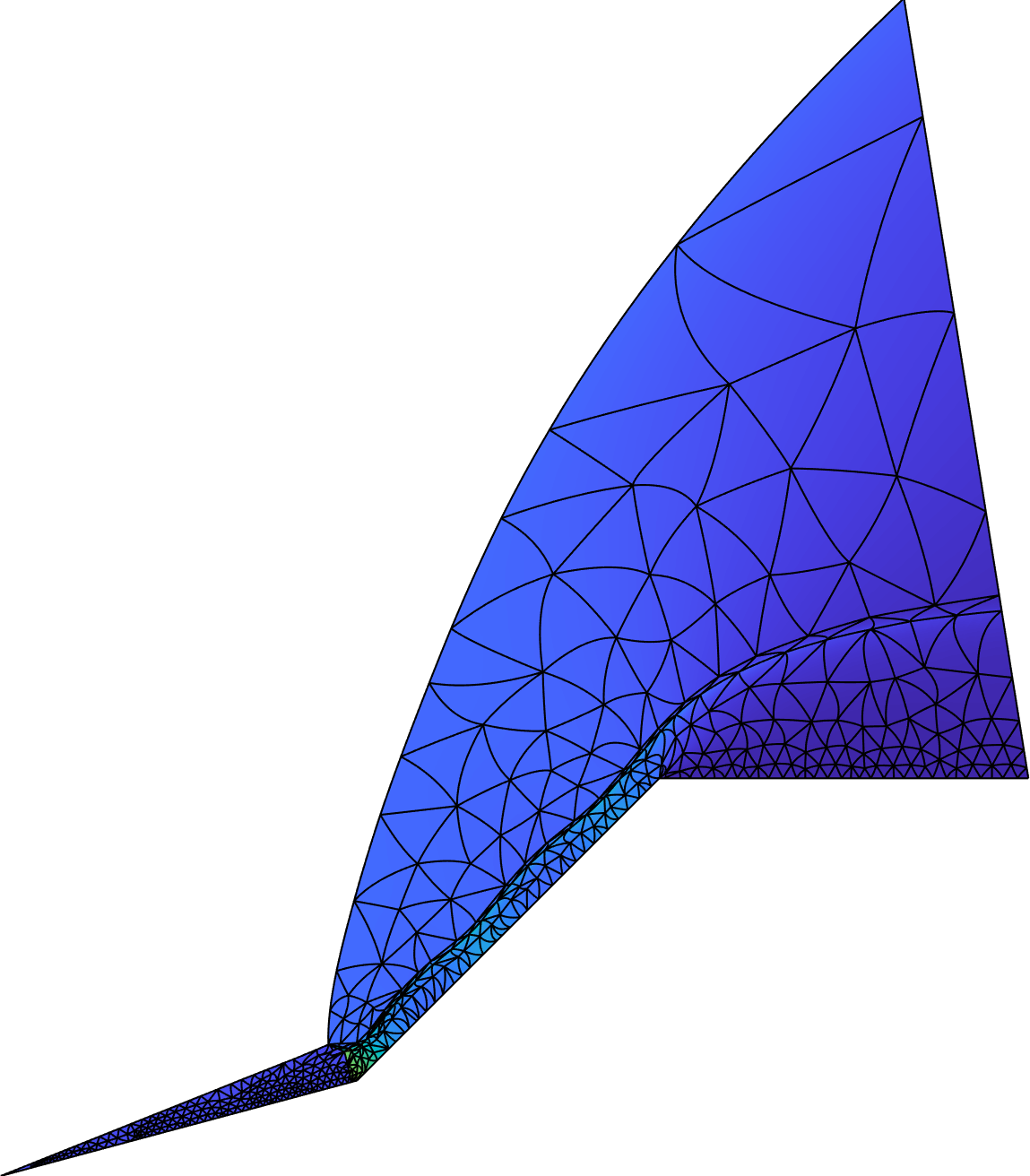} 
\includegraphics[height=0.3\textwidth]{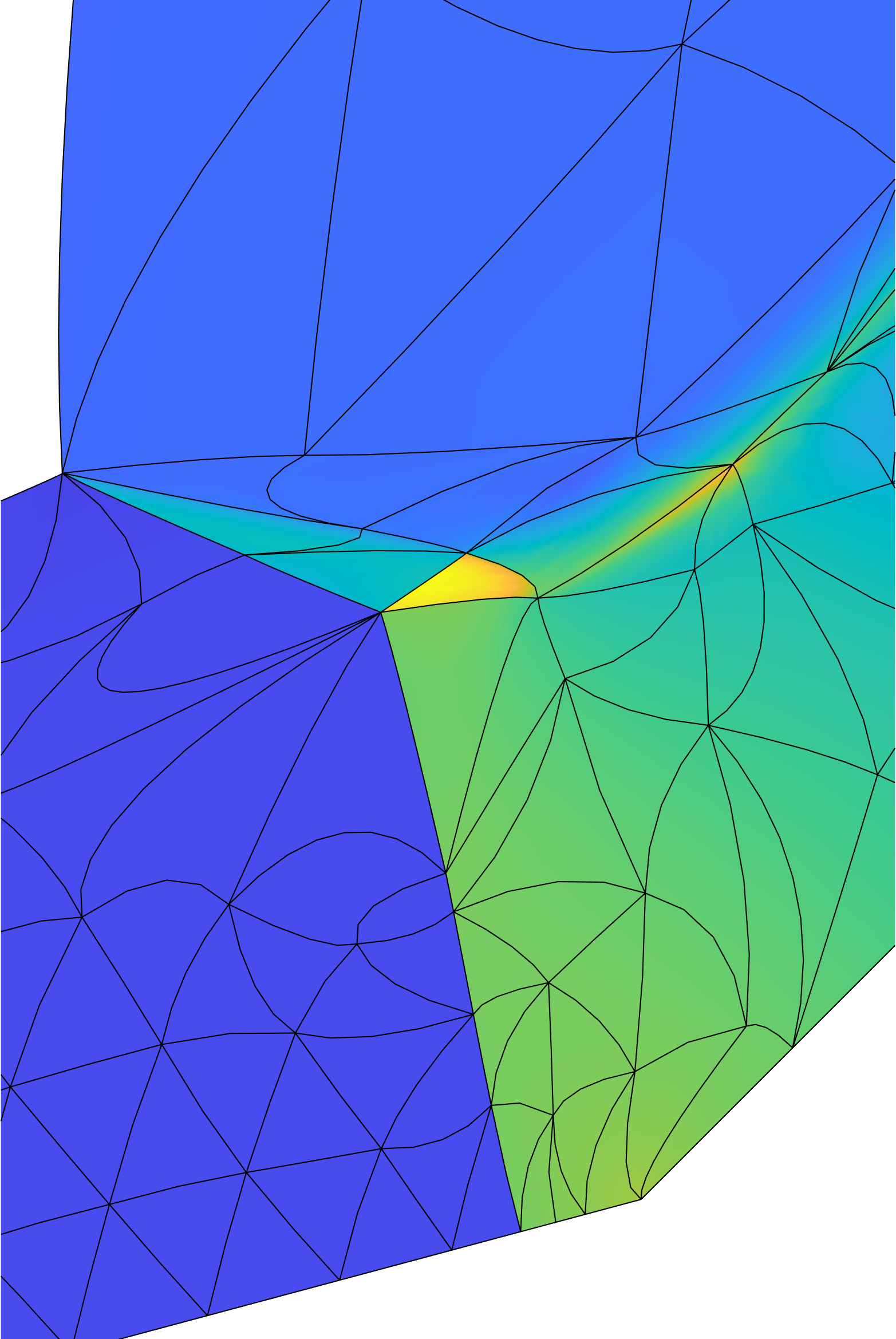} 
\includegraphics[height=0.3\textwidth]{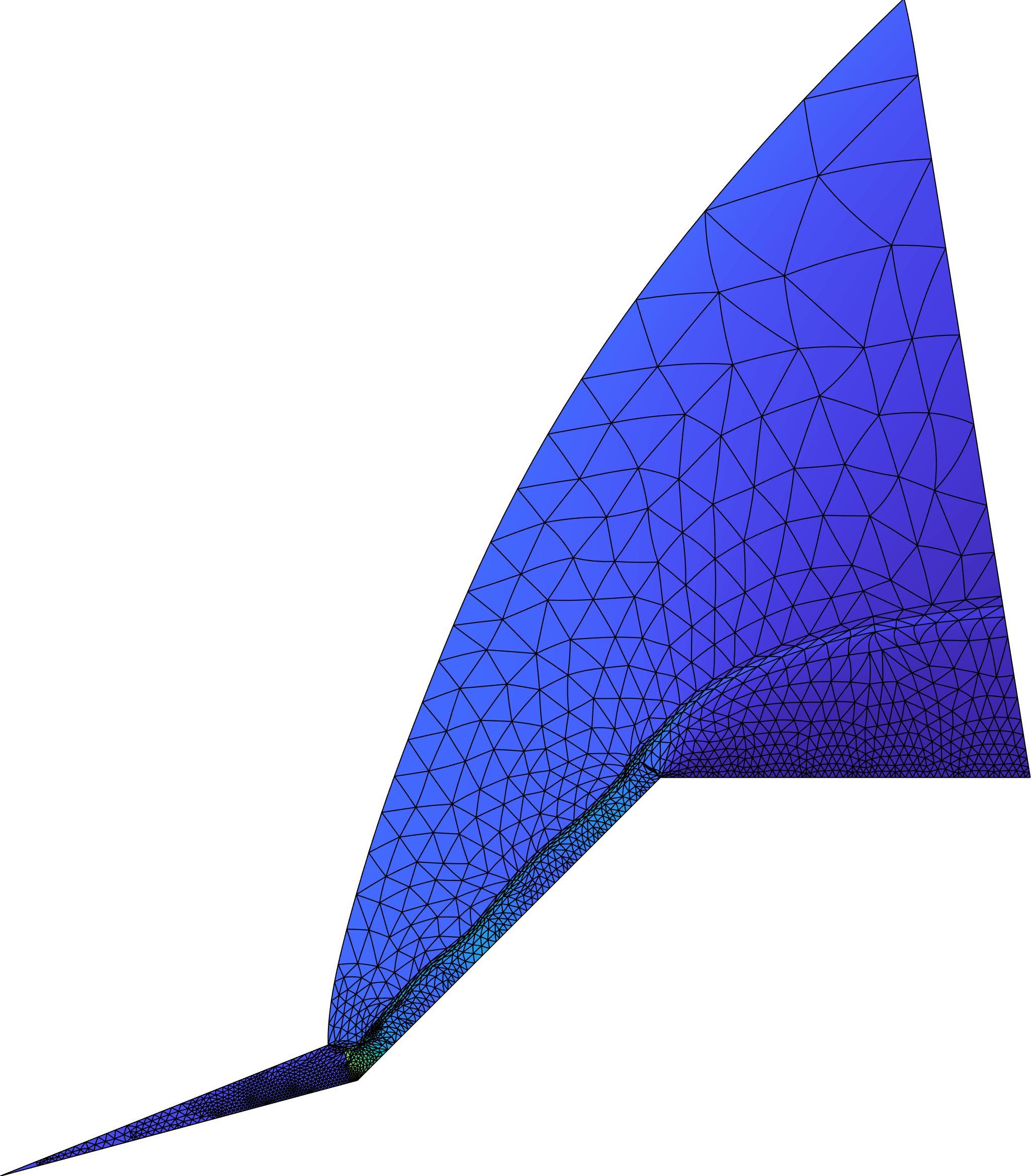} 
\includegraphics[height=0.3\textwidth]{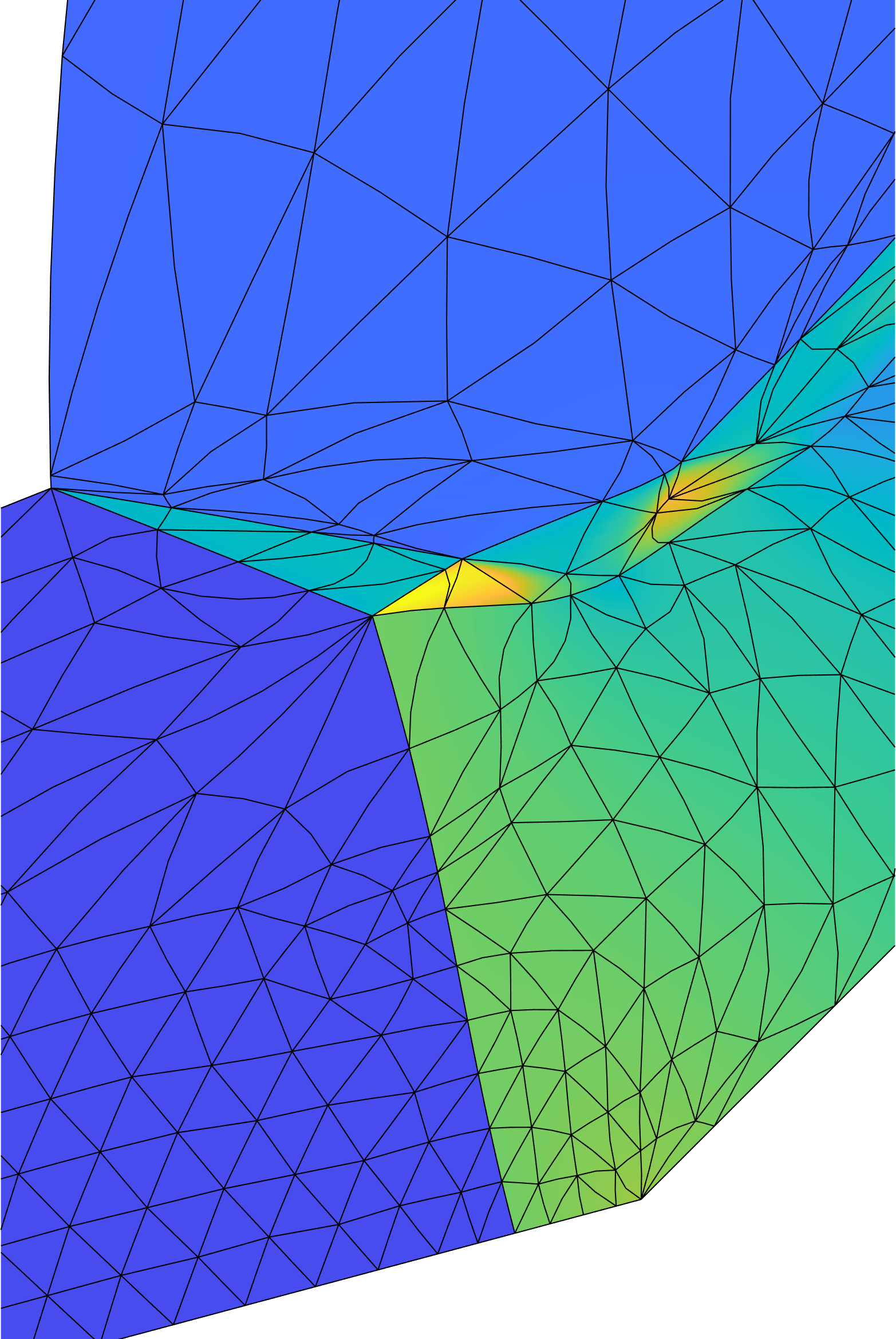}  \\
\includegraphics[height=0.3\textwidth]{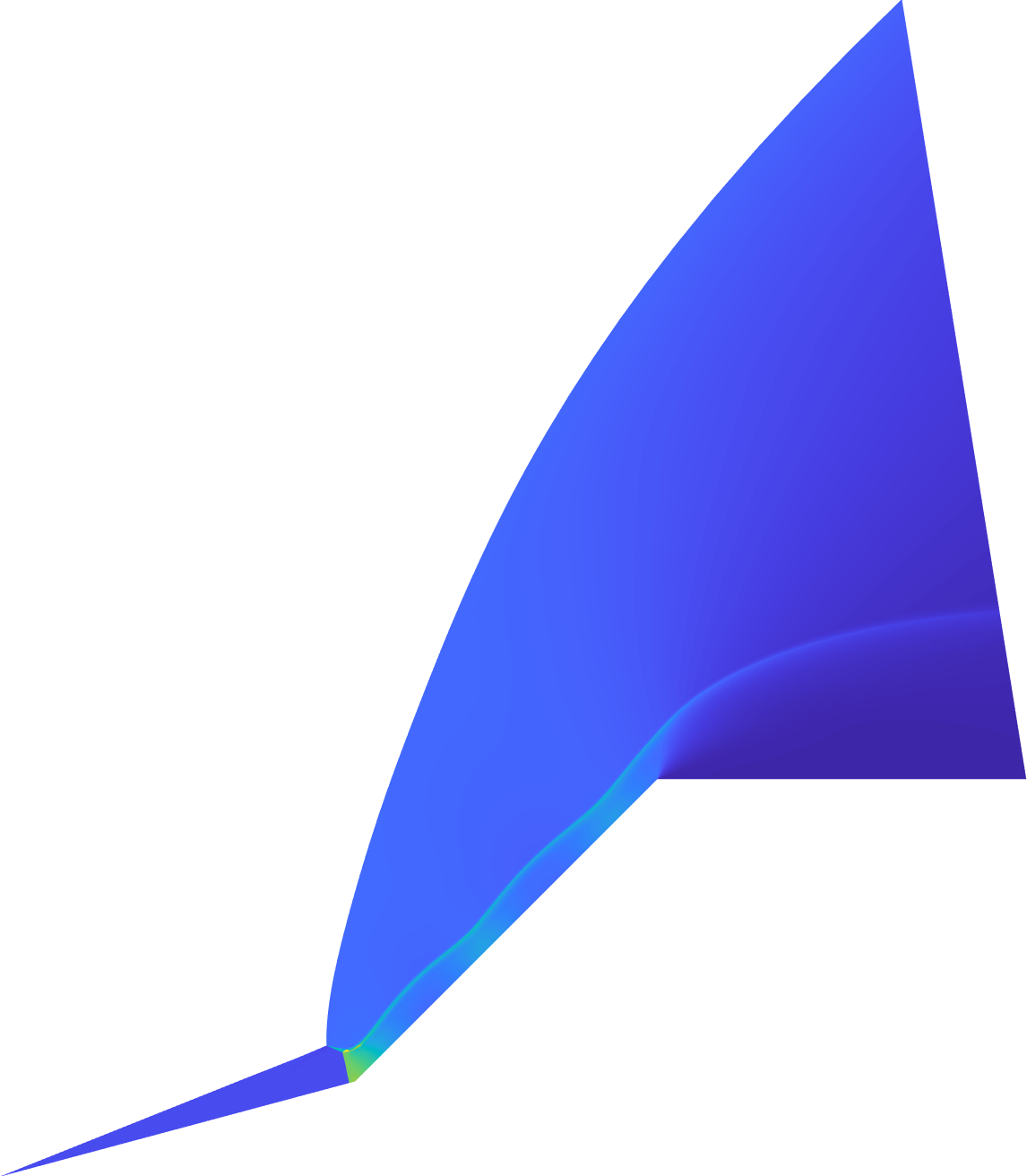} 
\includegraphics[height=0.3\textwidth]{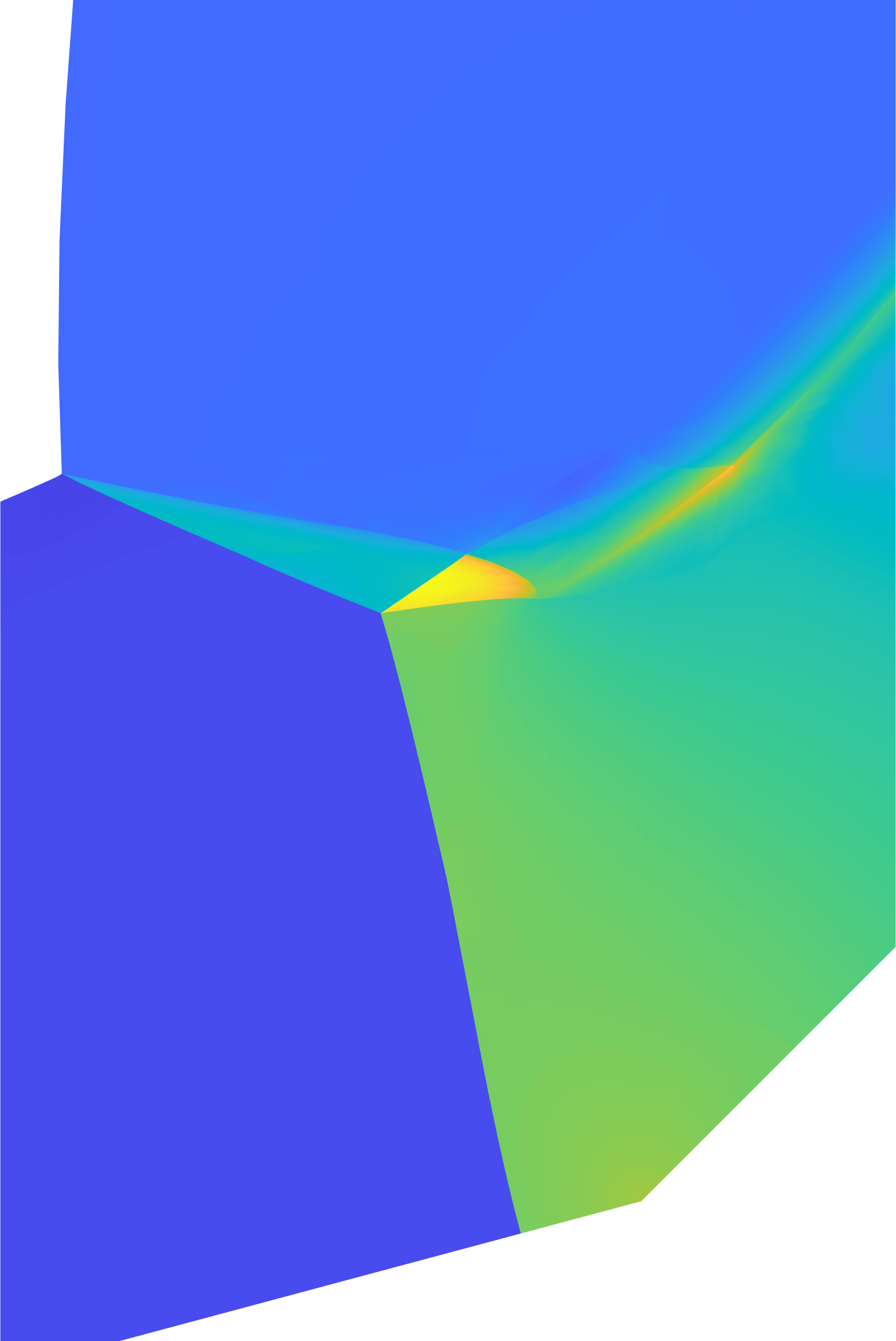} 
\includegraphics[height=0.3\textwidth]{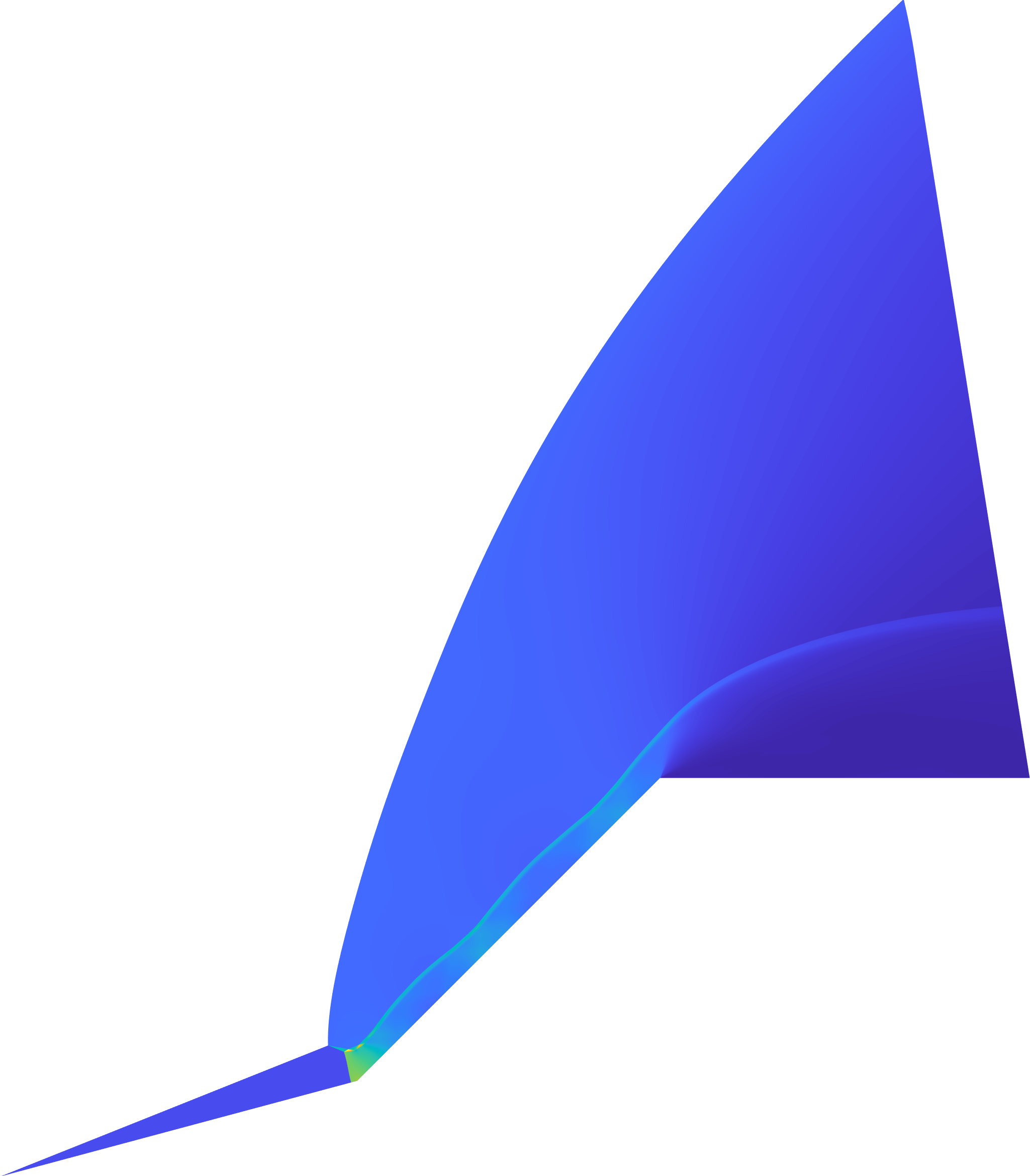}  
\includegraphics[height=0.3\textwidth]{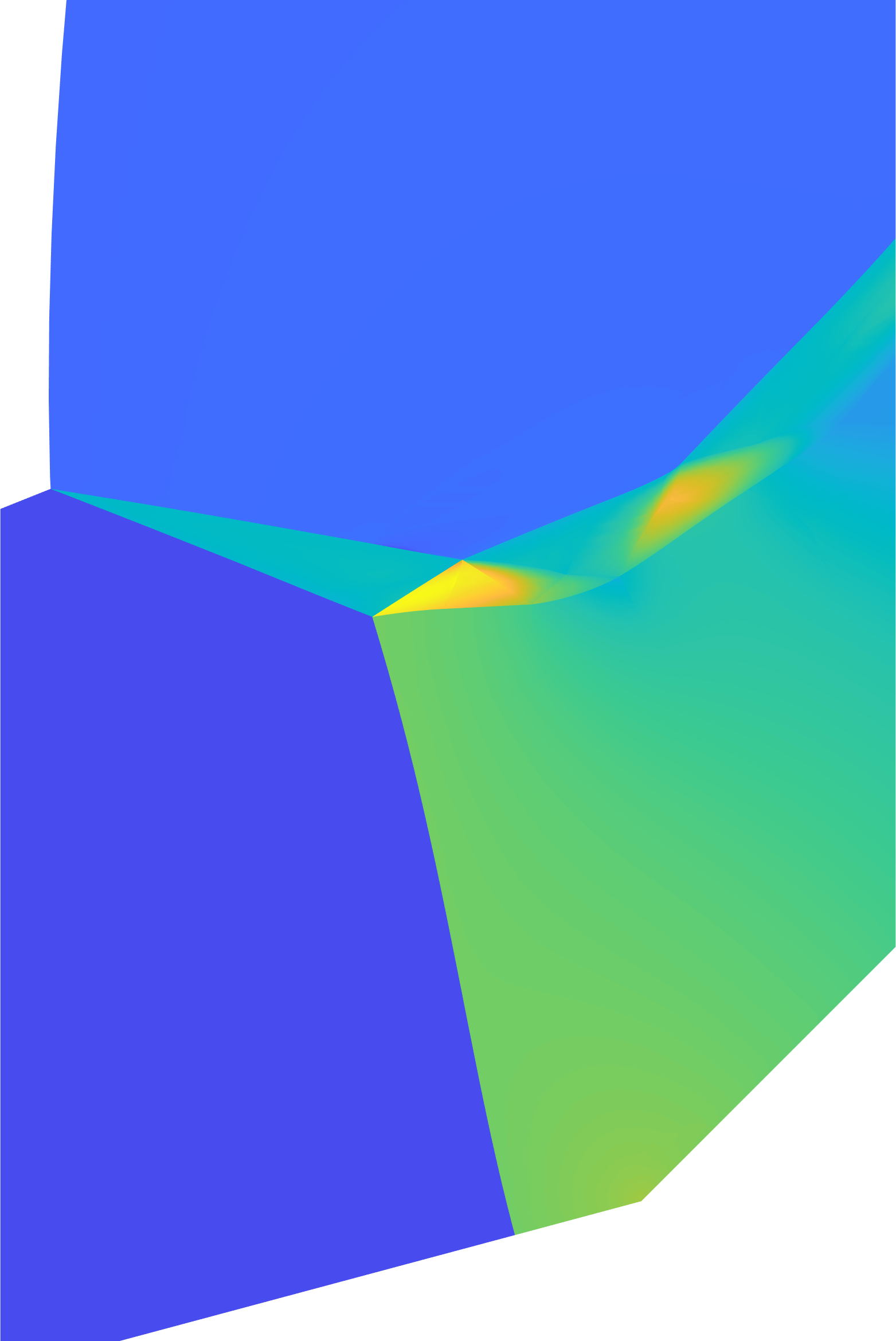}
}
 \caption{
 Density distribution at $M_\infty = 6.8$ for double-wedge simulation:
 coarse mesh (\textit{left}), coarse mesh with shock interaction blown up
 (\textit{middle-left}), refined mesh (\textit{right-middle}), and refined
 mesh with shock interaction blown up (\textit{right}). Colorbar in
 Figure~\ref{fig:dw:coarse_sweep}.
 }
\label{fig:dw:compare_cr}
\end{figure}

Despite the lead shock being accurately tracked, we observe the complex 
shock interactions have a low resolution and the contact discontinuity is not 
well-tracked on the coarse mesh as shown in Figure~\ref{fig:dw:compare_cr}. 
This is due to the low mesh element density in the shock interaction region.
As such, we refine the mesh at stage 60 (immediately after the re-mesh and solve),
which leads to a grid with 2846 elements, and use the HOIST method to compute the
corresponding refined solution (Figure~\ref{fig:dw:remesh}). From this mesh and solution,
the Mach continuation strategy continues until $M_\infty = 6.8$ is reach with
$n_1 = 30$. While the lead shock and solution are visually indistinguishable,
a close look at the shock interaction region shows the refined mesh tracks and
resolves the Type IV interaction and contact much better. Furthermore, the
refined mesh compresses more elements into the supersonic jet to better resolve it.
Finally, we show the positions of the lead shock at selected parameter stages
(Figure~\ref{fig:dw_shkpos}) and note the re-mesh at stage 60 did not disrupt
the lead shock location because the shock positions overlap (1) before and after
the re-mesh at stage 60 and (2) on the coarse and fine mesh at $M_\infty = 6.8$.

To close this study, we justify the reduced-mesh Mach continuation approach
relative to directly applying the HOIST method to the parameter of interest.
For this, we return to the full domain and use HOIST to directly solve the problem
at $M_\infty = 5.8$ (stage 91). In this setting, we observe poorly condition elements
upstream of the lead shock, numerous element collapses near the beginning of the first
wedge, an inaccurate lead shock position, and under-resolved shock interactions. These
issues can be avoided by either using a refined mesh or the proposed continuation approach
(Figure~\ref{fig:dw:coarse_sweep}).

\ifbool{fastcompile}{}{
\begin{figure}[!htbp]
\centering
\includegraphics{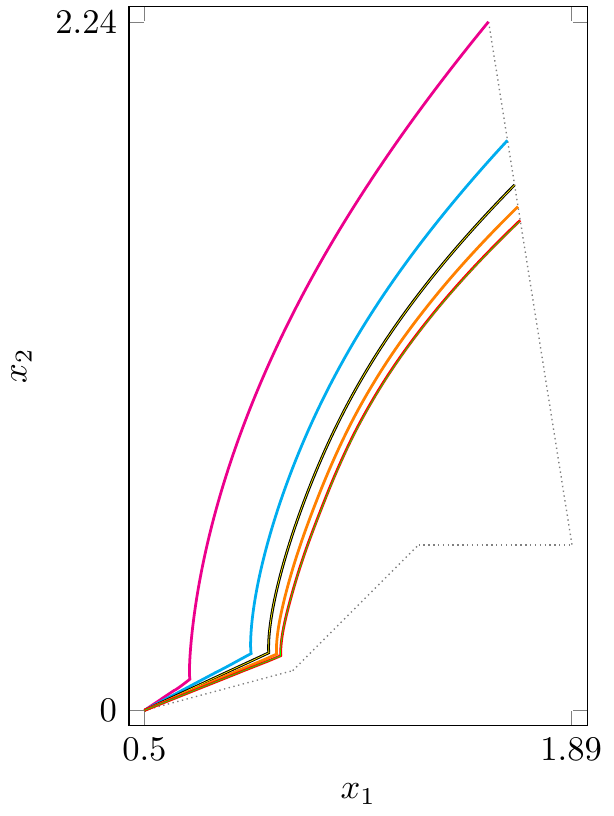}
\caption{
Lead shock positions for the double wedge simulation:
$M_\infty = 2.8$ 
(\raisebox{2pt}{\protect\tikz \protect\draw[magenta, thick] (0pt,100pt) -- (10pt,100pt);}), 
$M_\infty = 3.77$
(\raisebox{2pt}{\protect\tikz \protect\draw[cyan, thick] (0pt,100pt) -- (10pt,100pt);}), 
$M_\infty = 4.77$ (stage 60 before re-mesh )
(\raisebox{2pt}{\protect\tikz \protect\draw[black, thick] (0pt,100pt) -- (10pt,100pt);}), 
$M_\infty = 4.77$ (stage 60 after re-mesh)
(\raisebox{2pt}{\protect\tikz \protect\draw[yellow, very thin] (0pt,100pt) -- (10pt,100pt);}), 
$M_\infty = 5.77$
(\raisebox{2pt}{\protect\tikz \protect\draw[orange, thick] (0pt,100pt) -- (10pt,100pt);}), 
$M_\infty = 6.8$ (final stage, coarse mesh)
(\raisebox{2pt}{\protect\tikz \protect\draw[red, thick] (0pt,100pt) -- (10pt,100pt);}), 
$M_\infty = 6.8$ (final stage, fine mesh) 
(\raisebox{2pt}{\protect\tikz \protect\draw[green, thick] (0pt,100pt) -- (10pt,100pt);}).
 }
\label{fig:dw_shkpos}
\end{figure}
}

\begin{figure}[!htbp]
\centering
\ifbool{fastcompile}{}{
\includegraphics[width=0.35\textwidth]{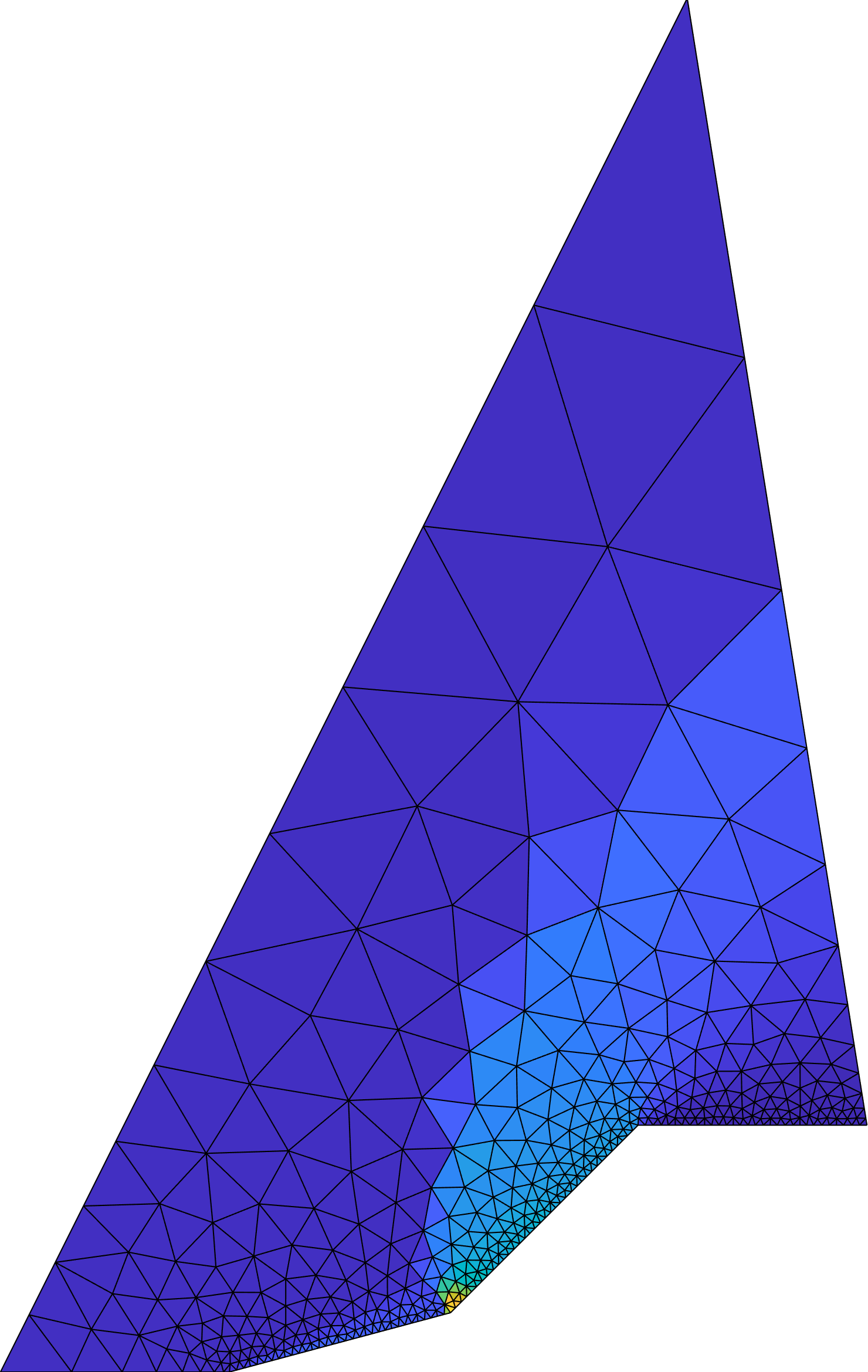} 
\includegraphics[width=0.35\textwidth]{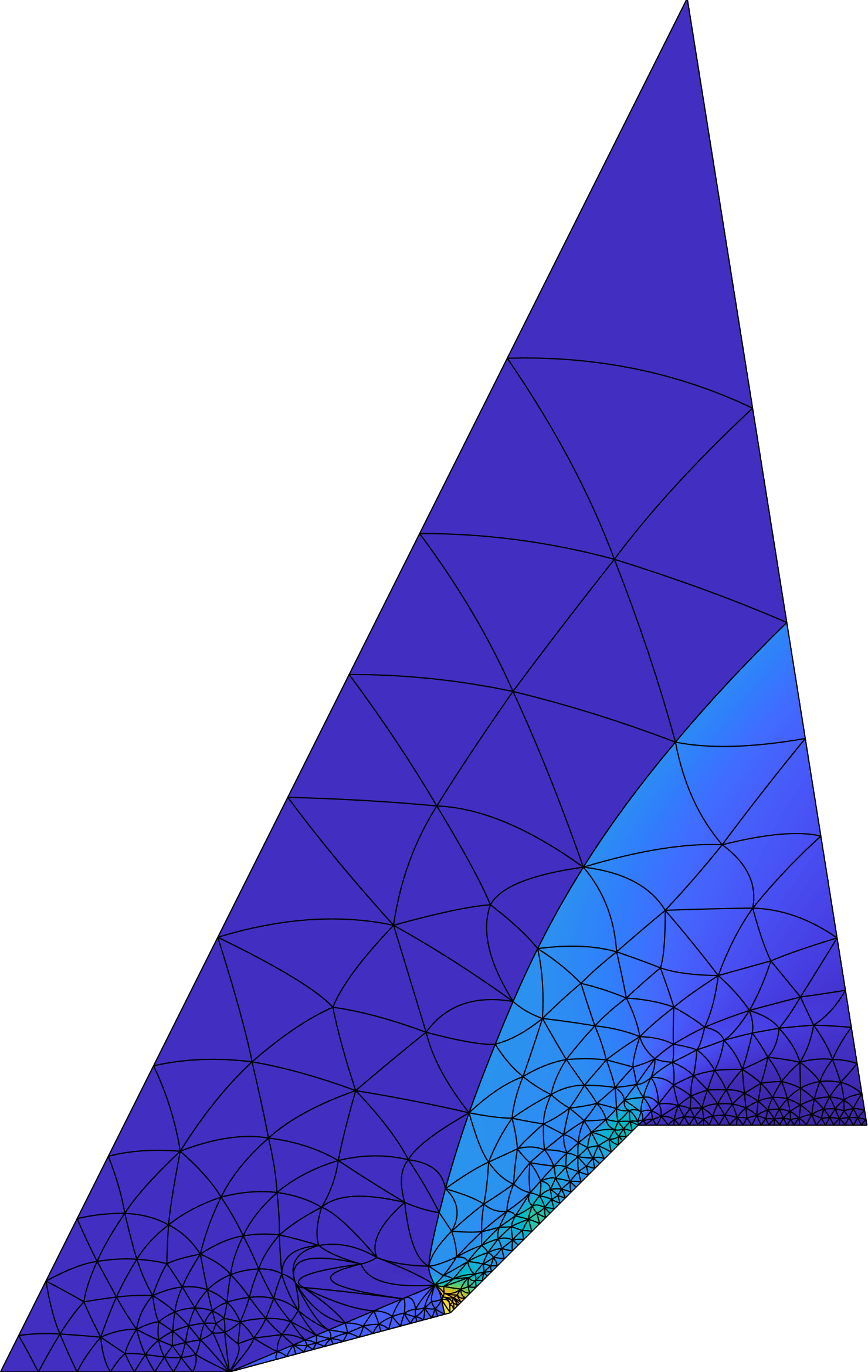} 
}
\colorbarMatlabParula{0.66}{5}{10}{15}{20.87}
 \caption{
 Initial guess (density) on a shock-agnostic mesh (full domain) of
 $M_\infty=5.8$ flow over the double-wedge geometry (\textit{left})
 and the HOIST solution (\textit{right}). Mesh motion from the
 shock-upstream portion of the domain introduces poorly conditioned
 elements.
 }
\label{fig:dw:direct_s91}
\end{figure}

\subsection{Initial condition sweep}
\label{sec:numexp:ic}
To demonstrate the flexibility of the proposed many-query framework beyond Mach
continuation, we close by considering a parametrized Riemann problem of the
Euler equations. The Euler equations model unsteady, compressible flow in a
one-dimensional domain $\Omega \subset\Rbb$ and read
\begin{equation} \label{eqn:ueuler}
\begin{split}
\pder{}{t}\rho(x,t) + \pder{}{x}\left(\rho(x,t) v(x,t)\right) &= 0, \\
\pder{}{t}\left(\rho(x,t)v(x,t)\right) + \pder{}{x}\left(\rho(x,t)v(x,t)^2+P(x,t)\right) &= 0, \\
\pder{}{t}\left(\rho(x,t)E(x,t)\right) + \pder{}{x}\left(\left[\rho(x,t)E(x,t)+P(x,t)\right]v(x,t)\right) &= 0,
\end{split}
\end{equation}
for all $x\in\Omega$ and $t\in\Tcal$, where $\Tcal\coloneqq\coloneqq(0,T]$
is the temporal domain and $T \in \Rbb_{>0}$ is the final time. The density 
$\func{\rho}{\Omega\times\Tcal}{\Rbb_{>0}}$, velocity
$\func{v}{\Omega\times\Tcal}{\Rbb}$, and total energy 
$\func{E}{\Omega\times\Tcal}{\Rbb_{>0}}$ of the fluid are
implicitly defined as the solution of (\ref{eqn:ueuler}). 
We assume the fluid is an ideal gas, which leads to the
relationship between pressure and energy in (\ref{eqn:EtoP}).
Finally, the Euler equations can be written as a general system of steady
conservation laws (\ref{eqn:claw-phys}) over the space-time domain
$\Omega\times\Tcal$ with space-time coordinate $z = (x,t)$ as
\begin{equation*}
 U(z) = \begin{bmatrix} \rho(x,t) \\ \rho(x,t) v(x,t) \\ \rho(x,t) E(x,t) \end{bmatrix}, \qquad
 F(U(z)) = \begin{bmatrix} \rho(x,t)v(x,t) & \rho(x,t) \\ \rho(x,t) v(x,t)^2 + P(x,t) & \rho(x,t)v(x,t) \\ (\rho(x,t)E(x,t)+P(x,t)) v(x,t) & \rho(x,t)E(x,t) \end{bmatrix}, \qquad
 S(U(z)) = 0.
\end{equation*}

We consider a Riemann problem of the Euler equations (\ref{eqn:ueuler}) with the
following parametrized initial condition
\begin{equation}
 (\rho(x,0), v(x,0), P(x,0)) =
 \begin{cases}
  (\rho_L, 0, P_L(\mu)) & x < 0.5 \\
  (\rho_R(\mu), 0, P_R(\mu)) & x \geq 0.5,
 \end{cases}
\end{equation}
where $\mu \in \Dcal \coloneqq [0, 1]$ and
\begin{equation}
 P_L(\mu) = 25 \mu + (1-\mu), \qquad
 \rho_R(\mu) = 0.35 \mu + 0.125 (1-\mu), \qquad
 P_R(\mu) = 0.075 \mu + 0.1 (1-\mu)
\end{equation}
over the spatial domain $\Omega \coloneqq (0, 1)$ and temporal domain $\Tcal = (0, 0.2]$.
At $\mu = 0$, this corresponds to the canonical Sod shock tube and possesses
similar features to the Woodward-Colella blast \cite{woodward1984numerical}
at $\mu = 1$. The discontinuity magnitudes and propagation speed of the
waves varies significantly for $\mu \in \Dcal$.
Because larger values of $\mu$ will have faster waves, we allow the final
time to vary with $\mu$ such that the waves travel approximately the same
distance in $\Tcal$. That is, we let $T = T(\mu)$ with $T(0) = 0.2$ and
final time for other parameters determined using the HOIST framework.

\begin{figure}[!htbp]
\centering
\ifbool{fastcompile}{}{
\begin{tikzpicture}
\begin{groupplot}[
  group style={
      group size=1 by 3,
      vertical sep=0.2cm
  },
  width=0.8\textwidth,
  axis equal image,
  xlabel={$x$},
  ylabel={$t$},
  xtick = {0, 0.2647, 0.8504, 1},
  ytick = {0, 0.2},
  xmin=0, xmax=1,
  ymin=0, ymax=0.2,
  yticklabel pos=left,
  xticklabel style={/pgf/number format/.cd,fixed zerofill,precision=2},
  y label style={at={(1.05, 0.5)}}
]
\nextgroupplot[xlabel={}, yticklabels={,,}, xticklabels={,,}]
\addplot graphics [xmin=0,xmax=1,ymin=0,ymax=0.2] {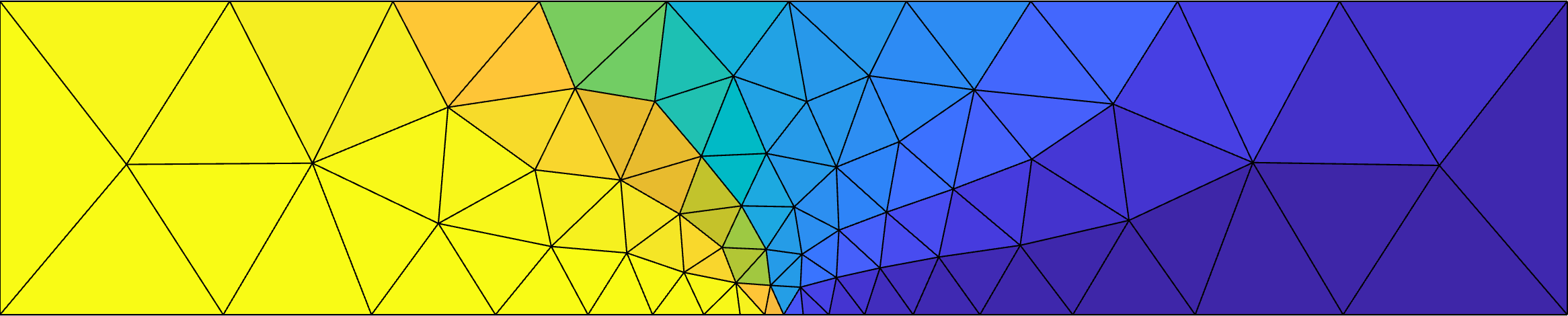};
\nextgroupplot[xlabel={}, yticklabels={,,}, xticklabels={,,}]
\addplot graphics [xmin=0,xmax=1,ymin=0,ymax=0.2] {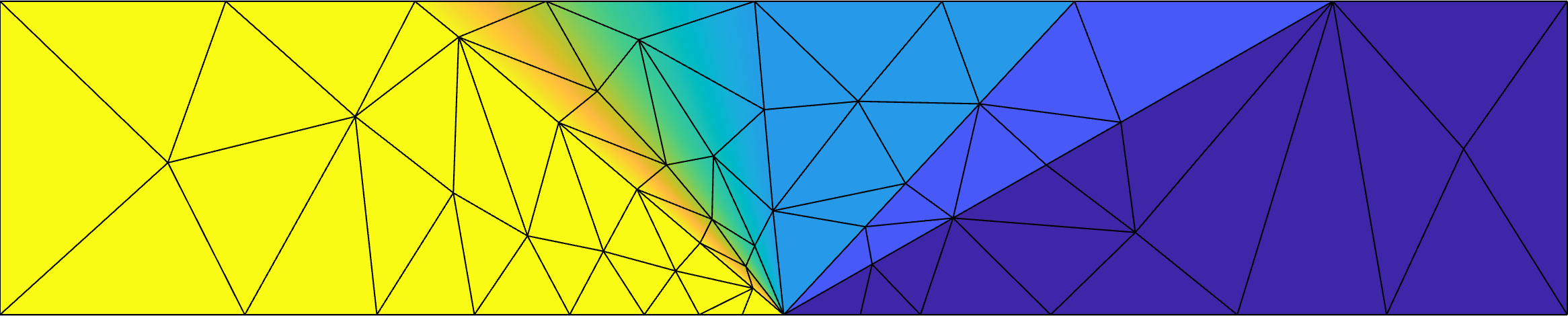};
\nextgroupplot[xticklabels={0, 0.2647, 0.8504, 1}]
\addplot graphics [xmin=0.2647,xmax=0.8504,ymin=0,ymax=0.2] {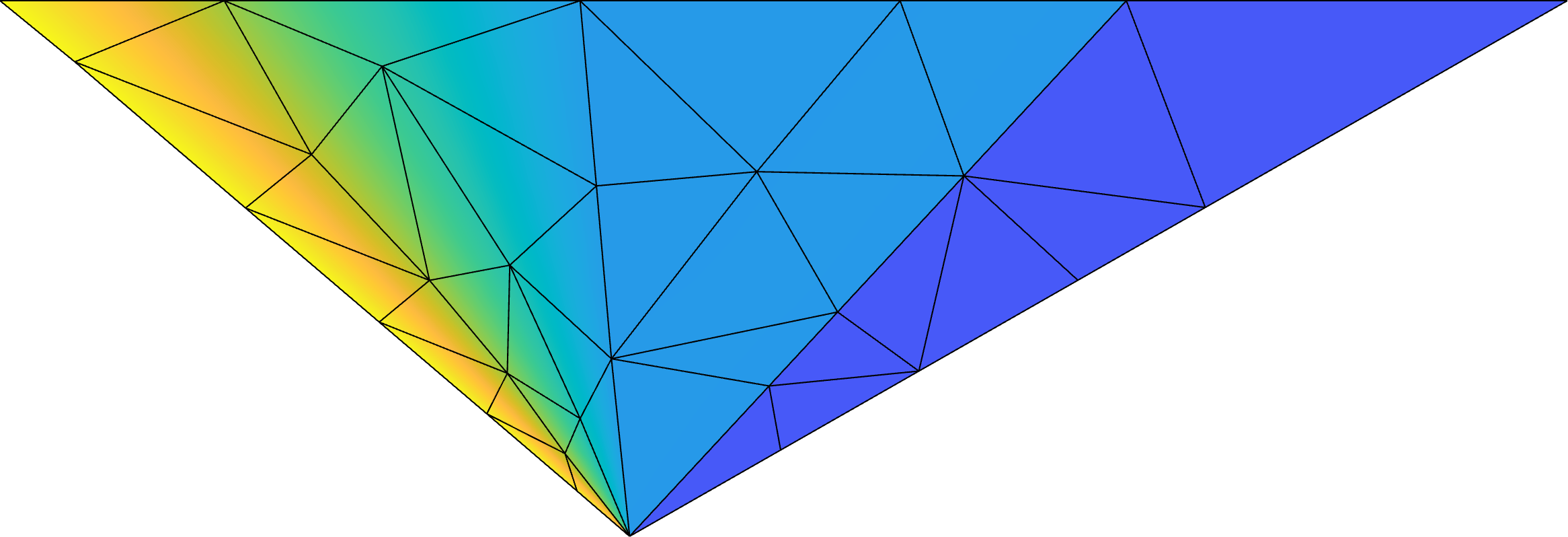};
\addplot [solid, gray, dashed]
coordinates {
( 0.2647,  0)
( 0.2647,  0.2)};
\addplot [solid, gray, dashed]
coordinates {
( 0.8504,  0)
( 0.8504,  0.2)};
\end{groupplot}
\end{tikzpicture}
\colorbarMatlabParula{0.125}{0.25}{0.5}{0.75}{1}
}
\caption{
  Initial guess (density) on a shock-agnostic space-time mesh for the Riemann problem
  at $\mu = 0$ (\textit{top}), the shock-aligned mesh and corresponding solution obtained
  using HOIST (\textit{middle}), and the corresponding solution and mesh extracted from
  the region to the left of the rarefaction head and right of the shock (\textit{bottom}).
}
\label{fig:riemann_prob_stage0}
\end{figure}

Following the proposed procedure, we begin with a shock-agnostic
mesh of the space-time domain consisting of $82$ triangles and apply the HOIST
method at $\mu=0$ (Sod) to create a shock-aligned mesh and the corresponding solution.
We use linear ($q=1$) mesh elements and quadratic ($p=2$) solution approximation because
all non-smooth features are straight-sided for any $\mu\in\Dcal$ and the solution
inside the rarefaction is nonlinear. At $\mu = 0$, the solution possesses a shock wave,
contact discontinuity, and rarefaction wave, with both the shock wave and head of the
rarefaction qualifying as lead ``shocks'' (non-smooth features separating a boundary
from the remainder of the flow). Finally, we remove all elements left of the rarefaction
and right of the shock wave, except one layer of upstream elements are retained (not shown
in figure for clarity) to improve boundary condition enforcment, to create a reduced mesh
of $42$ elements (Figure~\ref{fig:riemann_prob_stage0}).

\begin{figure}[!htbp]
\centering
\begin{tikzpicture}
\begin{groupplot} [
group style={group size = 2 by 2, horizontal sep = 0.5cm, vertical sep = 0.75cm},
title style={at={(current bounding box.north west)}, anchor=west}]
\nextgroupplot[width=0.4\textwidth, xtick={}, ytick={0, 0.2}, xticklabels={}, yticklabels={0, 0.2}, xmin=0.2647, xmax=0.8504, ymin=0, ymax=0.2]
\addplot []
graphics [xmin=0.2648,xmax=0.8504,ymin=0,ymax=0.2] { _img/riemann_prob/sod_extracted_domain.png};

\nextgroupplot[width=0.4\textwidth, xtick={0.2647, 0.8504}, ytick={}, xticklabels={0.26, 0.85}, yticklabels={}, xmin=0.2647, xmax=0.8504, ymin=0, ymax=0.2]
\addplot []
graphics [xmin=0.2648,xmax=0.8504,ymin=0,ymax=0.2] { 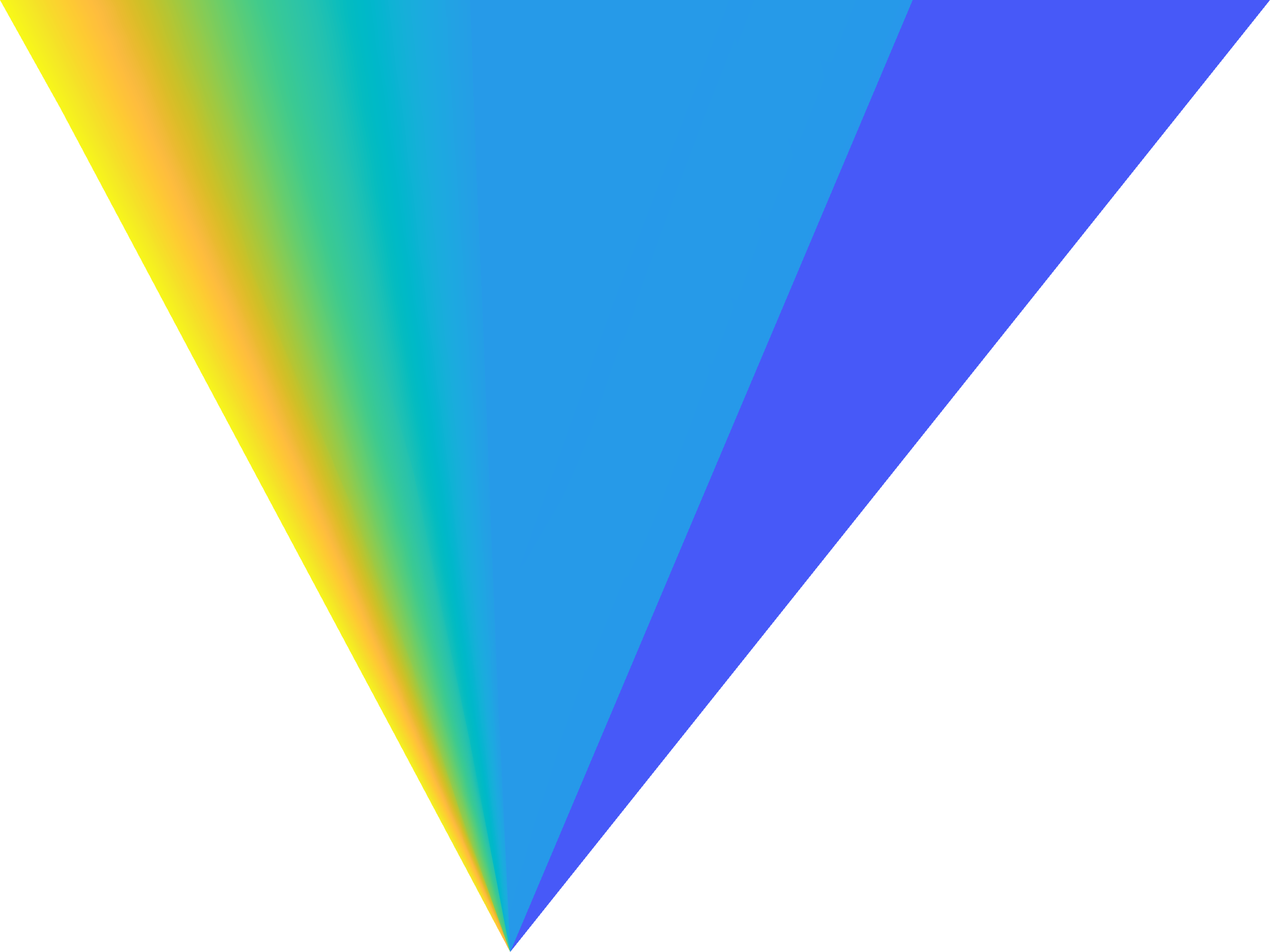};

\nextgroupplot[width=0.4\textwidth, xtick={}, ytick={0, 0.6}, xticklabels={}, yticklabels={0, 0.06}, xmin=0.2219, xmax=0.7744, ymin=0, ymax=0.6]
\addplot []
graphics [xmin=0.222,xmax=0.7744,ymin=0,ymax=0.6] { 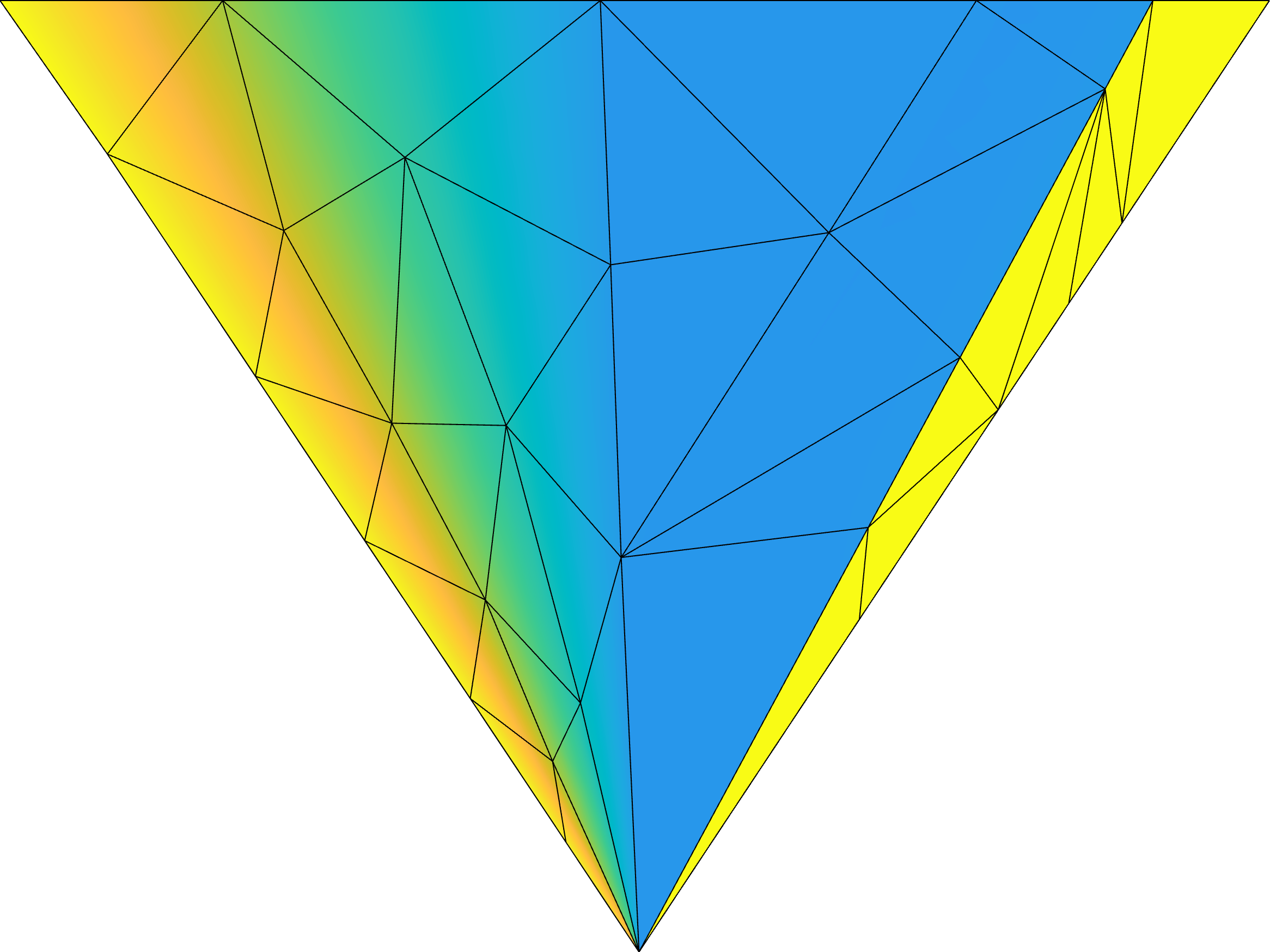};

\nextgroupplot[width=0.4\textwidth, xtick={0.2219, 0.7744}, ytick={}, xticklabels={0.22, 0.77}, yticklabels={}, xmin=0.2219, xmax=0.7744, ymin=0, ymax=0.6]
\addplot []
graphics [xmin=0.222,xmax=0.7744,ymin=0,ymax=0.6] { 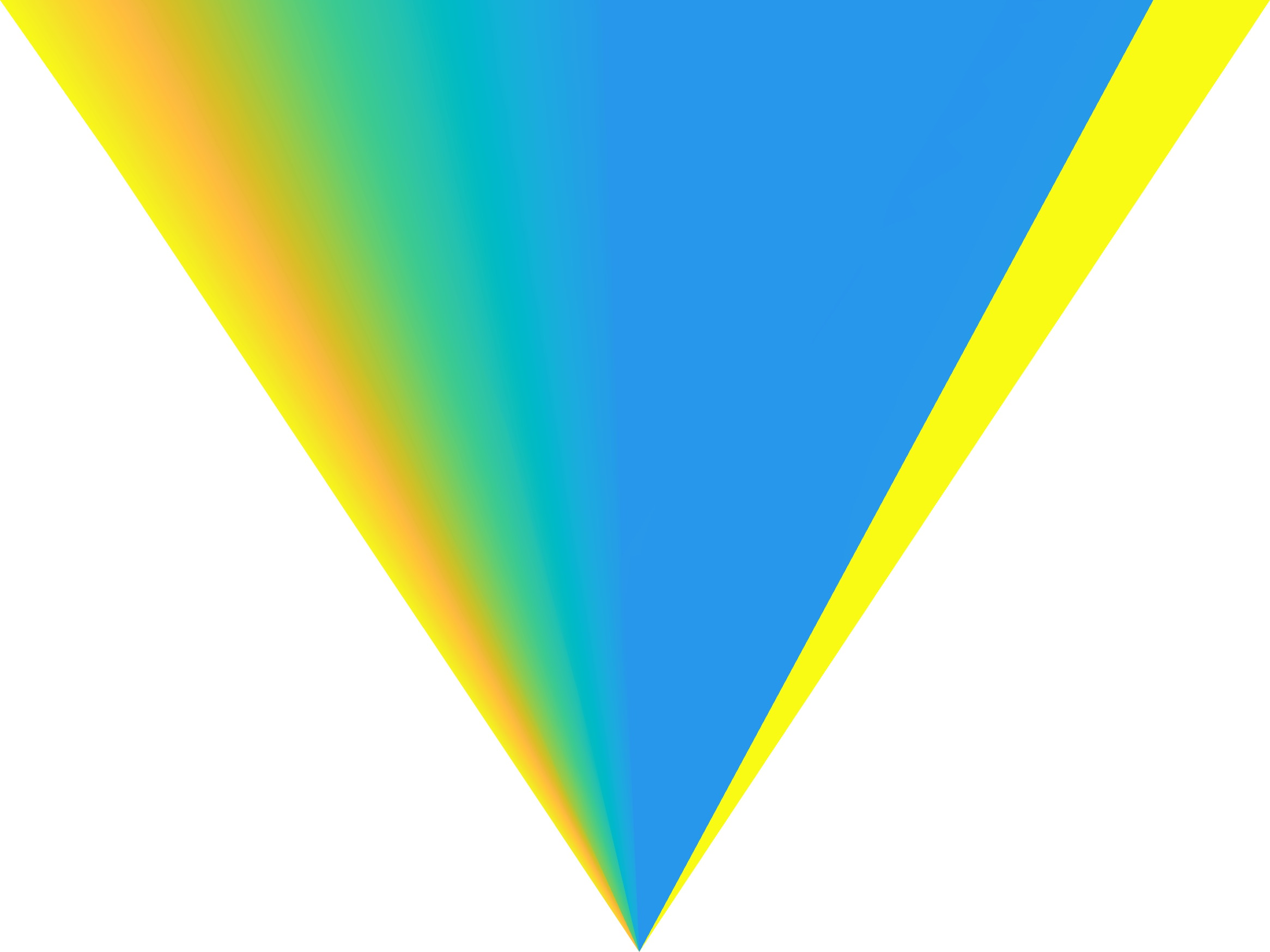};
\end{groupplot}\end{tikzpicture}
 \colorbarMatlabParula{0.125}{0.25}{0.5}{0.75}{1}

\vspace{3mm}
\begin{tikzpicture}
\begin{groupplot} [
group style={group size = 2 by 2, horizontal sep = 0.5cm, vertical sep = 0.75cm},
title style={at={(current bounding box.north west)}, anchor=west}]

\nextgroupplot[width=0.4\textwidth, xtick={0}, ytick={0, 0.477}, xticklabels={}, yticklabels={0, 0.0477}, xmin=0.2469, xmax=0.7374, ymin=0, ymax=0.477]
\addplot []
graphics [xmin=0.2469,xmax=0.7374,ymin=0,ymax=0.477] { 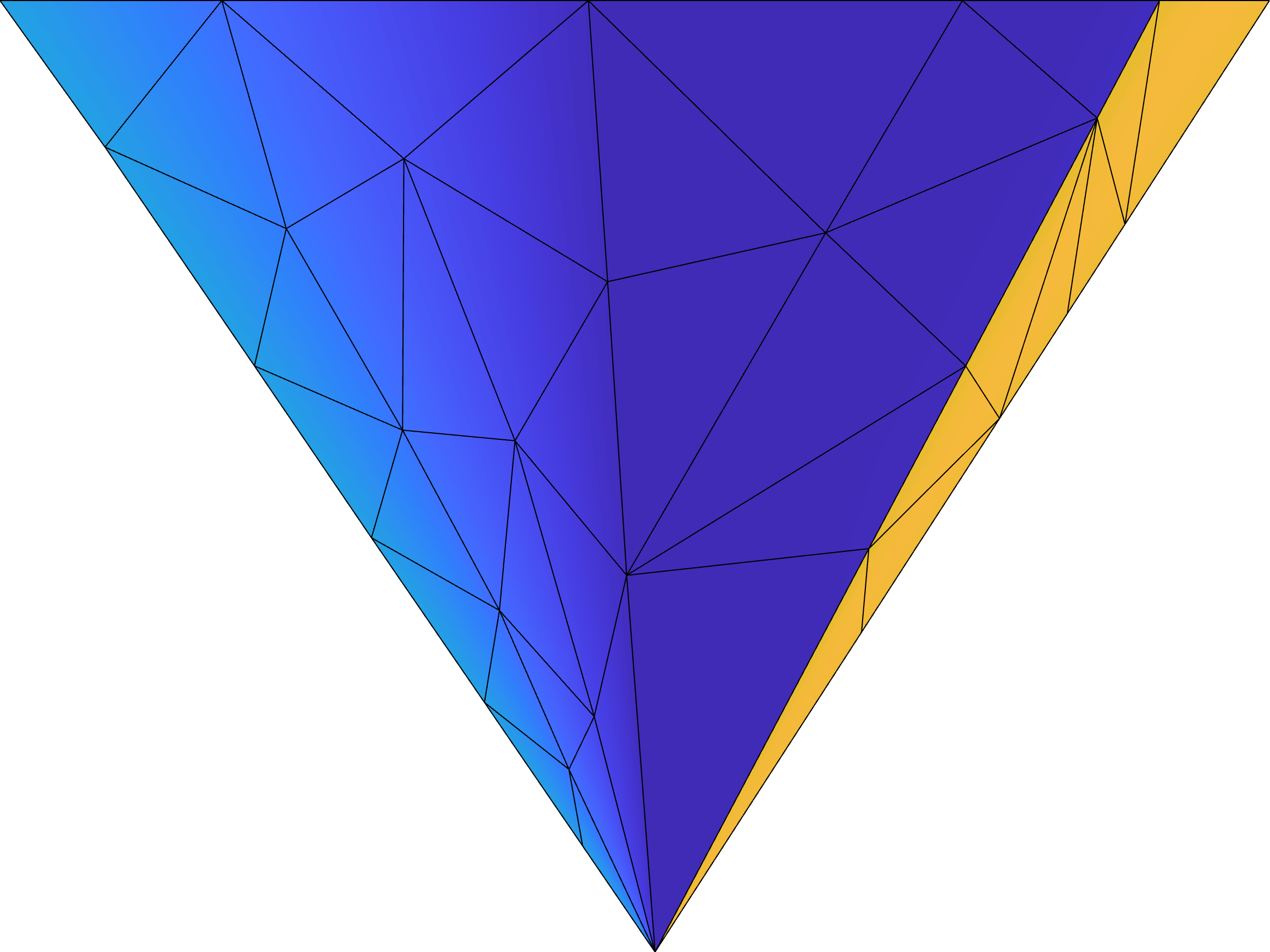};

\nextgroupplot[width=0.4\textwidth, xtick={0.2469, 0.7374}, ytick={0}, xticklabels={0.25, 0.74}, yticklabels={}, xmin=0.2469, xmax=0.7374, ymin=0, ymax=0.477]
\addplot []
graphics [xmin=0.2469,xmax=0.7374,ymin=0,ymax=0.477] { 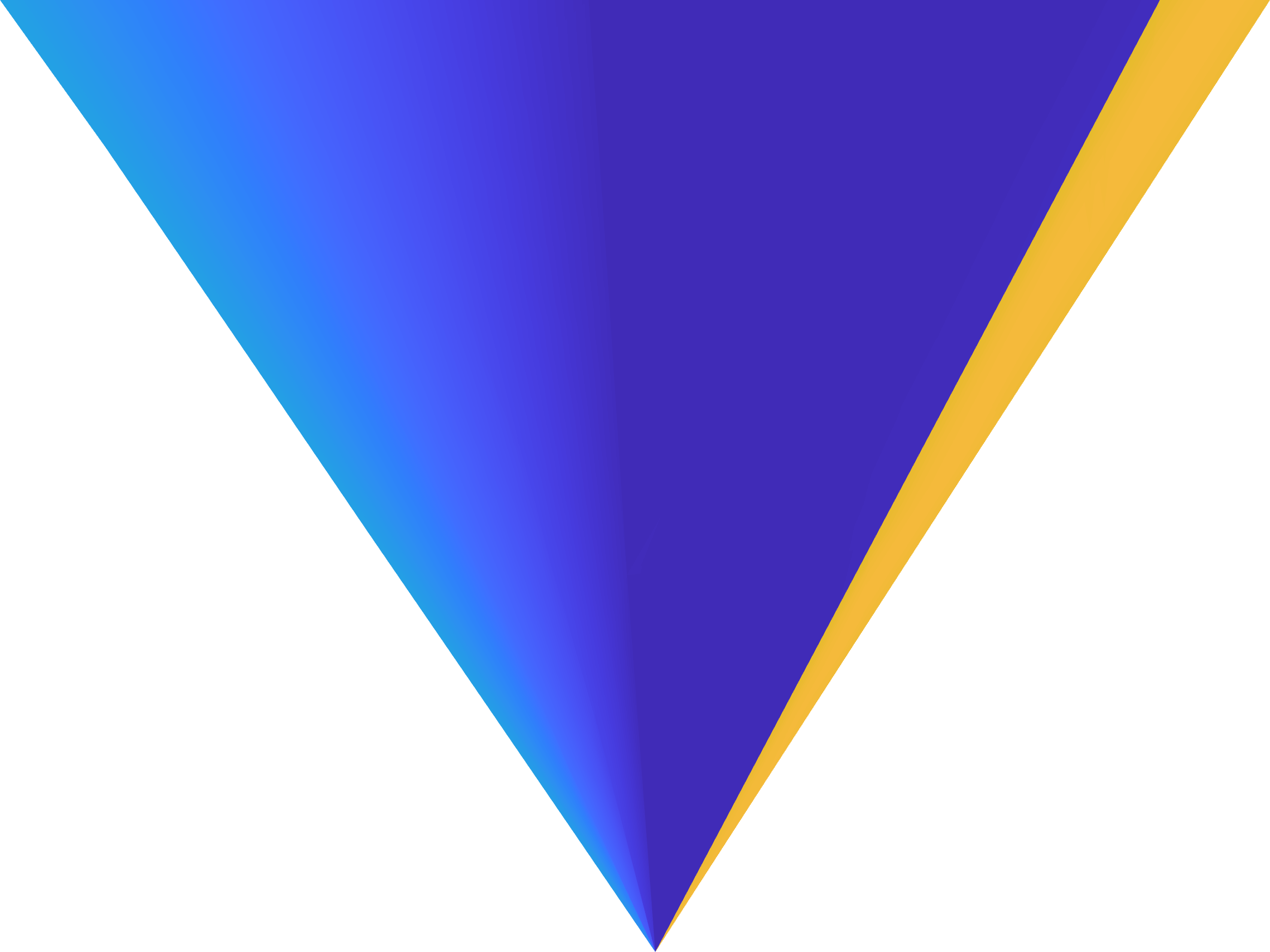};

\nextgroupplot[width=0.4\textwidth, xtick={}, ytick={0, 0.39}, xticklabels={}, yticklabels={0, 0.039}, xmin=0.2746, xmax=0.7066, ymin=0, ymax=0.39]
\addplot []
graphics [xmin=0.2746,xmax=0.7066,ymin=0,ymax=0.39] { _img/riemann_prob/final_State_extracted_domain_msh.png};

\nextgroupplot[width=0.4\textwidth, xtick={0.2746, 0.7066}, ytick={}, xticklabels={0.27, 0.7}, yticklabels={}, xmin=0.2746, xmax=0.7066, ymin=0, ymax=0.39]
\addplot []
graphics [xmin=0.2746,xmax=0.7066,ymin=0,ymax=0.39] { _img/riemann_prob/final_State_extracted_domain_nomsh.png};

\end{groupplot}\end{tikzpicture}
 \colorbarMatlabParula{0.4}{0.8}{1.2}{1.6}{2}
 \caption{Space-time density distribution at $\mu = 0, 0.6, 0.8, 1.0$
          (\textit{top to bottom}).}
 \label{fig:riemann:s0s75s100s125}
\end{figure}
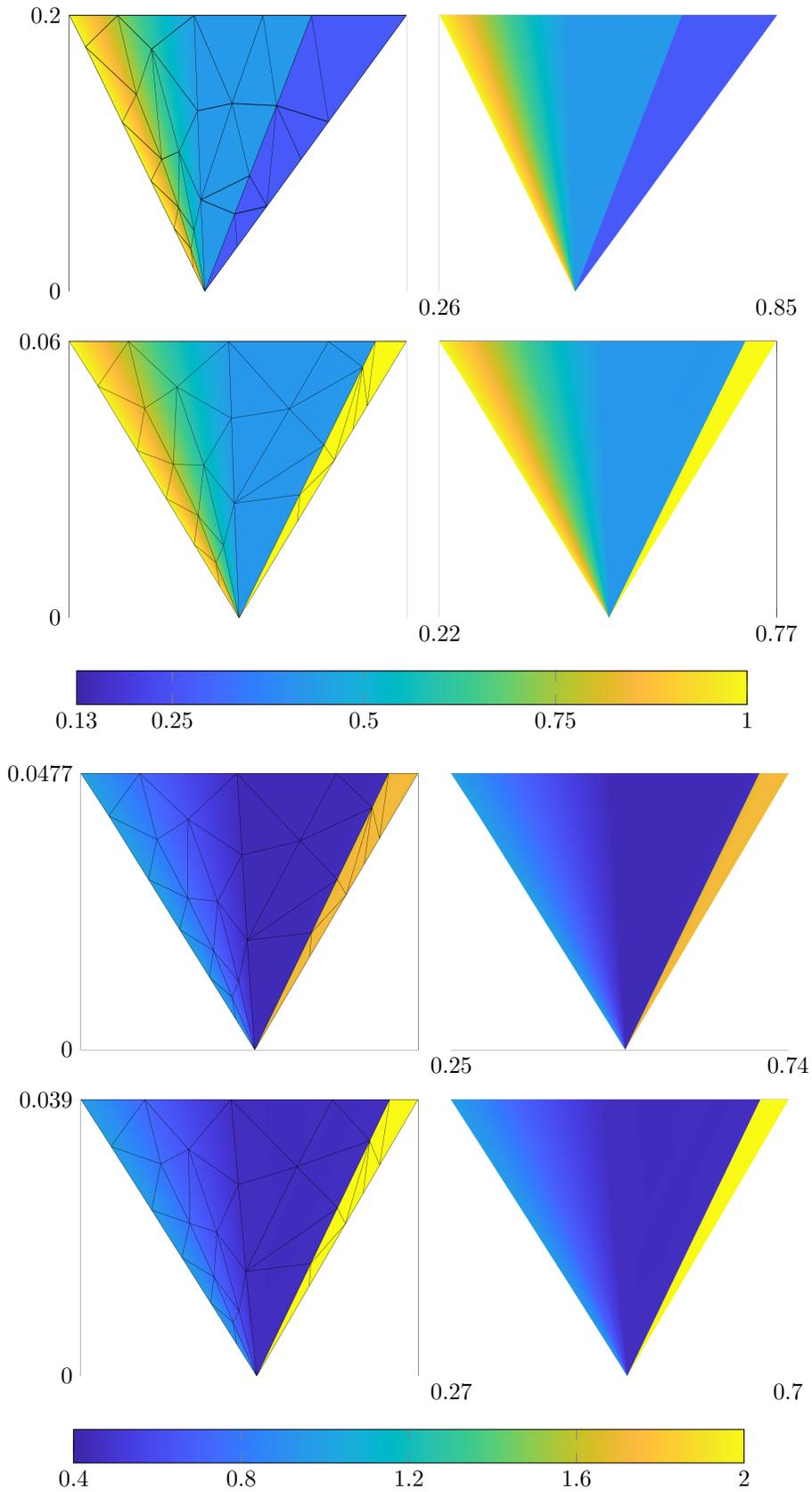

From this configuration, we perform a parameter sweep to reliably
compute the HOIST solution at $N = 125$ uniformly spaced samples
in $\Dcal$ (Figure~\ref{fig:riemann:s0s75s100s125}).
The HOIST solver parameters used for all $\Upsilon$ evaluations are
$\lambda=10^{-1}$, $\kappa_0=10^{-6}$, and $(\zeta, \upsilon)=(2, 0.5)$. 
To handle the parameter-dependent final time $T(\mu)$, we allow the
HOIST solver to determine the location of the top boundary ($t = T$).
Figure~\ref{fig:riemann:s0s75s100s125} shows four selected parameter
configurations. The upper temporal boundary adjusts to ensure the waves
travel roughly the same spatial distance. At $\mu = 0$, the speed of the
shock is approximately $1.75$ and the density jump magnitude is about $2$,
whereas at $\mu = 1$, the shock speed is $5.3$ and the density jump is $5.7$.
Finally, it can be seen that for larger values of $\mu$, the speed of the
shock and contact are close, which makes the constant wedge between the
shock and contact slender with a sharp angle. We have observed that thin
regions like this can be difficult to directly track from a shock-agnostic
mesh, often requiring a much finer mesh so elements can fit to the region
without being collapsed. The proposed many-query setting alleviates this
burden by initializing these difficult situations from nearly aligned grids
(from a HOIST solve at a different parameter configuration).

%
%

\begin{figure}
\centering
 \raisebox{-0.5\height}{\input{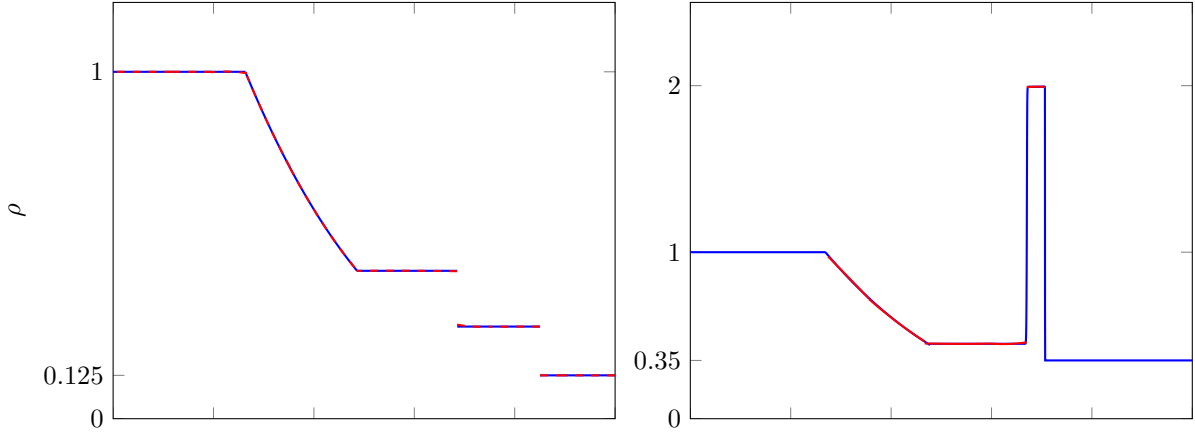}}
 \caption{Density slices at final time for $\mu = 0$ (\textit{left}) and $\mu = 1$
 (\textit{right}) for the HOIST solution (\ref{line:ist_sol_sod}) and reference solution
 (\ref{line:analytical_sol_sod}).}
 \label{fig:riemann_prob_slices}
\end{figure}

Finally, we verify the HOIST solution at the extremal parameters ($\mu = 0, 1$).
At $\mu = 0$ (Sod) we compare to the exact solution at the final time $T(\mu = 0)$
and at $\mu = 1$ we compare to a highly refined second-order finite volume simulation
at the final time $T(\mu = 1) = 0.039$. At both parameters, the HOIST solution matches
the reference solution well with the shock, contact, and head and tail of the
rarefaction tracked (Figure~\ref{fig:riemann_prob_slices}). This also shows the
dramatic variation the solution undergoes as the parameter varies throughout $\Dcal$.
 
\section{Conclusion}
\label{sec:conclude}
We introduce a specialized version of the high-order 
implicit shock tracking framework, originally developed in 
\cite{zahr_implicit_2020, huang2022robust}, for problems with
parametrized lead shocks. The approach applies implicit 
shock tracking (e.g., HOIST \cite{zahr_implicit_2020, huang2022robust}
or MDG-ICE \cite{corrigan_moving_2019}) on a shock-agnostic
mesh at one parameter configuration of interest to generate
a shock-aligned mesh and the corresponding flow field. All
elements upstream of the lead shock are removed to produce a
reduced mesh that will be used for all subsequent parameter
configurations, with the farfield boundary condition directly
applied on the shock boundary. In addition to the significant
reduction in elements (up to a factor of three in this work),
further reduction in the mesh DoFs is possible when the lead
shock is the only non-smooth feature in the domain because only
nodes on the shock boundary are optimized; the remainder are
determined from boundary constraints and PDE-based smoothing.
As a result, the shock boundary deforms as the parameters
change and the remaining nodes are positioned to improve element
quality. In addition to reducing the overall degrees
of freedom of the implicit shock tracking discretization, this also
improves robustness and accelerates convergence because the overall
shock-fitting problem is easier and a high-quality initial guess is
provided from the solution at previous parameter configurations.
The proposed framework can be used for most \textit{many-query}
applications involving parametrized lead shocks such as optimization,
uncertainty quantification, parameter sweeps, ``what-if'' scenarios,
or parameter-based continuation.

In this work, we use the abstract many-query setting for Mach number
continuation in steady inviscid flows and initial condition sweep
for one-dimensional Riemann problems. For continuation, we leverage
partially converged solves at intermediate stages to improve the efficiency
of the approach. A set of numerical experiments, which include two- and
three-dimensional supersonic and hypersonic flows, demonstrate the robustness
and flexibility of the proposed framework, in particular, the same HOIST solver
parameters can be used throughout the continuation process despite substantial
variations in Mach number (from $M_\infty=2$ to $M_\infty=10$) and initial state.
The framework was also shown to facilitate shock tracking for complex flow
features downstream of the lead shock, as demonstrated by the flow over the
double wedge where the shock interactions and supersonic jet were tracked
and well-resolved by the approach.

\section*{Acknowledgments}
This work is supported by AFOSR award numbers FA9550-20-1-0236, 
FA9550-22-1-0002, FA9550-22-1-0004, and ONR award number 
N00014-22-1-2299. The content of this publication does not necessarily 
reflect the position or policy of any of these supporters, and no official 
endorsement should be inferred.

\bibliographystyle{plain}
\bibliography{biblio,biblio_intro}

\end{document}